\documentclass{amsart}
\usepackage[latin1]{inputenc}
\usepackage[english,french]{babel}
\usepackage[T1]{fontenc}
\usepackage{graphicx}
\usepackage{amssymb}
\usepackage{amsthm}
\usepackage{amsmath}
\usepackage[all]{xy}
\usepackage{lipsum}

\begin{document}

\selectlanguage{english}  
\title{Newhouse phenomenon for automorphisms of low degree in $\mathbb{C}^{3}$}
\author{S\'ebastien Biebler}
\address{Sorbonne Universit\'e, 4 Place Jussieu, 75005 Paris}
\email{sebastien.biebler@imj-prg.fr}
\date{December 2019}
\thanks{This research was partially supported by the ANR project LAMBDA, ANR-13-BS01-0002.}

\begin{abstract}   
We show that there exists a polynomial automorphism $f$ of $\mathbb{C}^{3}$ of degree 2 such that for every automorphism $g$ sufficiently close to $f$, $g$ admits a tangency between the stable and unstable laminations of some hyperbolic set. As a consequence, for each $d \ge 2$, there exists an open set of polynomial automorphisms of degree at most $d$ in which the automorphisms having infinitely many sinks are dense. To prove these results, we give a complex analogous to the notion of blender introduced by Bonatti and D\'iaz.  
\end{abstract}

\maketitle

\small

\newtheorem* {mainTheorem}{Main Theorem}

\newtheorem {corollary}[subsubsection]{Corollary}
\newtheorem*{co1}{Corollary 1}
\newtheorem*{co2}{Corollary 2} 

\newtheorem {notation}[subsubsection]{Notation}
\newtheorem*{nt}{Notation}

\newtheorem{df}[subsubsection]{Definition}
\newtheorem{def2}[subsubsection]{Definition-Proposition}
\newtheorem*{opc}{Definition (Open Covering Property)}
\newtheorem*{ber}{Definition (Blender  Property)}

\newtheorem{prop}[subsubsection]{Proposition} 
\newtheorem{propo}[subsection]{Proposition}

\newtheorem{lemma}[subsubsection]{Lemma}

\newtheorem {remark} [subsubsection]{Remark}

\section{Introduction}

\subsection{Background}

Hyperbolic systems such as the horseshoe introduced by Smale were originally conjectured to be dense in the set of parameters in the 1960's. This was quickly discovered to be false in general for diffeomorphisms of manifolds of dimension greater than 2 (see \cite{po}). The discovery in the seventies of the so-called Newhouse phenomenon, i.e. the existence of residual sets of $C^{2}$-diffeomorphisms of compact surfaces with infinitely many sinks (periodic attractors) in \cite{n} showed it was false in dimension 2 too. The technical core of the proof is the reduction to a line of tangency between the stable and unstable foliations where two Cantor sets must have persistent intersections. This gives persistent homoclinic tangencies between the stable and unstable foliations, ultimately leading to infinitely many sinks. Indeed, it is a well known fact that a sink is created in the unfolding of a generic homoclinic tangency. \medskip

Palis and Viana showed in \cite{vi} an analogous result for real diffeomorphisms in higher dimensions. We say that a saddle periodic point of multipliers $| \lambda_{1}|  \le | \lambda_{2}|  < 1 <|  \lambda_{3}| $ is sectionally dissipative if the product of any two of its eigenvalues is less than 1 in modulus, that is, $| \lambda_{1}\lambda_{3}|  < 1$ and $| \lambda_{2}\lambda_{3}| < 1$ . More precisely, they proved that near any smooth diffeomorphism of $\mathbb{R}^{3}$ exhibiting a homoclinic tangency associated to a sectionally dissipative saddle periodic point, there is a residual subset of an open set of diffeomorphisms such that each of its elements displays infinitely many coexisting sinks.\medskip
 
In the complex setting, this reduction is not possible anymore and to get persistent homoclinic tangencies, we have to intersect two Cantor sets in the plane. Let us denote by $\text{Aut}_{d}(\mathbb{C}^{k})$ the space of polynomial automorphisms of $\mathbb{C}^{k}$ of degree $d$ for $d,k \ge 2$. Buzzard proved in \cite{bb1} that there exists an integer $d > 0$, an automorphism $G \in \text{Aut}_{d}(\mathbb{C}^{2})$ and a neighborhood $N \subset \mathrm{Aut}_{d}(\mathbb{C}^{2})$ of $G$ such that $N$ has persistent homoclinic tangencies. Buzzard gives an elegant criterion (see \cite{bb3}) which generates the intersection of two planar Cantor sets, hence leading to persistent homoclinic tangencies. In his article, Buzzard uses a Runge approximation argument to get a polynomial automorphism, which implies that the degree $d$ remains unknown and is supposedly very high.  \medskip

In the article \cite{bd1}, Bonatti and D\'iaz introduced a type of horseshoe they called blender horseshoe. The important property of such hyperbolic sets lies in the fractal configuration of one of their stable/unstable manifold which implies persistent intersections between any well oriented graph and this foliation. In some sense, the foliation behaves just as it had greater Hausdorff dimension than every individual manifold of the foliation. They find how to get robust homoclinic tangencies for some $C^{r}$-diffeomorphism of $\mathbb{R}^{3}$ using blenders in \cite{bdv1}. In the article \cite{dks}, one can find real polynomial maps of degree 2 with a blender. Other important studies of persistent tangencies using blenders include \cite{be} and \cite{rtem}.

\subsection{Results and outline} 

In this article, we generalize Buzzard's Theorem to dimension 3 and show that the degree can be controlled in this case. Here is our main result:
 
\begin{mainTheorem}
There exists a polynomial automorphism $f$ of degree 2 of $\mathbb{C}^{3}$ such that for every $g \in \mathrm{Aut}(\mathbb{C}^{3})$ sufficiently close to $f$, $g$ admits a tangency between the stable and unstable laminations of some hyperbolic set.  
\end{mainTheorem} 

Notice that in the previous result, $g$ is not assumed to be polynomial. 

\begin{co1}  
For each $d \ge 2$, there exists an open subset of $ \mathrm{Aut}_{d}(\mathbb{C}^{3})$ in which the automorphisms having a homoclinic tangency are dense. 
\end{co1} 

\begin{co2} 
For each $d \ge 2$, there exists an open subset $\mathrm{Aut}_{d}(\mathbb{C}^{3})$ in which the automorphisms having infinitely many sinks are dense. 
\end{co2}
   
Let us present the main ideas of the proof of this result. We consider the following automorphism of $\mathbb{C}^{3}$: 
\begin{equation} 
f_{0}  : (z_{1},z_{2},z_{3}) \mapsto (p_{c}(z_{1})+b  z_{2}+\sigma   z_{3}(z_{1}-\alpha),z_{1},\lambda z_{1} + \mu z_{3}+\nu) \label {e}
\end{equation}
where $p_{c}$ is a quadratic polynomial and the coefficients $b,\sigma,\alpha,\lambda,\mu,\nu$ are complex numbers. We prove that $f_{0}$ has a horseshoe $\mathcal{H}_{f_{0}}$ of index $(2,1)$: the first direction is strongly expanded, the second one is strongly contracted and the third one is moderately contracted by $f_{0}$. Informally speaking, the third projection restricted to $\mathcal{H}_{f_{0}}$ satisfies a special "open covering property" formalized in the following definition. This is an analogous in the complex setting of the notion of cs-blender in the sense of Bonatti and D\'iaz.  

\begin{ber} 
Let $f$ be a polynomial automorphism of $\mathbb{C}^{3}$, $D$ a tridisk of $\mathbb{C}^{3}$ and $\mathcal{H}_{f} = \bigcap_{- \infty}^{+ \infty}  f^{n}   (D)$ a horseshoe of index $(2,1)$. We will suppose that there exist $k>1$ and three cone fields $C^{u},C^{ss},C^{cs}$ such that in $D$: \begin{enumerate} \item $C^{u}$ is $f$-invariant, \item $C^{ss}$ is $f^{-1}$-invariant,  \item every vector in $C^{u}$ is expanded by a factor larger than $k$ under $f$, \item every vector in $C^{ss}$ is expanded by a factor larger than $k$ under $f^{-1}$, \item every vector in $C^{cs}$ is expanded by a factor larger than 1 under $f^{-1}$. \end{enumerate} We say that $\mathcal{H}_{f}$ is a blender if there exists a non empty open set $D' \Subset D$ such that every curve tangent to $C^{ss}$ intersecting $D'$ intersects the unstable set $W^{u}(\mathcal{H}_{f})$ of $\mathcal{H}_{f}$. 
\end{ber}

Besides, we show that $f_{0}$ has a periodic point which is sectionally dissipative. Once the blender is constructed, finding persistent tangencies is not trivial. We introduce manifolds with special geometry called folded manifolds. We prove that any folded manifold which is in good position has a tangency with the unstable manifold of a point of $\mathcal{H}_{f_{0}}$. We choose the parameters $c$, $b$ and $\sigma$ in order to create an initial heteroclinic tangency between the unstable manifold of a point of $\mathcal{H}_{f_{0}}$ and a folded manifold. This folded manifold is  in good position and included in the stable manifold of another point of $\mathcal{H}_{f_{0}}$. This enables us to produce persistent heteroclinic tangencies between stable and unstable manifolds of points of $\mathcal{H}_{f_{0}}$. This gives rise to homoclinic tangencies associated to the sectionally dissipative point. By a classical argument going back to Newhouse, this provides a subset of the set of automorphisms of degree 2 in which automorphisms displaying infinitely many sinks are dense. \medskip

An important point to notice is that the map $f_{0}$ defined in Eq. (1) is a perturbation of a skew product, with on the basis a H\'enon mapping (it is a skew product for $\sigma = 0$). The structure of H\'enon mapping will be important to create a horseshoe in Proposition \ref{ferm} and an initial fold in Proposition \ref{orr} (in particular, see Lemma \ref{utileresult}). The affine third coordinate is chosen so that the horseshoe displays the blender property (see Subsection 3.2). The perturbation term $\sigma   z_{3}(z_{1}-\alpha)$ allows to straighten the fold in a particular direction by iterating in Subsections 5.2 and 5.3. \medskip

The plan of the paper is as follows. In Section 2, we choose a family of quadratic polynomials and we fix complex coefficients $\lambda,\mu,\nu$. In Section 3, we introduce the map $f_{0}$ which depends on three parameters $c,b,\sigma$ and the associated horseshoe and we show that it has the blender property. Then, in Section 4, we introduce the formalism of folded manifolds and the mechanism which gives persistent tangencies. In Section 5, we prove that it is possible to choose $f_{0}$ in order to have a heteroclinic tangency. Finally, we prove the main Theorem in Section 6. In Appendix A, we explain how to construct a sink from a sectionally dissipative tangency. \newline
\newcommand{\lrpr}[1]{\left(#1\right)} \newline
\textbf{Note:} This article is a complete rewriting of a first version released on arXiv in November 2016. In that version the polynomial automorphism $f$ was of degree 5. To the best of the author's knowledge, the notion of blender was used there for the first time in holomorphic dynamics. Notice that blenders also appeared in complex dynamics in \cite{dm} and \cite{taflin}. 
\newline \newline \textbf{Acknowledgments :}
The author would like to thank his PhD advisor, Romain Dujardin as well as Pierre Berger and the anonymous referee for many invaluable comments.

\section{Preliminaries}
\subsection{Choice of a quadratic polynomial}

 In the following, we will consider the Euclidean norm on $\mathbb{C}^{n}$ for $n \in \{1,2,3\}$. 

\begin{notation}
 We denote by $\mathbb{D} \subset \mathbb{C}$ the open unit disk, and by $\mathbb{D}(0,r)$ the open disk centered at 0 of radius $r$ for any $r>0$. In particular, $\mathbb{D}(0,1) = \mathbb{D}$. \end{notation} 
\begin{notation} We will denote by $\mathrm{dist}$ the distance induced by the Euclidean norm on $\mathbb{C}^{n}$ for $n \in \{1,2,3\}$. 
\end{notation}

\begin{notation} 
For every $z \in \mathbb{C}^{3}$ and $i \in \{1,2,3\}$, we denote by $\mathrm{pr}_{i}(z) = z_{i}$ the $i^{th}$-coordinate of $z$.
\end{notation}

 In the following proposition, we carefully choose a family of quadratic polynomials with special properties.

\begin{prop} \label{marber}

For every integer $q >1$, there exists a disk $\mathcal{C} \subset \mathbb{C}$ of center $c_{0} \in \mathbb{C}$, a holomorphic family $(p_{c})_{c \in \mathcal{C}}$ of quadratic polynomials, two integers $m$ and $r$ (with $r$ independent of $q$), a constant $\chi>1$ and a disk $\mathbb{D}'$ with $\mathbb{D} \subset \mathbb{D}'$ such that: 

\begin{enumerate}
\item For every $c \in \mathcal{C}$, $p_{c}^{-r}( \mathbb{D})$ (resp. $p_{c}^{-r}( \mathbb{D}')$) admits two disjoint components $ \mathbb{D}_{1}, \mathbb{D}_{2}$ (resp. $ \mathbb{D}'_{1}, \mathbb{D}'_{2}$) included in $ \mathbb{D}$ (resp. $ \mathbb{D}$) such that $p_{c}^{r}$ is univalent on both $ \mathbb{D}_{1}$ and $ \mathbb{D}_{2}$ (resp. $ \mathbb{D}'_{1}$ and $ \mathbb{D}'_{2}$). Moreover $p_{c}^{r-1} (\mathbb{D}_{1}), p_{c}^{r-1}(\mathbb{D}_{2})  \Subset \mathbb{D}$ and $p_{c}^{r-1} (\mathbb{D}'_{1}), p_{c}^{r-1}(\mathbb{D}'_{2})  \Subset \mathbb{D}$.
\item  Denote by $\alpha_{c} = \bigcap_{n \ge 0} (p_{c}^{r})^{-n}(\mathbb{D}_{1})$ and $\gamma_{c} =  \bigcap_{n \ge 0} (p_{c}^{r})^{-n}(\mathbb{D}_{2})$ which are two fixed points of $p_{c}^{r}$. Then for every $c \in \mathcal{C}$, $\alpha_{c}$ is a repulsive fixed point of $p_{c}$, $|p_{c}'(\alpha_{c})|>\frac{6}{5}$ and we have:  
$$A : = r \alpha_{c_{0}} \neq \gamma_{c_{0}} +p_{c_{0}}(\gamma_{c_{0}})+ \cdots + p_{c_{0}}^{r-1}(\gamma_{c_{0}}) :  = B \text{ and } |A-B|>1 \,  .$$
\item We have $|p'_{c}|>\chi$ on $\mathbb{D}'$ and  $|(p^{r}_{c})'|>2$ on a neighborhood of $\mathbb{D}_{1} \cup \mathbb{D}_{2}$.
\item The critical point 0 is preperiodic for $c = c_{0}$ :  $p_{c_{0}}^{m}(0) = \alpha_{c_{0}} \neq 0$ with $p_{c_{0}}(0) \neq 0, \cdots, p_{c_{0}}^{m-1}(0) \neq 0$ and at $c= c_{0} $, we have:  $ \frac{d  }{dc}\big(  p_{c}^{m}(0) - \alpha_{c} \big) \neq 0$.
 \item There exists $R>0$ such that $\mathbb{D}' \subset \mathbb{D}(0,R)$ and such that the Julia set of $p_{c}$ is included in $\mathbb{D}(0,R)$ for every $c \in \mathcal{C}$. 
\item The polynomial $p_{c}$ has a periodic point $\delta_{c}$ of multiplier $\nu_{c}$ satisfying $1< |\nu_{c}|< (1+10^{-10})^{1/qr}$ for every $c \in \mathcal{C}$.
\end{enumerate}

\end{prop}

\begin{proof} 
We begin by working with the family of quadratic polynomials $p_{c}(z) = z^{2}+c$, we will rescale at the end of the proof. We begin by taking the only real quadratic polynomial $p_{a}(z) = z^{2}+a$ with one parabolic cycle $\delta_{a}$ of period 3. In particular, $a<-1$ and $a \in \mathbb{D}(0,2)$. For any $z \in \mathbb{C}$ such that $|z| \ge 10$, we have $|p_{a}(z)| = |z^{2}+a| \ge 10 |z|-|a| \ge 10 |z|-2$ and then $|p^{n}_{a}(z)| \rightarrow + \infty$. This shows that the Julia set of $p_{a}$ is strictly included in $\mathbb{D}(0,10)$. Simple calculations show that $z^{2}+a$ has two real fixed points $\alpha^{+}_{a} = \frac{1}{2}(1+\sqrt{1-4a})>\frac{1}{2}(1+\sqrt{5})>1$ and $\alpha^{-}_{a} = \frac{1}{2}(1-\sqrt{1-4a})<\frac{1}{2}(1-\sqrt{5})<-\frac{6}{10}$. We take two open disks $ \mathbb{B}'_{+} \subset \mathbb{D}(0,10)$ and $ \mathbb{B}'_{-} \subset \mathbb{D}(0,10)$ respectively centered around $\alpha_{a}^{+} $ and $\alpha^{-}_{a}$ which are both disjoint from the orbit of the critical point 0 of $z^{2}+a$ (this is possible since the critical orbit tends to the parabolic orbit of $z^{2}+a$). Since $\alpha_{a}^{+} $ and $\alpha^{-}_{a}$ are repulsive fixed points of $p_{a}$, there exists some $\chi>1$ such that  $|p'_{c}|>\chi$ on $ \mathbb{B}'_{+} \cup \mathbb{B}'_{-} $, up to reducing $ \mathbb{B}'_{+}$ and $\mathbb{B}'_{-} $ if necessary. \medskip

Since $\alpha_{a}^{+} $ and $\alpha^{-}_{a}$ are repulsive fixed points, still reducing $ \mathbb{B}'_{+}$ and $ \mathbb{B}'_{-}$ if necessary, we have that for every $r \ge 1$, there is a connected component of $p_{a}^{-r}(\mathbb{B}'_{+})$ (resp. $p_{a}^{-r}(\mathbb{B}'_{-})$) which contains $\alpha_{a}^{+} $ (resp.  $\alpha^{-}_{a}$) and whose $r$ first iterates are all included in $\mathbb{B}'_{+}$ (resp. $\mathbb{B}'_{-})$. We denote by $\tilde{ \mathbb{B}}_{+}$ and $\tilde{ \mathbb{B}}_{-}$  the respective connected components of $p_{a}^{-1}(\mathbb{B}'_{+})$ and $p_{a}^{-1}(\mathbb{B}'_{-})$ which contain $\alpha_{a}^{+} $ and $\alpha^{-}_{a}$ and are defined this way. Then we fix open disks $\mathbb{B}_{+}$ and $\mathbb{B}_{-}$ of respective centers $\alpha_{a}^{+}$ and  $\alpha^{-}_{a}$  such that $\tilde{ \mathbb{B}}_{+} \Subset  \mathbb{B}_{+} \Subset \mathbb{B}'_{+}$ and $\tilde{ \mathbb{B}}_{-} \Subset  \mathbb{B}_{-} \Subset \mathbb{B}'_{-}$. Since both $\mathbb{B}'_{+}$ and $\mathbb{B}'_{-}$ intersect the Julia set of $p_{a}$ and are disjoint from the critical orbit, we can find some integer $\overline{r}$ such that $p^{\overline{r}}_{a}(\mathbb{B}'_{+})$ contains $\mathbb{B}'_{-}$ and $p^{\overline{r}}_{a}(\mathbb{B}'_{-})$ contains $\mathbb{B}'_{+}$. Then we can find some open set $\overline{B}_{+} \Subset \tilde{ \mathbb{B}}_{+} \Subset \mathbb{B}_{+}$ satisfying $p_{a}(\overline{B}_{+}) \Subset \mathbb{B}'_{+}$, $p_{a}^{1+\overline{r}}(\overline{B}_{+}) \Subset  \tilde{ \mathbb{B}}_{-} $, $p_{a}^{2+\overline{r}}(\overline{B}_{+}) \Subset  \mathbb{B}'_{-} $ and $p_{a}^{2+2\overline{r}}(\overline{B}_{+}) = \mathbb{B}'_{+}$. Hence, denoting $r = 2+2\overline{r}$, we have $\overline{B}_{+} \Subset \mathbb{B}_{+}$ and $p_{a}^{r}$ sends $\overline{B}_{+}$ biholomorphically onto $\mathbb{B}'_{+}$. We denote by $\gamma_{a}$ the periodic point of $p_{a}$ of period $r$  which is the unique fixed point of the restriction of $p_{a}^{r}$ to $\overline{B}_{+}$. We notice that $\gamma_{a}  \neq \alpha^{+}_{a}$. Similarly, we can define $\overline{B}_{-} \Subset \mathbb{B}_{-}$ such that $p_{a}^{r}$ sends $\overline{B}_{-}$ biholomorphically onto $\mathbb{B}'_{-}$ and $p_{a}^{r/2}(\gamma_{a}) = p^{1+\overline{r}}_{a}(\gamma_{a})     \neq \alpha^{-}_{a}$ is the unique fixed point of the restriction of $p_{a}^{r}$ to $\overline{B}_{-}$. \medskip

Since $\alpha_{a}^{+} \neq \alpha^{-}_{a}$, it is not possible to satisfy simultaneously $\gamma_{a} +p_{a}(\gamma_{a})+ \cdots + p_{a}^{r-1}(\gamma_{a})  =  r \alpha_{a}^{+}$ and $ \gamma_{a} +p_{a}(\gamma_{a})+ \cdots + p_{a}^{r-1}(\gamma_{a})  = r \alpha_{a}^{-}$. In the following, we will denote by $\alpha_{a}$ a point in $\{\alpha_{a}^{+}, \alpha^{-}_{a}\}$ such that the inequality $  \gamma_{a} +p_{a}(\gamma_{a})+ \cdots + p_{a}^{r-1}(\gamma_{a})  \neq r \alpha_{a}$ is satisfied. We also denote by $\mathbb{B}$, $\mathbb{B}'$, $\tilde{\mathbb{B}}$ and $\overline{B}$ the sets corresponding to $\alpha_{a}$. Up to replacing $\gamma_{a}$ by $p_{a}^{r/2}(\gamma)$ if $\alpha_{a} = \alpha^{-}_{a}$, we can suppose that $\gamma_{a} \in \mathbb{B}$. The multiplier of $\alpha_{a}$ is of  modulus $|2\alpha_{a}|>\mathrm{min}(2,\frac{6}{5}) = \frac{6}{5}$. We take the component $ \mathbb{B}'_{1}$  of $p^{-r}_{a}(\mathbb{B}')$ containing $\alpha_{a}$ and where $p^{r}_{a} $ is univalent defined at the beginning of the last paragraph. We have $ \mathbb{B}'_{1} \Subset \mathbb{B} \Subset  \mathbb{B}'$. We also take the component $ \mathbb{B}'_{2}$  of $p^{-r}_{a} ( \mathbb{B}')$  containing $\gamma_{a}$ and where $p^{r}_{a} $ is univalent equal to $\overline{B}$. It holds $ \mathbb{B}'_{2}  \Subset \tilde{\mathbb{B}}\Subset \mathbb{B} \Subset  \mathbb{B}'$. Replacing $r$ by one of its multiples if necessary (still denoted by $r$), $ \mathbb{B}'_{1}$ and $ \mathbb{B}'_{2}$ are disjoint. We also take the respective components $ \mathbb{B}_{1}$ and $ \mathbb{B}_{2}$ of $p^{-r}_{a} ( \mathbb{B})$ included into those of $p^{-r}_{a} ( \mathbb{B}') $. Since $\tilde{ \mathbb{B}} \Subset  \mathbb{B}$, it holds $p_{a}^{r-1} (\mathbb{B}_{1}), p_{a}^{r-1}(\mathbb{B}_{2}) , p_{a}^{r-1} (\mathbb{B}'_{1}), p_{a}^{r-1}(\mathbb{B}'_{2}) \Subset \mathbb{B}$. Still replacing $r$ by a  multiple if necessary, we have  $ |r \alpha_{a} - (\gamma_{a} +p_{a}(\gamma_{a})+ \cdots + p_{a}^{r-1}(\gamma_{a}))|>10 $. Since $ \mathbb{B}'_{1}  \Subset \mathbb{B}$ and $ \mathbb{B}'_{2}  \Subset  \mathbb{B}$, by the Schwarz Lemma, there exists $\theta>1$ such that $|(p_{c}^{r})'|>\theta$ on a neighborhood of $\mathbb{B}_{1} \cup  \mathbb{B}_{2}$. Taking a  multiple of $r$ if necessary, $|(p_{c}^{r})'|>2$ on a neighborhood of $\mathbb{B}_{1} \cup  \mathbb{B}_{2}$.\medskip
 
Let us fix $q>1$. By continuity, for $c $ in some neighborhood $\mathcal{C}_{a}$ of $a$ in $\mathbb{C}$, it holds: 

\begin{enumerate}
\item $p_{c}^{-r}( \mathbb{B})$ (resp. $p_{c}^{-r}( \mathbb{B}')$) admits two components $ \mathbb{B}_{1}, \mathbb{B}_{2}$ (resp. $ \mathbb{B}'_{1}, \mathbb{B}'_{2}$) included in $ \mathbb{B}$ (resp. $ \mathbb{B}$) containing the continuations $\alpha_{c}$ and $\gamma_{c}$ and  such that $p_{c}^{r}$ is univalent on both $ \mathbb{B}_{1}$ and $ \mathbb{B}_{2}$ (resp. $ \mathbb{B}'_{1}$ and $ \mathbb{B}'_{2}$). Moreover  $p_{c}^{r-1} (\mathbb{B}_{1})$, $p_{c}^{r-1}(\mathbb{B}_{2}),$  $p_{c}^{r-1} (\mathbb{B}'_{1})$, $p_{c}^{r-1}(\mathbb{B}'_{2})$ $ \Subset \mathbb{B}$,
 \item the continuation $\alpha_{c}$ of $\alpha_{a}$ is a repulsive fixed point of $p_{c}$ such that $|p'_{c}(\alpha_{c})|>\frac{6}{5}$,
\item $ r \alpha_{c} \neq \gamma_{c} +p_{c}(\gamma_{c})+ \cdots + p_{c}^{r-1}(\gamma_{c}) $ and $|r \alpha_{c} -(  \gamma_{c} +p_{c}(\gamma_{c})+ \cdots + p_{c}^{r-1}(\gamma_{c}))|>10$,
\item $|p'_{c}|>\chi$ on $\mathbb{B}'$ and $|(p_{c}^{r})'|>2$ on a neighborhood of $\mathbb{B}_{1} \cup  \mathbb{B}_{2}$,
\item  the Julia set of $p_{c}$ is included in $\mathbb{D}(0,10)$, 
\item  the continuation $\delta_{c}$ of $\delta_{a}$ is of multiplier $\nu_{c}$ such that $(1-10^{-10})^{1/qr} < |\nu_{c}| < (1+10^{-10})^{1/qr}$.
\end{enumerate} 

The parameter $a$ belongs to the Mandelbrot set. Misiurewicz parameters are dense inside the Mandelbrot set so it is possible to find a parameter $\tilde{c}$ inside the interior of $\mathcal{C}_{a}$ such that the critical point 0 is preperiodic for $p_{\tilde{c}}$. The critical point 0 is sent after a finite number of iterations of $p_{\tilde{c}}$ on a periodic orbit. This periodic orbit is accumulated by preimages of $\alpha_{\tilde{c}}$ by iterates of $p_{\tilde{c}}$. Then by the Argument Principle it is possible to take a new Misiurewicz parameter $c_{0} $ in the interior of $\mathcal{C}_{a}$ such that 0 is still preperiodic but with associated orbit the fixed point $\alpha_{c_{0}}$. There exists an integer $m$ such that $p_{c_{0}}^{m}(0) = \alpha_{c_{0}}$  with $p_{c_{0}}(0) \neq 0, \cdots, p_{c_{0}}^{m-1}(0) \neq 0$. The inequality $ \frac{d  }{dc}  \big( p_{c}^{m}(0)-\alpha_{c}  \big) \neq 0$ at $c = c_{0}$ is a direct consequence of Lemma 1, Chapter 5 of \cite{hudo}. For the parameter $c_{0}$, $\delta_{c_{0}}$ is repulsive of multiplier $\nu_{c_{0}}$ such that $1 < | \nu_{c_{0}}| < (1+10^{-10})^{1/qr}$. We pick some ball $\mathcal{C} \subset \mathcal{C}_{a}$ of center $c_{0}$ where this is still true. \medskip

For each $c \in \mathcal{C}$, we do a rescaling by an affine map so that after rescaling $\mathbb{B} \subset \mathbb{D}(0,10)$ is sent on $\mathbb{D} = \mathbb{D}(0,1)  $. Properties 1, 2, 4 and 6 are still true. Property 5 is still true with a disk $\mathbb{D}(0,R)$ with a fixed $R>0$ instead of $\mathbb{D}(0,10)$. Since $r\alpha_{c} \neq \gamma_{c} +p_{c}(\gamma_{c})+ \cdots + p_{c}^{r-1}(\gamma_{c})$ and $|r \alpha_{c} -(  \gamma_{c} +p_{c}(\gamma_{c})+ \cdots + p_{c}^{r-1}(\gamma_{c}))|>10$ before rescaling, we have $A \neq B$ and $|A-B|>1$ after and then Property 3 is true.  Then Properties 1, 2, 3, 4, 5 and 6 are satisfied for every $c \in \mathcal{C}$. In the following, after rescaling, we will denote $\mathbb{B}, \mathbb{B}', \mathbb{B}_{1}, \mathbb{B}_{2}, \mathbb{B}'_{1}, \mathbb{B}'_{2}$ by $\mathbb{D}, \mathbb{D}', \mathbb{D}_{1}, \mathbb{D}_{2}, \mathbb{D}'_{1}, \mathbb{D}'_{2}$. For simplicity, we will still denote by $p_{c}$ the polynomial after rescaling.
\end{proof}

\subsection{Choice of an IFS}

\begin{notation} \label{defh1h}
For every $c \in \mathcal{C}$, we denote by $h_{1}$ and $h_{2}$ the two inverse branches of $p_{c}^{r}$ on
 $\mathbb{D}'$ given by Proposition \ref{marber} such that $\alpha_{c} = \bigcap_{n \ge 0} h_{1}^{n}( \mathbb{D})$ and $\gamma_{c} = \bigcap_{n \ge 0} h^{n}_{2}( \mathbb{D})$.    
\end{notation}

\begin{notation} \label{h}
We denote $\mu_{0} = (1-10^{-4})^{\frac{1}{qr}} \cdot e^{i \cdot \frac{\pi}{2qr}}$ which depends on the integer $q$. In particular, we have the following equality:  $\mu_{0}^{qr} = (1-10^{-4}) \cdot e^{i \cdot \frac{\pi}{2}}$.
\end{notation} 
 
In the following result, we iterate $q$ times the maps $h_{1}$ and $h_{2}$ with a specific choice for the integer $q$. Remind that $A$, $B$ and $R$ were defined in Proposition \ref{marber}.

\begin{prop} \label{fo}

There exists an integer $q \ge 100$ such that, after reducing $\mathcal{C}$ if necessary, the following holds for every $c \in \mathcal{C}$:

\begin{enumerate}
\item $|(h_{j}^{q})'| <10^{-10}$ for $j \in \{1,2\}$ on a neighborhood $ \mathbb{D}''$ of $ \mathbb{D}$ with $ \mathbb{D} \subset  \mathbb{D}'' \subset  \mathbb{D}'$,
\item $\mathrm{diam} \big( h^{q}_{j} ( \mathbb{D}' ) \big) \le  10^{-11}  \cdot \mathrm{dist} (  h^{q}_{j} ( \mathbb{D}' )  ,  \partial \mathbb{D} )$ for $j \in \{1,2\}$, 
\item  for every $z \in h_{1}^{q}( \mathbb{D}')$ and $0 \le n \le qr  \big( 1-10^{-10} r^{-1}  R^{-1} \min(1, |A-B|) \big)  $:  $$  \mu_{0}^{0}   p_{c}^{n+r-1}(z)+ \cdots + \mu^{r-1}_{0}  p_{c}^{n}(z)  \in \mathbb{D}(A, 10^{-10} \cdot  |A-B| )  \, , $$ 
\item for every $z \in h^{q}_{2}( \mathbb{D}')$ and $0 \le n \le qr \big( 1-10^{-10}  r^{-1} R^{-1} \min(1,|A-B|) \big)$: $$  \mu_{0}^{0} p_{c}^{n+r-1}(z)+ \cdots + \mu^{r-1}_{0}  p_{c}^{n}(z)  \in \mathbb{D}(B, 10^{-10} \cdot |A-B|  )  \, . $$
\end{enumerate}

\end{prop}

\begin{proof} 
We first show the result for $c = c_{0}$. According to property (3) of Proposition \ref{marber},  $|(p^{r}_{c})'|>2$ on a neighborhood of $\mathbb{D}_{1} \cup \mathbb{D}_{2}$. Then, taking $q \ge 100$ such that $2^{q}> 10^{10}$, we have $|(h_{j}^{q})'|<10^{-10}$ on some disk $ \mathbb{D}''$ with $ \mathbb{D} \subset  \mathbb{D}'' \subset  \mathbb{D}'$. Since $h_{j}$ is a contraction such that $\bigcap_{n \ge 0} h_{1}^{n}(\mathbb{D}') = \{\alpha_{c}\}$ and $\bigcap_{n \ge 0} h^{n}_{2}(\mathbb{D}') = \{\gamma_{c}\}$, increasing the value of $q$ if necessary, we have that $\mathrm{diam}\big( (h^{q}_{j} ( \mathbb{D}' )\big) \le  10^{-11}  \cdot \text{dist} (  h^{q}_{j} ( \mathbb{D}' )  ,  \partial \mathbb{D} )$. When $q \rightarrow + \infty$, we both have $\mu_{0}^{k} \rightarrow 1$ and $p_{c_{0}}^{n+k}(z) \rightarrow \alpha_{c_{0}}$ uniformly in $0 \le k < r$, $0 \le n \le qr(1-10^{-10}  r^{-1} R^{-1} \min(1, |A-B|)) $ and $z \in h_{1}^{q}( \mathbb{D}')$. Then, increasing the value of $q$ if necessary, we have that $ \mu_{0}^{0}  p_{c_{0}}^{n+r-1}(z)+ \cdots + \mu^{r-1}_{0}  p_{c_{0}}^{n}(z)  \in \mathbb{D}(A, 10^{-10} |A-B|)$. The proof of the last item is similar. Since all these conditions are open, reducing the ball $\mathcal{C}$ of center $c_{0}$ if necessary, they  remain true for every $c \in \mathcal{C}$.
\end{proof} 

\begin{notation} \label{rbound}
Since $r$ is independent of $q$ (see Proposition \ref{marber}), we can increase $q$ so that $r \le 10^{-10} qR^{-1}  \min(1, |A-B|)$. From now on, we fix such a value of $q$ and the associated value $\mu_{0}$. \end{notation}

\subsection{Choice of the parameters $\lambda$ and $\nu$}

In this Subsection, we introduce two new coefficients $\lambda$ and $\nu$. These constants will apppear on the third coordinate of the polynomial automorphisms of $\mathbb{C}^{3}$ we are going to work with. This will be used to create a horseshoe in Proposition \ref{ferm} and to show that this horseshoe displays the blender property in Subsection 3.2.
 
\begin{notation}
We denote by $A' = (\mu_{0}^{r-1}   \alpha_{c_{0}}  +\cdots +  \mu_{0}^{0} p_{c_{0}}^{r-1}(\alpha_{c_{0}}) )$ and $B' = (  \mu_{0}^{r-1}  \gamma_{c_{0}}+\cdots +  \mu_{0}^{0} p_{c_{0}}^{r-1}(\gamma_{c_{0}})   )$. \end{notation}

By Proposition \ref{fo}, $A' \in \mathbb{D}(A,10^{-10}|A-B|)$ and $B' \in \mathbb{D}(B,10^{-10}|A-B|)$. According to item 2 of Proposition \ref{marber}, this implies that:
\begin{equation} \label{ineqb}
|A'-B'|> \frac{1}{2}|A-B|>\frac{1}{2} \, . 
\end{equation}

\begin{prop}\label{munu}
There exist two constants $\lambda,\nu$ such that $|\lambda|<1$ and satisfying: $$\lambda A'(1+ \mu_{0}^{r}+ \cdots + \mu_{0}^{qr-r}) + \nu  (1+\mu_{0}+ \cdots +\mu^{qr-1}_{0}) =\frac{9}{10} \cdot 10^{-4} \, ,$$ 
$$\lambda B' (1+ \mu_{0}^{r}+ \cdots + \mu_{0}^{qr-r}) + \nu  (1+\mu_{0}+ \cdots +\mu^{qr-1}_{0}) = - \frac{9}{10} \cdot 10^{-4} \, .$$
\end{prop}

\begin{proof} 
We have: $1+ \mu_{0}^{r}+ \cdots + \mu_{0}^{qr-r} = (1-   \mu_{0}^{qr})/(1-   \mu^{r}_{0} ) $. By Notation \ref{h}, we have $\mu_{0}^{qr} = (1-10^{-4}) \cdot e^{i \cdot \frac{\pi}{2}} \neq 1$ and then $1+ \mu_{0}^{r}+ \cdots + \mu_{0}^{qr-r} \neq 0$. Similarly we have $1+ \mu_{0}+ \cdots + \mu_{0}^{qr-1} = (1-   \mu_{0}^{qr})/(1-   \mu_{0} ) \neq 0$. Since $A' \neq B'$, it is possible to pick two coefficients $\lambda$ and $\nu$ so that the images of these two complex numbers by the affine map $z \mapsto \lambda  (1+ \mu_{0}^{r}+ \cdots + \mu_{0}^{qr-r})  z + \nu (1+ \mu_{0}+ \cdots + \mu_{0}^{qr-1}) $ are respectively equal to $\frac{9}{10} \cdot 10^{-4}$ and $-\frac{9}{10} \cdot 10^{-4}$. It remains to show that $|\lambda|<1$. To this end, we will need the following technical lemma: 

\begin{lemma} \label{mirber}
The complex number $\mu_{0}$ satisfies the following inequality: 
$$ \frac{q}{2} \le 1+|\mu_{0}^{r}|+|\mu_{0}|^{2r}+ \cdots + |\mu_{0}|^{qr-r} \le 10 \cdot | 1+\mu_{0}^{r}+\mu_{0}^{2r}+ \cdots +\mu_{0}^{qr-r}| \, .$$
\end{lemma}

\begin{proof}
We have $\frac{1}{2} \le 1$, $\frac{1}{2} \le |\mu_{0}^{r}|$, $\cdots$, $\frac{1}{2} \le |\mu_{0}^{qr-r}|$ so the first inequality is trivial. Since every term $\mu_{0}^{nr}$ ($0 \le n < q$) has a positive real part and since this real part is larger than $\frac{1}{2}$ for $0 \le n \le \frac{1}{2}(q-1)$, we have $\frac{1}{2}(q-1) \cdot \frac{1}{2} \le \text{Re}(  1+\mu_{0}^{r}+\mu_{0}^{2r}+ \cdots +\mu_{0}^{qr-r}    )$ and then $1+|\mu_{0}^{r}|+|\mu_{0}|^{2r}+ \cdots + |\mu_{0}|^{qr-r} \le q \le 10 \cdot  \frac{1}{2}(q-1) \cdot \frac{1}{2} \le 10 \cdot | 1+\mu_{0}^{r}+\mu_{0}^{2r}+ \cdots +\mu_{0}^{qr-r}|$. The proof is complete.
\end{proof}

We are now in position to end the proof of Lemma \ref{munu}. By definition of $\lambda$ and $\nu$, we  have $|\lambda| |A'-B'| |1+ \mu_{0}^{r}+ \cdots + \mu_{0}^{qr-r}| = 2 \cdot \frac{9}{10} \cdot 10^{-4}$. We already proved that $|A'-B'|>\frac{1}{2}$ in Eq. (\ref{ineqb}) and by Lemma \ref{mirber} we also have $|1+ \mu_{0}^{r}+ \cdots + \mu_{0}^{qr-r}| \ge q/20 \ge 100/20  \ge 1 $. This implies that $|\lambda|<1$ and so the result is proven.
 \end{proof} 

\begin{corollary} \label{hyu}
Reducing $\mathcal{C}$ if necessary, there exists a neighborhood $\mathcal{B}_{\mu}$ of $\mu_{0}$ such that for every $c \in \mathcal{C}$ and $\mu \in \mathcal{B}_{\mu}$ it holds:
 
\begin{enumerate} 
\item for every $z \in h_{1}^{q}( \mathbb{D}')$, we have: $$ \nu + \lambda p_{c}^{qr-1}(z)+ \mu  (\nu+ \cdots +\mu  (\nu+  \lambda z )) \in \mathbb{D}(\frac{9}{10} \cdot 10^{-4}, 10^{-10}) \, , $$ 
\item  for every $z \in h_{2}^{q}( \mathbb{D}')$, we have: $$ \nu + \lambda  p_{c}^{q-1}(z)+ \mu  (\nu+ \cdots +\mu (\nu+  \lambda  z )) \in  \mathbb{D}(- \frac{9}{10} \cdot 10^{-4},10^{-10}) \, . $$ \end{enumerate}

\end{corollary}
 
\begin{proof} 
We first prove the result for $c = c_{0}$ and $\mu = \mu_{0}$. According to Proposition \ref{munu}, we have: 
$$ \nu + \lambda   p^{qr-1}_{c_{0}}(z)+ \mu_{0} (\nu+ \cdots +\mu_{0}  (\nu+  \lambda  z )) - \frac{9}{10} \cdot 10^{-4} =$$ $$ \lambda  \sum_{n= 0}^{l-1} (p_{c_{0}}^{qr-1-n}(z)-\alpha_{c_{0}})\mu_{0}^{n} + \lambda  \sum_{n= l}^{qr-1} (p_{c_{0}}^{qr-1-n}(z)-\alpha_{c_{0}})\mu_{0}^{n} \, ,$$ 
where $l$ is the smallest integer such that $l \ge 10^{-10} qR^{-1}  \min(1, |A-B|)$ and which is a multiple of $r$. By Notation \ref{rbound}, we have $l \le 2 \cdot 10^{-10} qR^{-1}  \min(1, |A-B|)$. In particular, $qr-l$ is a multiple of $r$. Using the third item of Proposition \ref{fo}, it holds: 
$$  |\lambda  \sum_{n= l}^{qr-1} (p^{qr-1-n}_{c_{0}}(z)-\alpha_{c_{0}})\mu_{0}^{n}| \le  |\lambda|   \cdot 10^{-10}   |A- B|  \cdot (1+|\mu_{0}|^{r}+|\mu_{0}|^{2r}+ \cdots + |\mu_{0}|^{qr-r}) \, .$$ 
We already proved that $|A'-B'| > \frac{1}{2} |A-B|$ in Eq. (\ref{ineqb}). In particular, this implies that $|\lambda| \cdot |A-B| \cdot | 1+\mu^{r}_{0}+\mu_{0}^{2r}+ \cdots +\mu_{0}^{qr-r}| < |\lambda| \cdot 2 |A'-B'| \cdot |1+ \mu_{0}^{r}+ \cdots + \mu_{0}^{qr-r}| $. Then, by Proposition \ref{munu}, this yields $|\lambda| \cdot |A-B| \cdot | 1+\mu^{r}_{0}+\mu_{0}^{2r}+ \cdots +\mu_{0}^{qr-r}| < 2 \cdot  2 \cdot \frac{9}{10} \cdot 10^{-4}$. By Lemma \ref{mirber}, it also holds: $1+|\mu_{0}^{r}|+|\mu_{0}|^{2r}+ \cdots + |\mu_{0}|^{qr-r} \le 10 \cdot | 1+\mu_{0}^{r}+\mu_{0}^{2r}+ \cdots +\mu_{0}^{qr-r}|$. All this together implies the following: 

\begin{equation}
\label{ineqmu1}
 |\lambda  \sum_{n= l}^{qr-1} (p^{qr-1-n}_{c_{0}}(z)-\alpha_{c_{0}})\mu_{0}^{n}| \le 10^{-10} \cdot 10 \cdot \big(  2 \cdot  2 \cdot \frac{9}{10} \cdot 10^{-4} \big) \le 10^{-11} \, . 
\end{equation}

Since both $\mathbb{D}'$ and the Julia set of $p_{c_{0}}$ are included in $\mathbb{D}(0,R)$ (see item 5 of Proposition \ref{marber}), we also have: 
$$  | \lambda  \sum_{n= 0}^{l-1} (p^{qr-1-n}_{c_{0}}(z)-\alpha_{c_{0}})\mu_{0}^{n}| \le   |\lambda  | \cdot 2R  \cdot (1+|\mu_{0}|+|\mu_{0}|^{2}+ \cdots + |\mu_{0}|^{l-1})  \le   |\lambda  | \cdot  2R \cdot  l   \, . $$ 
Since $l \le 2 \cdot 10^{-10} qR^{-1}  \min(1, |A-B|)$ and then using the inequality $|A-B| < 2|A'-B'|$ from Eq. (\ref{ineqb}), the latter is smaller than:
$$ |\lambda| \cdot 2R \cdot 2 \cdot 10^{-10}\frac{q}{R} \cdot \min(1, |A-B|) \le 10^{-8} \cdot |\lambda| \cdot \frac{q}{2} \cdot \frac{|A-B|}{2} \le 10^{-8} \cdot |\lambda| \cdot \frac{q}{2} \cdot |A'-B'| \, .$$  
Using successively the inequality $ \frac{q}{2} \le 10 \cdot | 1+\mu_{0}^{r}+\mu_{0}^{2r}+ \cdots +\mu_{0}^{qr-r}| $ from Lemma \ref{mirber} and then Proposition \ref{munu}, the latter is finally smaller than $10^{-8}  \cdot |\lambda| \cdot |A'-B'| \cdot 10 \cdot | 1+\mu_{0}+\mu_{0}^{2}+ \cdots +\mu_{0}^{qr-r}| \le 10^{-7} \cdot 2 \cdot \frac{9}{10} \cdot 10^{-4}$ and so:

\begin{equation}
\label{ineqmu3}
| \lambda  \sum_{n= 0}^{l-1} (p^{qr-1-n}_{c_{0}}(z)-\alpha_{c_{0}})\mu_{0}^{n}| \le 2 \cdot 10^{-11} \, .
\end{equation}

Then we just have to sum the two inequalities of Eq. (\ref{ineqmu1}) and Eq. (\ref{ineqmu3}) to prove the result for $c= c_{0}$ and $\mu = \mu_{0}$. By continuity and since the inequality is open, it remains true for every $\mu$ in some ball $ \mathcal{B}_{\mu}$ of center $\mu_{0}$ and $c \in \mathcal{C}$ after reducing $\mathcal{C}$ if necessary. Then item 1 is true and the proof of item 2 is similar. The result is proven.\end{proof}

\begin{remark} \label{rmu} 
We reduce $\mathcal{B}_{\mu}$ so that we both have $|\mu| < (1-  10^{-4} +10^{-10}     )^{\frac{1}{qr}}$, $|\mu|^{2qr}>1-2 \cdot 10^{-4}$ and $\mu^{qr} \subset \mathbb{D}(\mu_{0}^{qr},10^{-10})$ for every $\mu \in \mathcal{B}_{\mu}$.  
\end{remark}

\subsection{Adjusting the parameter $\mu$}

In this subsection, we slightly perturb the coefficient $\mu_{0}$ into a new value $\mu$ in order to satisfy some equality for a product of matrices. Notice that this choice has nothing to do with the next section and the blender property, it will be useful in Section 5.

\begin{notation} \label{bet}
We denote $\beta_{0} = 0, \beta_{1} = p_{c_{0}}(0), \beta_{2} = p^{2}_{c_{0}}(0), \ldots, \beta_{m} = p^{m}_{c_{0}}(0)= \alpha_{c_{0}}$ the points of the orbit of the critical point 0 before landing onto the fixed point $\alpha_{c_{0}}$.
\end{notation}   
 
\begin{notation}
For every $\mu \in \mathcal{B}_{\mu}$, we define: $ w_{0} =0 + 1 \cdot (\beta_{0}-\beta_{m})$,  $w_{1} = p'_{c_{0}}(\beta_{1})  w_{0} + \mu (\beta_{1} -\beta_{m})$ $ \cdots $ and $w_{m-1} = p'_{c_{0}}(\beta_{m-1})  w_{m-2} + \mu^{m-1}  (\beta_{m-1} -\beta_{m})$ where $p'_{c_{0}}(\beta_{0}) =  0$, $p'_{c_{0}}(\beta_{1}) \neq 0, \cdots, p'_{c_{0}}(\beta_{m-1}) \neq 0$.
\end{notation}

\begin{df} \label{K}
Since $\beta_{m-1} -\beta_{m} \neq 0 $, $w_{m-1}$ is a polynomial of degree $(m-1)$ in the variable $\mu$ so we fix some $\mu \in \mathcal{B}_{\mu}$ such that $w_{m-1} = w_{m-1}(\mu) \neq 0$.
\end{df}
 
\begin{notation} 
We denote for every $\sigma \in \mathbb{C}$, $0 \le n \le m-1$:
 $$M_{n}^{\sigma} =  \begin{pmatrix}
   p'_{c_{0}}(\beta_{n}) & 0 & \sigma (\beta_{n} -\beta_{m})  \\
   1 & 0& 0 \\
   \lambda & 0 & \mu
\end{pmatrix} \, . $$
\end{notation}

\begin{prop}
 We have $M^{\sigma}_{m-1} \cdots M^{\sigma}_{0} \cdot  ( 0 , 0 , 1 )  =  ( \zeta_{1}(\sigma), \zeta_{2}(\sigma) , \zeta_{3}(\sigma) )$, where $\zeta_{1},\zeta_{2},\zeta_{3}$ are holomorphic functions such that $\zeta_{1}(\sigma) = w_{m-1} \cdot \sigma+O(\sigma^{2})$ and $\zeta_{3}(\sigma) = \mu^{m}+O(\sigma)$.  \end{prop}

\begin{proof} 
It is a straightforward consequence of Definition \ref{K}.
\end{proof}

Simple calculations yield the following corollary (the important fact here is that $ p'_{c_{0}}(\beta_{0})=0$ since $\beta_{0} = 0$ is the critical point of $p_{c_{0}}$). 

\begin{corollary} \label{reffin}
Let $\epsilon^{n}_{1}, \epsilon^{n}_{2}, \epsilon^{n}_{3}$ be three holomorphic functions such that $\epsilon^{n}_{1}(\sigma) = O(\sigma)$, $\epsilon^{n}_{2}(\sigma) = O(\sigma^{2})$ and $\epsilon^{n}_{3}(\sigma) = O(\sigma)$ for $0 \le n \le m-1$. Let us denote for every $\sigma \in \mathbb{C}$, $0 \le n \le m-1$: $$N_{n}^{\sigma} =  \begin{pmatrix}
p'_{c_{0}}(\beta_{n}) + \epsilon^{n}_{1}(\sigma) & \epsilon^{n}_{2}(\sigma) & \sigma  (\beta_{n} +\epsilon^{n}_{3}(\sigma)-\beta_{m})  \\
1 & 0& 0 \\
\lambda & 0 & \mu
\end{pmatrix} \, . $$
For every holomorphic maps $\xi_{1},\xi_{2}$ such that $\xi_{1}(\sigma) = O(\sigma)$ and $\xi_{2}(\sigma) = O(\sigma)$, we get:  
$$N^{\sigma}_{m-1} \cdots N^{\sigma}_{0} \cdot  ( \xi_{1}(\sigma) , \xi_{2}(\sigma) , 1)  = ( \zeta_{1}(\sigma), \zeta_{2}(\sigma) , \zeta_{3}(\sigma) ) \, ,$$ where $\zeta_{1},\zeta_{2},\zeta_{3}$ are holomorphic functions such that $\zeta_{1}(\sigma) = w_{m-1} \cdot \sigma+O(\sigma^{2})$ and $\zeta_{3}(\sigma) = \mu^{m}+O(\sigma)$.  
\end{corollary}

\section{Construction of a blender}

In this section, we construct a polynomial automorphism $f_{0}$ of $\mathbb{C}^{3}$. We show that $f_{0}$ has a horseshoe $\mathcal{H}_{f_{0}}$ and that $\mathcal{H}_{f_{0}}$ is a complex blender.

\subsection{Three complex dimensions: the map $f_{0}$}

We recall that $\mathcal{C}$ and $p_{c}$ were defined in Proposition \ref{marber}. We consider now the 3 dimensional map $f_{0}  (z_{1},z_{2},z_{3}) = (p_{c}(z_{1})+b z_{2}+\sigma   z_{3} (z_{1}-\alpha_{c_{0}}),z_{1},\lambda z_{1} + \mu z_{3}+\nu)  $  introduced in Eq. \eqref{e}. It is clear that it is a polynomial automorphism for $c \in \mathcal{C}$ and $b \neq 0$. In the following, we will see that the first direction is expanded by $f_{0}$ and corresponds to the direction of the unstable manifolds of a hyperbolic set we are going to describe. The second and third directions are contracted by $f_{0}$ and correspond to the directions of the stable manifolds of this hyperbolic set.

\begin{notation} \label{con}
We define the following constant cone fields: $C^{u} = \{v = (v_{1},v_{2},v_{3}) \in \mathbb{C}^{3} : \max(|v_{2}|,|v_{3}|)  \le  \chi^{-1}  \cdot |v_{1}|\}$, $C^{ss} = \{v = (v_{1},v_{2},v_{3}) \in \mathbb{C}^{3} : \max(|v_{1}|,|v_{3}|) \le 10^{-6} \cdot |v_{2}|\}$ and
$C^{cs} = \{v = (v_{1},v_{2},v_{3}) \in \mathbb{C}^{3} : \max(|v_{1}|,|v_{2}|) \le 10^{-6} \cdot |v_{3}|\}$, where the constant  $\chi>1$ was defined in Proposition \ref{marber}.
\end{notation}

We now give a non general definition of a horseshoe which is specific to our context.

\begin{df} \label{df7}
Given an automorphism $F : \mathbb{C}^{3} \rightarrow \mathbb{C}^{3}$, a tridisk $D = D_{1} \times D_{2} \times D_{3} \subset \mathbb{D}^{3}$ and an integer $p \ge 1$, we say that $\mathcal{H}_{F} =\bigcap_{n \in \mathbb{Z}} F^{n}(\overline{D})$ is a $p$-branched horseshoe for $F$ if: 
\begin{enumerate} 
\item $F(D) \cap D$ has $p$ components $D^{j,u}$ which do not intersect $D_{1} \times \partial (D_{2} \times  D_{3})$,
 \item $F^{-1}( D) \cap D$ has $p$ components $D^{j,s}$ which do not intersect $\partial D_{1} \times D_{2} \times D_{3}$,
 \item  on $\bigcup_{1 \le j \le p} D^{j,s}$, the cone field $C^{u}$ is $F$-invariant, and on $\bigcup_{1 \le j \le p} D^{j,u}$ the cone field $C^{ss}$ is $F^{-1}$-invariant. Moreover there exists $\Xi>1$ such that the cone field $\{ ( v_{1},v_{2},v_{3}) : ||(v_{2},v_{3})|| > \Xi ||v_{1}||  \}$ contains $C^{cs}$ and is $F^{-1}$-invariant on $\bigcup_{1 \le j \le p} D^{j,u}$,
\item  there exists $C_{F}>1$ such that at every point of $\bigcup_{1 \le j \le p} D^{j,s}$, for every non zero  $v \in C^{u}$, we have  $||DF(v)||>C_{F}  ||v||$, and at every point of $\bigcup_{1 \le j \le p} D^{j,u}$, for every non zero $v \in C^{ss} \cup \{ ( v_{1},v_{2},v_{3}) : ||(v_{2},v_{3})|| > \Xi||v_{1}||  \}$, we have $ ||D(F^{-1})(v)||>C_{F}  ||v||$. \end{enumerate}
\end{df}

\begin{prop}
If $\mathcal{H}_{F} = \bigcap_{n \in \mathbb{Z}} F^{n}(  \overline{D} )$ is a $p$-branched horseshoe, then it is a horseshoe in the classical meaning of this term, that is a compact, invariant, transitive, hyperbolic set. \end{prop}

\begin{proof}
The set $\bigcap_{n \in \mathbb{Z}} F^{n}( \overline{ D} )$ is compact as an intersection of compact sets and $F$-invariant by definition. Moreover, one can take the (non necessarily invariant) decomposition $\mathbb{C}^{3} \simeq \mathbb{R}^{6} =  \mathbb{C} \bigoplus \mathbb{C}^{2} \simeq \mathbb{R}^{2} \bigoplus \mathbb{R}^{4}$ and the associated constant cone fields $\mathcal{C}^{u}_{\mathbb{R}} = \{ ( v_{1},v_{2},v_{3}) : ||v_{1}||>\chi||(v_{2},v_{3})||\}$ and $\mathcal{C}^{s}_{\mathbb{R}} = \{ ( v_{1},v_{2},v_{3}) : ||(v_{2},v_{3})|| > \Xi ||v_{1}||  \}$. The definition above implies that both $\mathcal{C}^{u}_{\mathbb{R}}$ is $F$-invariant and $\mathcal{C}^{s}_{\mathbb{R}}$ is $F^{-1}$-invariant. Moreover, they are expanded by a factor $C_{F}$ larger than 1 respectively under $F$ and $F^{-1}$. Besides, the sets  $D^{j,u}$ do not intersect $D_{1} \times \partial (D_{2} \times  D_{3})$ and the sets $D^{j,s}$ do not intersect $\partial D_{1} \times D_{2} \times D_{3}$. Then $\bigcap_{n \in \mathbb{Z}} F^{n}( \overline{ D} )$ is a horseshoe in the sense of Definition 6.5.2 of \cite{caga}. According to the discussion following this definition, $\bigcap_{n \in \mathbb{Z}} F^{n}( \overline{ D} )$ is hyperbolic (this is also a straightforward application of the cone field criterion, Corollary 6.4.8 in \cite{caga}) and is topologically conjugate to a shift. In particular, it is transitive. This ends the proof of the proposition.
\end{proof}

\begin{remark} \label{remarkhorsehoesqure} 
For our definition, a $p$-branched horseshoe has one unstable direction and two stable directions. It is also straightforward that if $F$ has a $p$-branched horseshoe, then $F^{2}$ has a $p^{2}$-branched horseshoe. \end{remark}

\begin{df}
We say that a saddle periodic point of multipliers $| \lambda_{1}|  \le | \lambda_{2}|  < 1 <|  \lambda_{3}| $ is sectionally dissipative if the product of any two of its eigenvalues is less than 1 in modulus, that is, $| \lambda_{1}\lambda_{3}|  < 1$ and $| \lambda_{2}\lambda_{3}| < 1$.
\end{df}

In the next proposition, we prove that if $b$ and $\sigma$ are sufficiently small, then some iterate of $f_{0}= f_{c,b,\sigma}$ has a $2$-branched horseshoe. Moreover, we introduce a neighborhood $\mathcal{F}$ of $f_{0}$ where this property persists. In Section 5, we will make a particular choice of $c$, $b$ and $\sigma$ so that the stable manifold of a periodic point of $f_{0}$ will have special properties, which will persist in a new neighborhood $\mathcal{F}' \Subset \mathcal{F}$ of $f_{0}$ in $\mathrm{Aut}_{2}(\mathbb{C}^{3})$.

\begin{prop} \label{ferm}
Let $q \ge 100$ and $\mathcal{C}$,$r$ given by Proposition  \ref{marber}. Let $f_{0} = f_{c,b,\sigma}$ be the polynomial automorphism of $\mathbb{C}^{3}$ introduced in Eq. (1). Then, there exists $10^{-10}>b_{0}>0$ and $10^{-10}>\sigma_{0}>0$ independent of $c \in \mathcal{C}$  such that if $0<|b| < b_{0}$ and $0 \le |\sigma| < \sigma_{0}$, then $\mathcal{H}_{f_{0}} = \bigcap_{n \in \mathbb{Z}} (f^{qr}_{0})^{n}( \overline{ \mathbb{D}^{3}} )$ is a $2$-branched horseshoe. Moreover, $f_{0}$ has a periodic point $\delta_{f_{0}}$ that is sectionally dissipative and belongs to the homoclinic class of the continuation $\alpha_{f_{0}}$ of $\alpha_{c}$. These results remain true for any $f$ in some neighborhood $\mathcal{F} = \mathcal{F}(c,b,\sigma)$ of $f_{0}$ in $\mathrm{Aut}_{2}(\mathbb{C}^{3})$.
\end{prop}

\begin{proof} 
When $\sigma = b = 0$, $f_{c,0,0} : (z_{1},z_{2},z_{3}) \mapsto (p_{c}(z_{1}),z_{1},\lambda  z_{1} + \mu  z_{3}+\nu)$ satisfies:
 $$f^{qr}_{c,0,0}(z_{1},z_{2},z_{3})=(p^{qr}_{c}(z_{1}),p^{qr-1}_{c}(z_{1}), \nu + \lambda p_{c}^{qr-1}(z_{1})+ \mu  (\nu+ \cdots +\mu  (\nu+  \lambda z_{1} )) + \mu^{qr} z_{3} )  .$$
By Proposition \ref{marber} $(1)$,  $p_{c}^{r-1} (\mathbb{D}'_{1}), p_{c}^{r-1}(\mathbb{D}'_{2})  \Subset \mathbb{D}$ and so  $p_{c}^{qr-1} (h_{1}^{q}(\mathbb{D}')),p_{c}^{qr-1}(h_{2}^{q}(\mathbb{D}'))  \Subset \mathbb{D}$. Then for any $(z_{1},z_{2},z_{3}) \in  \big(   h_{1}^{q}(\mathbb{D}')   \cup h_{2}^{q}(\mathbb{D}') \big) \times \overline{\mathbb{D}}^{2}$, the second coordinate of $f_{c,0,0}(z_{1},z_{2},z_{3})$ is in $\mathbb{D}$. According to Corollary \ref{hyu}, it holds  $\nu + \lambda p_{c}^{qr-1}(z_{1})+ \mu  (\nu+ \cdots +\mu  (\nu+  \lambda z_{1} )) \in \mathbb{D}( \pm \frac{9}{10} \cdot 10^{-4}, 10^{-10})$ for every  $(z_{1},z_{2},z_{3}) \in  \big(   h_{1}^{q}(\mathbb{D}')   \cup h_{2}^{q}(\mathbb{D}') \big) \times \overline{\mathbb{D}}^{2}$. Since  $|\mu|^{qr} < 1-  10^{-4}+10^{-10}$ (see Remark \ref{rmu}), the third coordinate of $f_{c,0,0}(z_{1},z_{2},z_{3})$ lies in $\mathbb{D}$ too. Then the intersection  $f^{qr}_{c,0,0} \big(   \big(   h_{1}^{q}(\mathbb{D}')   \cup h_{2}^{q}(\mathbb{D}') \big) \times \overline{\mathbb{D}}^{2}  \big) \cap \overline{\mathbb{D}^{3}}$ has two components which are two graphs over the coordinate $z_{1} \in \overline {\mathbb{D}}$ and which do not intersect $\overline{\mathbb{D}} \times \partial (  \overline{\mathbb{D}} \times \overline{\mathbb{D)}}$. \medskip

By continuity, there exists $10^{-10}>\sigma_{0}>0$ such that for every $0 \le |\sigma| < \sigma_{0}$, $f^{qr}_{c,0,\sigma} \big(   \big(   h_{1}^{q}(\mathbb{D}')   \cup h_{2}^{q}(\mathbb{D}') \big) \times \overline{\mathbb{D}}^{2}  \big) \cap \overline{\mathbb{D}^{3}}$ has still two components which are graphs over $z_{1} \in  \overline{\mathbb{D}}$. The distance to $\overline{\mathbb{D}} \times \partial (  \overline{\mathbb{D}} \times \overline{\mathbb{D)}}$ of each of these two graphs is bounded from below by a strictly positive constant  independent of $0 \le |\sigma| < \sigma_{0}$. Thus there exists $10^{-10}>b_{0}>0$ such that for every  $0 \le |\sigma| < \sigma_{0}$ and  $0<|b| < b_{0}$, the intersection $f^{qr}_{c,b,\sigma} \big(    \big(  h_{1}^{q}(\mathbb{D}')   \cup h_{2}^{q}(\mathbb{D}') \big) \times \overline{\mathbb{D}}^{2} \big)  \cap  \overline{\mathbb{D}}^{3}$ has two connected components $D^{1,u}$ and $D^{2,u}$. These are two tridisks which do not intersect $\overline{\mathbb{D}} \times \partial (  \overline{\mathbb{D}} \times \overline{\mathbb{D)}}$ and are also connected components of $f^{qr}_{c,b,\sigma}( \overline{  \mathbb{D}^{3}}) \cap \overline{  \mathbb{D}^{3}}$. This proves item 1 of Definition \ref{df7}. We denote now $f_{c,b,\sigma} = f_{0}$. The proof of item 2 is similar and gives two connected components $D^{1,s}$ and $D^{2,s}$ in $f^{-qr}_{0}( \overline{  \mathbb{D}^{3}}) \cap \overline{  \mathbb{D}^{3}}$. Reducing $\sigma_{0}$ and $b_{0}$ if necessary, the cone field $C^{u}$ is $f_{0}^{qr}$-invariant on $D^{1,s} \cup D^{2,s}$ by Property (3) of Proposition \ref{marber} and dilated under $f_{0}^{qr}$ by a factor larger than  $10^{10}$ because $|h_{j}'|<10^{-10}$ on $\mathbb{D}''$ (see Property 1 of Proposition \ref{fo}). The modulus of the derivative of $p_{c}$ is bounded from below by a strictly positive constant on $\mathbb{D}' \cup p_{c}(\mathbb{D}') \cup \cdots \cup p_{c}^{qr-1}(\mathbb{D}')$. Thus, still reducing  $\sigma_{0}$ and $b_{0}$ if necessary, the cone field $C^{ss}$ is $f_{0}^{-qr}$-invariant on $D^{1,u} \cup D^{2,u}$ and dilated under $f_{0}^{-qr}$ by a factor at least $10^{10}$. For $\sigma \rightarrow  0$  and $b \rightarrow 0$, the directions of the two expanding eigenvectors of $D_{z}f^{-qr}_{0}$ respectively tend to the second and third directions (uniformly in $z$) while the contracting eigenvector of $D_{z}f^{-qr}_{0}$ stays in $C^{u}$. Since the associated rates of dilatation/contraction remain uniformly distant to 1 when $\sigma,b \rightarrow 0$, this implies that some cone of the form  $\{ ( v_{1},v_{2},v_{3}) : ||(v_{2},v_{3})|| > \Xi  ||v_{1}||  \}$ (for some  $\Xi>1$) is $f^{-qr}_{0}$-invariant on $D^{1,u} \cup D^{2,u}$ and dilated by a factor larger than 1. This finishes to prove both items 3 and 4 of Definition \ref{df7}. Then $\mathcal{H}_{f_{0}} = \bigcap_{n \in \mathbb{Z}} (f^{qr}_{0})^{n}( \overline{\mathbb{D}^{3} })$ is a $2$-branched horseshoe. \medskip

The multiplier of the periodic point $\delta_{c}$ has its modulus between 1 and $(1+10^{-10})^{1/qr}$ by Property 6 of Proposition \ref{marber}. Then its continuation $\delta_{f_{0}}$ is a saddle point of expanding eigenvalue between 1 and $(1+2 \cdot 10^{-10})^{1/qr}$ after reducing $b_{0}$ and $\sigma_{0}$ if necessary. It has two contracting eigenvalues. When $b=\sigma= 0$, one is equal to 0 and the other one is a power of $\mu$ with $|\mu| < (1-10^{-4}+10^{-10})^{\frac{1}{qr}}$ (see Remark \ref{rmu}). Then, by continuity, reducing $b_{0}$ and $\sigma_{0}$ if necessary, $\delta_{f_{0}}$ is sectionally dissipative. We denote by $W^{u}_{\text{loc}}(\alpha_{f_{0}})$ (resp. $W^{s}_{\text{loc}}(\alpha_{f_{0}})$) the connected component of $W^{u}(\alpha_{f_{0}}) \cap \mathbb{D}^{3}$ (resp. $W^{s}(\alpha_{f_{0}}) \cap \mathbb{D}^{3}$) which contains $\alpha_{f_{0}}$, and we use the same notation for other periodic points. For small values of $b$ and $\sigma$, both $W^{u}_{\text{loc}}(\delta_{f_{0}})$ and $W^{u}_{\text{loc}}(\alpha_{f_{0}})$ are graphs over $z_{1} \in \mathbb{D}$ and both $W^{s}_{\text{loc}}(\delta_{f_{0}})$ and $W^{s}_{\text{loc}}(\alpha_{f_{0}})$ are graphs over $(z_{2},z_{3}) \in  \mathbb{D}^{2}$. Reducing $b_{0}$ and $\sigma_{0}$ if necessary, $W^{u}_{\text{loc}}(\alpha_{f_{0}})$ and $W^{s}_{\text{loc}}(\delta_{f_{0}})$ intersect and $W^{s}_{\text{loc}}(\alpha_{f_{0}})$ and $W^{u}_{\text{loc}}(\delta_{f_{0}})$ intersect so $\delta_{f_{0}}$ is in the homoclinic class of $\alpha_{f_{0}}$ . \medskip

Finally, since all the conditions of the proposition are open, these results remain true for any $f$ in some neighborhood $\mathcal{F} = \mathcal{F}(c,b,\sigma)$ of $f_{0}$ in $\text{Aut}(\mathbb{C}^{3})$. This concludes the proof of the proposition.
\end{proof}

\subsection{Blender property}

In this subsection, we show that the third coordinate of $f \in \mathcal{F}$ has some kind of open-covering property. Since $f^{qr} \in \mathcal{F}$ has a 2-branched horseshoe, $f^{2qr}$ has a 4-branched horseshoe by Remark  \ref{remarkhorsehoesqure}. We are interested in the geometric properties of the third coordinate of $f^{2qr}$. Here, we use it to show that the 4-branched horseshoe associated to $f^{2qr}$ is a blender. Remind that the maps $h_{1}$ and $h_{2}$ were defined in Notation \ref{defh1h}.

\begin{df} \label{uu}
For $f \in \mathcal{F}$, we denote by $f^{qr}[j]$ the restriction of $f^{qr}$ on $h^{q}_{j}( \mathbb{D}') \times  \mathbb{D}^{2}$ for $j \in \{1,2\}$. We put $V_{1} = f^{qr}[1]   (h^{q}_{1}( \mathbb{D}') \times  \mathbb{D}^{2}) \cap \mathbb{D}^{3}$ and $V_{2} =f^{qr}[2]   (h^{q}_{2}( \mathbb{D}') \times  \mathbb{D}^{2}) \cap  \mathbb{D}^{3}$. We define $U_{1} = f^{qr}[1] \big(   V_{1}  \cap   (h^{q}_{1}( \mathbb{D}') \times  \mathbb{D}^{2})   \big) \cap \mathbb{D}^{3}$, $U_{2} = f^{qr}[2]   \big(   V_{1} \cap   (h^{q}_{2}( \mathbb{D}') \times  \mathbb{D}^{2})       \big) \cap \mathbb{D}^{3} $, $U_{3} = f^{qr}[1]   \big(   V_{2} \cap    (h^{q}_{1}( \mathbb{D}') \times  \mathbb{D}^{2})     \big) \cap \mathbb{D}^{3}$ and $U_{4} = f^{qr}[2]   \big( V_{2} \cap    (h^{q}_{2}( \mathbb{D}') \times  \mathbb{D}^{2})         \big) \cap \mathbb{D}^{3}$. We define $g_{j} = f^{-2qr}_{|U_{j}}$ for $1 \le j \le 4$. 
\end{df}

\begin{notation}
We will denote $c_{1} = \frac{9}{10} \cdot 10^{-4} \cdot (1+i)$, $c_{2} = \frac{9}{10} \cdot 10^{-4} \cdot (-1+i)$, $c_{3} =  \frac{9}{10} \cdot 10^{-4} \cdot (1-i)$ and  $c_{4} =  \frac{9}{10} \cdot 10^{-4} \cdot (-1-i)$.
\end{notation}

\begin{lemma} \label{esti}
Reducing $b_{0}$, $\sigma_{0}$ and $\mathcal{F}$ if necessary, we have: 
\begin{enumerate} \item  $ \forall z \in f^{-qr}(V_{1}), \mathrm{pr}_{3}( f^{qr}(z))   \in  \mathbb{D}( \mu^{qr} z_{3} + \frac{9}{10} \cdot 10^{-4},10^{-9})$, \item  $ \forall z \in f^{-qr}(V_{2}), \mathrm{pr}_{3}( f^{qr}(z))   \in  \mathbb{D}( \mu^{qr} z_{3} - \frac{9}{10} \cdot 10^{-4},10^{-9})$. \end{enumerate}
\end{lemma}

\begin{proof}
For $\sigma = b = 0$, it is a simple consequence of Corollary \ref{hyu}. Then, by continuity we just have to take sufficiently small values of $\sigma_{0},b_{0},\mathcal{F}$ to get the bound $10^{-9}$.
\end{proof}

Since we both have $\mu_{0}^{qr} = (1-10^{-4}) \cdot e^{i \cdot \frac{\pi}{2}}$ (by Notation \ref{h}) and  $\mu^{qr} \subset \mathbb{D}(\mu_{0}^{qr},10^{-10})$ (by Remark \ref{rmu}), by iterating two times the previous result, we get the following: 

\begin{corollary} \label{estit}
For every $z=(z_{1},z_{2},z_{3}) \in U_{j}$, $ \mathrm{pr}_{3}( g_{j}(z)) \in  \mathbb{D}( \frac{1}{\mu^{2qr}}   (z_{3} - c_{j}), 10^{-6})$.
\end{corollary}

We now show an open covering property for the four affine maps $z \mapsto \frac{1}{\mu^{2qr} }  (z - c_{j})$.

\begin{prop} \label{rp9}
For every $z \in  \mathbb{D}(0,\frac{1}{10})$, there exists $j \in \{1,2,3,4\}$ such that: $$\frac{1}{\mu^{2qr} }  (z - c_{j}) \in  \mathbb{D}(0,\frac{1}{10}-10^{-5}) \, .$$
\end{prop}

\begin{proof}
We check that the union of the images of $ \mathbb{D}(0,\frac{1}{10}-10^{-5})$ under the four affine maps $z \mapsto \mu^{2qr} z+  c_{j} $ contains $ \mathbb{D}(0,\frac{1}{10})$. We begin by showing that for every point $z$ of the set $\{z : |z| = \frac{1}{10} \text{ and } 0 \le \text{Arg}(z) \le \frac{\pi}{2}\}$, the point $z-c_{1}$ belongs to the disk $\mathbb{D} \big( 0,|\mu|^{2qr}  (\frac{1}{10} -10^{-5}) \big)$. Let us point out that $|\mu|^{2qr}  > 1-2 \cdot 10^{-4}$ (see Remark \ref{rmu}). Denoting $z = x+ iy$ with $(x,y) \in \mathbb{R}_{+}^{2}$, we have $ |z-c_{1}|^{2} = (x-\frac{9}{10} \cdot 10^{-4})^{2} +  (y-\frac{9}{10} \cdot 10^{-4})^{2}$. Since $x^{2}+y^{2} = (\frac{1}{10})^{2}$, at least one between $x$ or $y$ is larger than $\frac{1}{\sqrt{2}} \cdot \frac{1}{10}$. Then $|z-c_{1}|^{2}$ is smaller than: 
  $$(\frac{1}{10})^{2} -2 \frac{1}{\sqrt{2}}  \frac{1}{10}  \frac{9}{10} 10^{-4} +2   (\frac{9}{10}  10^{-4})^{2}<\big( ( 1-2 \cdot 10^{-4})  (\frac{1}{10} -10^{-5})\big)^{2}< \big( |\mu|^{2qr} (\frac{1}{10} -10^{-5}) \big)^{2}$$ Then  for every point $z$ of the set $\{z : |z| = \frac{1}{10} \text{ and } 0 \le \text{Arg}(z) \le \frac{\pi}{2}\}$, we have that $z-c_{1} \in \mathbb{D} \big( 0,|\mu|^{2qr}(\frac{1}{10} -10^{-5}) \big)$. We also have $0-c_{1} \in  \mathbb{D} \big( 0,|\mu|^{2qr}(\frac{1}{10} -10^{-5}) \big)$. Thus by convexity the image of $ \mathbb{D}(0,\frac{1}{10}-10^{-5})$ by the affine map $z \mapsto \mu^{2qr}  z+  c_{1} $ contains the first quadrant of $ \mathbb{D}(0,\frac{1}{10})$. Then the result follows by symmetry. 
\end{proof}

\includegraphics[width=11cm]{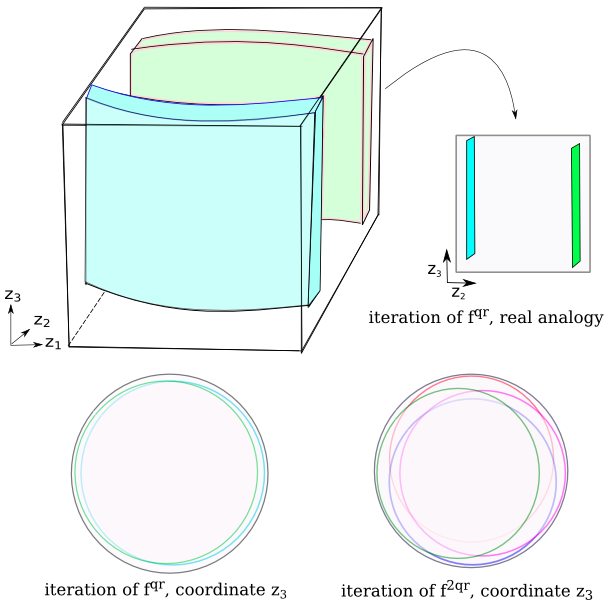}

\begin{center} 
Figure 1: complex blender. The top figure shows the sets $V_{1}$ and $V_{2}$ in a real analogy (that is in $\mathbb{R}^{3}$). The two bottom figures represent the respective images by the third projection map  $\pi_{3} :  \mathbb{C}^{3} \rightarrow \mathbb{C}$ of $V_{1}$,$V_{2}$ (on the left) and $U_{1}$,$U_{2}$,$U_{3}$,$U_{4}$ (on the right). 
\end{center}

We are finally in position to prove the main result of this section.

\begin{df} \label{curves}
A $ss$-curve $\Gamma$ is a holomorphic graph over  $z_{2}$ which has all its tangent vectors in $C^{ss}$ (recall that this cone was defined in Notation \ref{con}).
\end{df} 

\begin{prop} \label{blender}
Let $\Gamma $ be any $ss$-curve intersecting $\mathbb{D}^{2} \times  \mathbb{D}(0,\frac{1}{10})$. Then for every $f \in \mathcal{F}$, $\Gamma$ intersects the unstable manifold of a point of the horseshoe $\mathcal{H}_{f}$: in other words, $\mathcal{H}_{f}$ is a blender.
\end{prop}

\begin{proof}  
Let $\Gamma = \Gamma_{0}$ be a $ss$-curve intersecting $\mathbb{D}^{2} \times  \mathbb{D}(0,\frac{1}{10})$ and $f \in \mathcal{F}$. We show that there exists $j \in \{1,2,3,4\}$ such that $g_{j} (\Gamma)$ contains a $ss$-curve $\Gamma^{1}$ intersecting $\mathbb{D}^{2} \times  \mathbb{D}(0,\frac{1}{10})$. Let $Z = (Z_{1},Z_{2},Z_{3})$ be a point of $\Gamma \cap (\mathbb{D}^{2} \times  \mathbb{D}(0,\frac{1}{10}))$. By Lemma \ref{rp9} there exists $j \in \{1,2,3,4\}$ such that $\frac{1}{\mu^{2qr} }(Z_{3}  -c_{j}) \in  \mathbb{D}(0,\frac{1}{10}-10^{-5})$. Since $ \mathbb{D}(0,\frac{1}{10}) \subset \text{pr}_{3}(U_{j})$, by continuity, $\Gamma$ intersects $U_{j}$. The cone $C^{ss}$ is $g_{j}$-invariant since $\mathcal{H}_{f}$ is a 2-branched horseshoe, and so $\Gamma^{1} = g_{j}(\Gamma \cap U_{j})$ is a $ss$-curve. Then it holds  $\text{diam}(\text{pr}_{3}(\Gamma^{1}))   \le 10^{-6}$. By Corollary \ref{estit}, $|  \text{pr}_{3} ( g_{j}(Z) ) - \frac{1}{\mu^{2qr} }(Z_{3}-c_{j})|<10^{-6} $. This implies that $\text{pr}_{3}(  \Gamma^{1} )  \subset     \mathbb{D}(0,\frac{1}{10}-10^{-5}+2 \cdot 10^{-6})    \subset  \mathbb{D}(0,\frac{1}{10})$. Then  $\Gamma^{1}$ intersects $\mathbb{D}^{2} \times  \mathbb{D}(0,\frac{1}{10})$. By iteration, we can construct a sequence of $ss$-curves $\Gamma^{n} \subset  (g_{j_{n}}  \circ ... \circ g_{j_{1}})(\Gamma)$, each of them  intersecting $\mathbb{D}^{3}$. The set $\bigcap_{n \ge 1} (f^{2qr})^{n}(\Gamma^{n}) \subset \Gamma$ is non empty since it contains the intersection of a decreasing sequence of non empty compact sets. But any point in this intersection is a point of the unstable manifold of a point of the horseshoe $\mathcal{H}_{f}$. This ends the proof.
\end{proof}

\section{Mechanism to get persistent tangencies}

In this section, we explain how the blender property obtained in the last section leads to persistent tangencies with certain "folded" surfaces.

\subsection{Some definitions}

\begin{df}
A submanifold (or an analytic set) $W \subset \mathbb{D}^{n}$ is horizontal relatively to a decomposition $\mathbb{D}^{n} = \mathbb{D}^{k} \times \mathbb{D}^{n-k}$ if $W$ does not intersect $ \mathbb{D}^{k} \times  \partial \mathbb{D}^{n-k}$. We will also say it is horizontal relatively to the projection $\mathbb{D}^{n} \rightarrow \mathbb{D}^{k}$. If $\dim(W) = k$, then the natural projection on $\mathbb{D}^{k}$ is a branched covering of degree $d$. We will say that $W$ is of degree $d$. We similarly define vertical submanifolds in $\mathbb{D}^{k} \times \mathbb{D}^{n-k}$ and their degree.
\end{df}

The two following propositions are classical. For a proof, one can refer to \cite{dsfdm}. One can also refer to \cite{henam}.
  
\begin{prop} \label{fold55} 
Let $W$ be a horizontal curve of degree $d$ and $W'$ be a vertical curve of degree $d'$ in $\mathbb{D}^{2} = \mathbb{D} \times \mathbb{D}$. Then $W$ and $W'$ intersect in $dd'$ points with multiplicity.  
\end{prop}

\begin{prop} \label{fold77}
Let $W$ be a horizontal curve (resp. surface) of degree $d$ and $W'$ be a vertical surface (resp.curve) of degree $d'$ in $\mathbb{D}^{3} = \mathbb{D}^{1} \times \mathbb{D}^{2}$ (resp. $\mathbb{D}^{3} = \mathbb{D}^{2} \times \mathbb{D}^{1}$). Then $W$ and $W'$ intersect in $dd'$ points with multiplicity.  
\end{prop}
   
\begin{remark}  In particular, in the two previous propositions, if all the intersections are transverse, there are exactly $dd'$ distinct points of intersection. If it is not the case, there is at least one point of tangency. \end{remark}
 
We now introduce several definitions specific to our context.
   
\begin{notation}
Let $i,j$ be two distinct integers in $\{1,2,3\}$. We will denote by $\pi_{i}$ the projection over the $i^{th}$ coordinate and by $(\pi_{i},\pi_{j})$ the projection over the $i^{th}$ and the $j^{th}$ coordinates.
\end{notation}

\begin{df}  
Let $i,j$ be two distinct integers in $\{1,2,3\}$. A $(i,j)$-surface $\mathcal{S}$ is a complex surface horizontal relatively to the projection $(\pi_{i},\pi_{j})$.
\end{df} 

\begin{remark}
A $(i,j)$-surface $\mathcal{S}$ is a ramified covering of degree $d$ over $\mathbb{D}^{2}$. In the rest of this article, we will only consider ramified coverings of degree 1 or 2.
\end{remark}

\begin{df} \label{curves}
A $u$-curve $\Gamma$ is a holomorphic graph over $z_{1}$ which has all its tangent vectors in $C^{u}$ (recall that this cone was defined in Notation \ref{con}).
\end{df}

\begin{df}  
Let $k \in \{2,3\}$. A $(1,k)$-quasi plane $\mathcal{V}$ over a bidisk $\mathbb{D}_{1} \times \mathbb{D}_{k}$ is a graph $\{z_{k'} = v(z_{1},z_{k})\} $, where $k'$ is the integer defined by $\{1,2,3\} = \{1,k,k'\}$ and $v : \mathbb{D}_{1} \times \mathbb{D}_{k} \rightarrow \mathbb{D}$ is holomorphic, such that the $C^{0}$-norm of $Dv$ is bounded by 1 and $\mathcal{V}$ is foliated by $u$-curves $\mathcal{V}_{x}$.   \end{df}

\begin{remark}
In particular, every $(1,k)$-quasi plane over $\mathbb{D}^{2}$ is a $(1,k)$-surface.
\end{remark}
   
\begin{df}
Let $k \in \{2,3\}$. A $k$-folded curve is a holomorphic curve horizontal relatively to $\pi_{k}$ which is a ramified 2-covering over $z_{k}$ with exactly one point of ramification $z_{k,ram}$. We denote $\mathrm{fold}(\Gamma) = z_{k,ram}$ the \emph{fold} of $\Gamma$.
\end{df}
   
\begin{df}
Let $k \in \{2,3\}$. A $k$-folded $(2,3)$-surface $\mathcal{W}$ is a complex surface of degree 2 that is horizontal relatively to $(\pi_{2},\pi_{3})$ and such that for every $(1,k)$-quasi plane $\mathcal{V}$ over $\mathbb{D}^{2}$, $\Gamma = \mathcal{V} \cap \mathcal{W}$ is a $k$-folded curve. We denote: $$\mathrm{Fold}(\mathcal{W}) = \{ \mathrm{fold}(\Gamma) : \Gamma \text{ is a } k\text{-folded curve included in }  \mathcal{W}\}$$ which is a subset of $\mathbb{D}$. We say that $\mathcal{W}$ is concentrated if $\mathrm{diam}(  \mathrm{Fold}(\mathcal{W}) ) \le 10^{-5}$.    
\end{df}

\subsection{Preparatory lemmas}

In this Subsection we gather some simple technical results that will be useful in the following Subsections.

\begin{lemma} \label{lemb}
Let $\Gamma = \gamma(\mathbb{D})$ be a $k$-folded curve (with $k \in \{2,3\}$) included inside some  $(1,k)$-quasi plane over $\mathbb{D} \times \mathbb{D}_{k}$ with $\mathbb{D}(0,1/2) \subset \mathbb{D}_{k} \subset \ \mathbb{D}$ and $\mathrm{diam} \big(  \mathrm{pr}_{1}(\Gamma) \big) \le 10^{-10}$. Then for every disk $\mathbb{D}_{\Gamma} \Subset \mathbb{D}_{k}$ of radius $10^{-7}$ distant of at least $10^{-7}$ from $\partial \mathbb{D}_{k} \cup \{\mathrm{fold}(\Gamma)\}$, it holds: for every $Z$ with $\mathrm{pr}_{k}(\gamma(Z)) \in \mathbb{D}_{\Gamma}$ we have $|\gamma'_{1}(Z)| < 10^{-3} |\gamma'_{k}(Z) |$. In particular, $\Gamma \cap (\mathbb{D}^{2} \times \mathbb{D}_{\Gamma})$ is the union of two graphs upon $z_{k} \in \mathbb{D}_{\Gamma}$. 
\end{lemma}
   
\begin{proof}
Let $ \mathbb{D}_{\Gamma} \Subset \mathbb{D}_{k}$ be a disk of radius $10^{-7}$ distant of at least $10^{-7}$ from $\partial \mathbb{D}_{k} \cup \{\mathrm{fold}(\Gamma)\}$. We notice that $\Gamma$ is the union of two graphs over $z_{k}$ varying in the $10^{-7}$-neighborhood of $\mathbb{D}_{\Gamma}$. Let us denote $z_{1} = \xi(z_{k})$ one of them. Every point of $ \mathbb{D}_{\Gamma}$ is the center of a ball of radius $10^{-7}$ where $|\xi(z_{k})|<10^{-10}$. Hence, by the Cauchy inequality, we have $|\xi'(z_{k})|<10^{7} \cdot 10^{-10} = 10^{-3}$. Thus  $|\gamma'_{1}(Z)|< 10^{-3} |\gamma'_{k}(Z)|$ at every $Z$ such that $\mathrm{pr}_{k}(\gamma(Z)) \in \mathbb{D}_{\Gamma}$  and the result follows.
\end{proof}
   
Here are some consequences of Propositions \ref{fold55} and \ref{fold77}. 
   
\begin{lemma} \label{rq5}
Let $\mathcal{V}$ be a $(1,k)$-quasi plane over $\mathbb{D}_{1} \times \mathbb{D}_{k}$ with $\mathbb{D}(0,1/2) \subset \mathbb{D}_{k} \subset \ \mathbb{D}$. Let $\Gamma$ be a $u$-curve and $\tilde{\Gamma}$ be a graph over $z_{k}$ with $|\tilde{\gamma}'_{1}| < 10^{-3} |\tilde{\gamma}'_{k}|$, both included in $\mathcal{V}$. Then $\Gamma \cap \tilde{\Gamma}$ is a singleton.
\end{lemma}

\begin{prop} \label{rq11}
Let $\mathcal{W}$ be a $k$-folded $(2,3)$-surface with $k \in \{2,3\}$. Let $\Gamma$ be a $u$-curve. Then $\mathcal{W} \cap \Gamma$ has one or two points.
\end{prop}

\subsection{Main result}

Here is the main result of this section: we show that any concentrated 3-folded $(2,3)$-surface having its fold in good position has a point of tangency with the unstable manifold of a point of the horseshoe $\mathcal{H}_{f}$. 
   
\begin{prop} \label{rrrrr}   
Let $\mathcal{W}$ be a concentrated 3-folded $(2,3)$-surface  such that we both have $\mathrm{Fold}(\mathcal{W}) \subset  \mathbb{D}(0,\frac{1}{10})$ and $\mathrm{diam}(\mathrm{pr}_{1}(\mathcal{W})) \le 10^{-10}$. Then there exists a point $\kappa_{f}$ of the horseshoe $\mathcal{H}_{f}$ such that $\mathcal{W}$ has a point of tangency with the unstable manifold of $\kappa_{f}$.
\end{prop}

We begin by a lemma showing that the image of a folded curve contains a folded curve. Remind that the sets $U_{j}$ and the maps $g_{j}$ were defined in Definition \ref{uu}.

\begin{lemma} \label{lemb5}
Let $\tilde{\Gamma}$ be a $3$-folded curve included in a $(1,3)$-quasi plane $\tilde{\mathcal{V}}$ over $\mathbb{D}_{1} \times \mathbb{D}_{k}$ with $\mathbb{D}(0,1/2) \subset \mathbb{D}_{k} \subset \ \mathbb{D}$. We suppose that $\tilde{\Gamma}$ and $\tilde{\mathcal{V}}$ are included in $U_{j}$ for some  $j \in  \{1,2,3,4\}$. We also suppose $\mathrm{diam}(\mathrm{pr}_{1}(\tilde{\Gamma})) \le 10^{-10}$ and that $g_{j}(\tilde{\mathcal{V}})$ is included in a $(1,3)$-quasi plane $\mathcal{V}$ over $\mathbb{D}^{2}$. Then $\Gamma = g_{j}(\tilde{\Gamma})$ is a $3$-folded curve satisfying: $$\mathrm{fold}(\Gamma) \in \mathbb{D}( \frac{1}{\mu^{2qr}} \big(  \mathrm{fold}(\tilde{\Gamma})- c_{j} ), \frac{3}{2} \cdot 10^{-6} \big) \, .$$
\end{lemma}

\begin{proof}
By Lemma \ref{lemb}, for every disk $\mathbb{D}_{\tilde{\Gamma}} \Subset \mathbb{D}_{k}$ of radius $10^{-7}$  distant of at least $10^{-7}$ from $\partial \mathbb{D}_{k} \cup \{\mathrm{fold}(\tilde{\Gamma})\}$, it holds: for every $Z$ with $\mathrm{pr}_{k}(\tilde{\gamma}(Z)) \in \mathbb{D}_{\tilde{\Gamma}}$ we have $|\tilde{\gamma}'_{1}(Z)| < 10^{-3} |\tilde{\gamma}'_{k}(Z) |$. Moreover $\tilde{\Gamma} \cap (\mathbb{D}^{2} \times \mathbb{D}_{\tilde{\Gamma}})$ is the union of two graphs upon the coordinate $z_{k} \in \mathbb{D}_{\tilde{\Gamma}}$. We take the foliation $(\mathcal{V}_{t})_{t \in \mathbb{D}}$  of $\mathcal{V}$  by the $u$-curves $\mathcal{V}_{t} = \mathcal{V} \cap \{z_{3} = t\}$. For any $t \in \mathbb{D}$,  $\tilde{\mathcal{V}}_{t} =  f^{2qr}(\mathcal{V}_{t} \cap g_{j}(U_{j})  ) $ is a $u$-curve by invariance of $C^{u}$. Suppose that $\tilde{\mathcal{V}}_{t}$ intersects $\tilde{\Gamma}$ at some point $\tilde{\gamma}(Z)$ such that $\mathrm{pr}_{k}(\tilde{\gamma}(Z))$ belongs to the disk of same center as $\mathbb{D}_{\tilde{\Gamma}}$ and half radius. Since the tangent spaces of $ \tilde{\mathcal{V}}_{t}$ are directed by vectors in $C^{u}$ and since $\mathrm{diam}(\mathrm{pr}_{1}(\tilde{\Gamma})) \le 10^{-10}$, $\tilde{\mathcal{V}}_{t}$ intersects $\tilde{\Gamma}$  in exactly two points by Lemma \ref{rq5} and then it is also the case for $\mathcal{V}_{t}$ and  $\Gamma$. It is clear there exists infinitely many $t \in \mathbb{D}$ satisfying this property. This implies that $\text{pr}_{3}:\Gamma \rightarrow \mathbb{D}$ is a 2-covering. By the Riemann-H\"urwitz formula, it is a 2-covering with only one point of ramification. Moreover this shows that $\mathrm{fold} ( \Gamma  ) \in \text{pr}_{3}(g_{j}(\mathbb{D}^{2} \times \mathbb{D}(\text{fold}(\tilde{\Gamma}), 2 \cdot 10^{-7}))$. Then by Corollary \ref{estit}, we have $\mathrm{fold}(\Gamma') \in \mathbb{D}( \frac{1}{\mu^{2qr}} ( \text{fold}(\tilde{\Gamma})- c_{j} ), \frac{3}{2} \cdot 10^{-6})$.
\end{proof}
   
From now on, we prove lemmas that will show that the image under $g_{j}$ of a 3-folded $(2,3)$-surface having its fold in $ \mathbb{D}(0,\frac{1}{10})$, is concentrated and that it is possible to choose $j \in \{1,2,3,4\}$ so that the new fold is still in $ \mathbb{D}(0,\frac{1}{10})$. The following is an intermediate result to prove Lemma \ref{rp3}:

\begin{lemma} \label{presqueenfin}
Reducing $b_{0}$, $\sigma_{0}$ and $\mathcal{F}$ if necessary, there exists some constant  $\chi'>\chi>1$ and some compact neighborhood $\tilde{\mathbb{D}} $ of $\mathbb{D}$ with  $\mathbb{D} \Subset \tilde{\mathbb{D}} \Subset \mathbb{D}'$ such that for any $f \in \mathcal{F}$, $j \in \{1,2,3,4\}$ and $(1,3)$-quasi plane $\mathcal{V}'$ over $\mathbb{D}^{2}$, the set $f^{2qr} ( \mathcal{V}' )$ contains a  $(1,3)$-quasi plane $\mathcal{V} = \{z_{2} = v(z_{1},z_{3})\}$ over $\tilde{\mathbb{D}} \times \mathbb{D}(0,2/3)$. Moreover $\mathcal{V}$ is foliated by $u$-curves whose tangent vectors are all included in the subcone $ \tilde{C}^{u} = \{v = (v_{1},v_{2},v_{3}) \in \mathbb{C}^{3} : \max(|v_{2}|,|v_{3}|)  \le  \chi'^{-1}  \cdot |v_{1}|\} $  of $C^{u}$.
\end{lemma}

\begin{proof}
Let us fix $j \in \{1,2,3,4\}$. By invariance of $C^{ss}$, the set $f^{2qr} \big( \mathcal{V}' \cap g_{j}(U_{j}) \big)$ contains a holomorphic graph $\mathcal{V} = \{z_{2} = v(z_{1},z_{3})\}$ over $\mathbb{D} \times \mathbb{D}(0,2/3)$ included in $U_{j}$. When $b =\sigma = 0$, it is obvious that the $C^{0}$-norm of $Dv$ is bounded by 1 so it is still the case for every $f \in \mathcal{F}$, reducing  $b_{0}$, $\sigma_{0}$ and $\mathcal{F}$ if necessary. By property 1 of Proposition \ref{marber}, $p_{c}^{qr}$ is univalent on some neighborhood of $h_{1}^{q}(\mathbb{D})$ and on some neighborhood of $h_{2}^{q}(\mathbb{D})$. Then, reducing $b_{0}$, $\sigma_{0}$ and $\mathcal{F}$ if necessary, there exists some compact neighborhood $\tilde{\mathbb{D}}$  of $\mathbb{D}$ independent of $\mathcal{V}'$ such that $\mathbb{D} \Subset \tilde{\mathbb{D}} \Subset \mathbb{D}'$ and $v$ can be extended from $\tilde{\mathbb{D}}   \times \mathbb{D}(0,2/3)$ into $\mathbb{D}$. The extended graph is included in $f^{2qr} \big( \mathcal{V}')$. By property 3 of Proposition \ref{marber} we have $|p'_{c}|>\chi$ on $\mathbb{D}'$. Since $p'_{c}$ is continuous, this implies that $|p'_{c}| \ge\chi''$ for some $\chi''>\chi$ on $\tilde{\mathbb{D}}$ (reducing $\tilde{\mathbb{D}}$ if necessary). We denote $\chi' = \frac{1}{2}(\chi+\chi')$. Then, for a single point $z $ in a compact neighborhood of $(h_{1}^{q}(\mathbb{D})  \cup  h_{2}^{q}(\mathbb{D})) \times \mathbb{D}^{2}$, reducing $b_{0}$, $\sigma_{0}$ and $\mathcal{F}$ if necessary, the differential of $f^{2qr}$ at $z$ sends the closure of $C^{u}$ into $\tilde{C}^{u}$. Since $b_{0}$, $\sigma_{0}$ and $\mathcal{F}$ can be taken locally constant in $z$, by compactness, we can reduce $b_{0}$, $\sigma_{0}$ and $\mathcal{F}$ such that this remains true for every $z$. In particular, $\mathcal{V}$ is a $(1,3)$-quasi plane over $\mathbb{D}^{2}$ foliated by $u$-curves whose tangent vectors are all included in $ \tilde{C}^{u}$. Finally we take the minimal values of $b_{0}$, $\sigma_{0}$ and $\mathcal{F}$ over $j \in \{1,2,3,4\}$. This concludes the proof.  
\end{proof}

Since $\mathbb{D} \times \mathbb{D}(0,1/2) \Subset \tilde{\mathbb{D}} \times \mathbb{D}(0,2/3) $, by the Cauchy inequality, there exists  $\mathcal{C}>0$ such that for any holomorphic map $w$ from $\tilde{\mathbb{D}} \times \mathbb{D}(0,2/3)$ into $\mathbb{D}$, the $C^{0}$-norm of $Dw$ on  $\mathbb{D} \times \mathbb{D}(0,1/2)$ is bounded by $\mathcal{C}$.

\begin{lemma}\label{est}
Reducing $b_{0}$, $\sigma_{0}$ and $\mathcal{F}$ if necessary, for every $f \in \mathcal{F}$, $j \in \{1,2,3,4\}$ and $z \in \mathbb{D}$, we both have $\mathrm{diam} \big( \mathrm{pr}_{2}(U_{j}  \cap \{z_{1} = z \}    ) \big) < 10^{-7} \cdot \mathrm{dist} (U_{j}, \mathbb{D} \times \partial \mathbb{D} \times \mathbb{D})$ and $\mathrm{diam} \big( \mathrm{pr}_{2}(U_{j}  \cap \{z_{1} = z \}    ) \big) < 10^{-1} \cdot \mathcal{C}^{-1} \cdot (\chi'-\chi)$.
\end{lemma}

\begin{proof}
We notice that  for the map $(z_{1},z_{2},z_{3}) \mapsto (p_{c}(z_{1})+bz_{2},z_{1}, \lambda z_{1} + \mu z_{3} + \nu)$, the diameter $\mathrm{diam} \big( \mathrm{pr}_{2}(U_{j}  \cap \{z_{1} = z \}    ) \big)$ tends to 0 uniformly in $j \in \{1,2,3,4\}$ and  $z \in \mathbb{D}$ when $b$ tends to 0 so we can decrease $b_{0}$ so that the two inequalities are satisfied for every $j \in \{1,2,3,4\}$ and  $z \in \mathbb{D}$. Since both inequalities are open, reducing $\sigma_{0}$ and $\mathcal{F}$ if necessary, both remain true for every  $f \in \mathcal{F}$, $j \in \{1,2,3,4\}$ and $z \in \mathbb{D}$.
\end{proof}

\begin{lemma} \label{rp3}
Let $\mathcal{V}'^{0}$ and $\mathcal{V}'^{1}$ be two $(1,3)$-quasi planes over $\mathbb{D}^{2}$ and let $\mathcal{V}^{0}  =f^{2qr} \big( \mathcal{V}'^{0} \cap g_{j}(U_{j}) \big)$ and $\mathcal{V}^{1} =f^{2qr} \big( \mathcal{V}'^{1} \cap g_{j}(U_{j}) \big)$ for some $j \in \{1,2,3,4\}$. Then there exists a holomorphic family $(\mathcal{V}_{t})_{t \in \mathbb{D}(0,10^{6})}$ where $\mathcal{V}^{0} \subset \mathcal{V}_{0}$ and $\mathcal{V}^{1} \subset  \mathcal{V}_{1}$, and for every $t \in \mathbb{D}(0,10^{6})$, $\mathcal{V}_{t}$ is a $(1,3)$-quasi plane over $\mathbb{D} \times \mathbb{D}(0,1/2)$.
\end{lemma}

\begin{proof}
By the proof of Lemma \ref{presqueenfin}, both $\mathcal{V}^{0}  $ and $\mathcal{V}^{1}$ contain $(1,3)$-quasi planes $\mathcal{V}_{0}  $ and $\mathcal{V}_{1}$ over $\mathbb{D} \times \mathbb{D}(0,1/2)$ included in $U_{j}$. The $(1,3)$-quasi plane $\mathcal{V}_{0}$ can be written as a graph $(z_{1},z_{3}) \mapsto v_{0}(z_{1},z_{3})$ over $\mathbb{D} \times \mathbb{D}(0,1/2)$ and the $(1,3)$-quasi plane $\mathcal{V}_{1}$ can be written as a graph $(z_{1},z_{3}) \mapsto v_{1}(z_{1},z_{3})$ over $\mathbb{D} \times \mathbb{D}(0,1/2)$. For every $t \in \mathbb{D}(0,10^{6})$, we denote $v_{t}(z_{1},z_{3}) =  v_{0}(z_{1},z_{3}) + t \cdot ( v_{1}(z_{1},z_{3}) - v_{0}(z_{1},z_{3}) ) $, which defines a graph $\mathcal{V}_{t} $ over $(z_{1},z_{3}) \in \mathbb{D} \times \mathbb{D}(0,1/2)$. By the first inequality of Lemma \ref{est}, for every $t \in \mathbb{D}(0,10^{6})$, $\mathcal{V}_{t} $ does not  intersect $\mathbb{D} \times \partial \mathbb{D} \times \mathbb{D}$. By Lemma \ref{presqueenfin}, $\mathcal{V}_{0}$ is foliated by $u$-curves whose tangent vectors are all included in $ \tilde{C}^{u} $. We notice that $ | v_{1}(z_{1},z_{3}) - v_{0}(z_{1},z_{3})|$ is bounded by $\sup_{z \in \mathbb{D}} \mathrm{diam} \big( \mathrm{pr}_{2}(U_{j}  \cap \{z_{1} = z \}    ) \big) < 10^{-1} \cdot \mathcal{C}^{-1} \cdot (\chi'-\chi)$ by the second inequality of Lemma \ref{est}. Still using  Lemma \ref{presqueenfin}, $(z_{1},z_{3}) \mapsto v_{1}(z_{1},z_{3}) - v_{0}(z_{1},z_{3})$ can be extended on  $\tilde{\mathbb{D}} \times \mathbb{D}(0,2/3)$. Reducing $b_{0}$, $\sigma_{0}$ and $\mathcal{F}$, its image is also included in $\mathbb{D}(0, 10^{-1} \cdot  \mathcal{C}^{-1} \cdot (\chi' - \chi))$. Then by definition of $\mathcal{C}$, the $C^{0}$-norm of $D (v_{1}-v_{0})$ on  $\mathbb{D} \times \mathbb{D}(0,1/2)$ is bounded by $\mathcal{C} \cdot 10^{-1} \cdot \mathcal{C}^{-1}  \cdot (\chi'-\chi)= 10^{-1} \cdot (\chi'- \chi)$. Since $\mathcal{V}^{0}$ is foliated by $u$-curves whose tangent vectors are all included in $ \tilde{C}^{u}$,  $\mathcal{V}^{t}$ is then foliated by curves with tangent vectors in $C^{u}$, that is $u$-curves, for every $t \in \mathbb{D}(0,10^{6})$. Thus  $\mathcal{V}_{t}$ is a $(1,3)$-quasi plane over $\mathbb{D} \times \mathbb{D}(0,1/2)$ for every $t \in \mathbb{D}(0,10^{6})$. 
\end{proof}

\begin{lemma} \label{rp5}
Let $\mathcal{V}_{0}$ and $\mathcal{V}_{1}$ be two $(1,3)$-quasi planes over $\mathbb{D} \times \mathbb{D}(0,1/2)$ and a holomorphic family $(\mathcal{V}_{t})_{t \in \mathbb{D}(0,10^{6})}$  containing  $\mathcal{V}_{0}$ and $\mathcal{V}_{1}$ and such that $\mathcal{V}_{t}$ is a $(1,3)$-quasi plane over $\mathbb{D} \times \mathbb{D}(0,1/2)$ for $t \in \mathbb{D}(0,10^{6})$. Let $\mathcal{W}$ be a 3-folded $(2,3)$-surface such that $\mathrm{Fold}(\mathcal{W}) \subset \mathbb{D}(0,1/10)$. Then $|\mathrm{fold}(\mathcal{V}_{1} \cap \mathcal{W}) -  \mathrm{fold}(\mathcal{V}_{0} \cap \mathcal{W}) |<  10^{-6}$.
\end{lemma}

\begin{proof}
We consider the function $t \mapsto  \text{fold}(\mathcal{V}^{t} \cap \mathcal{W})$ defined on $\mathbb{D}(0,10^{6})$. It is holomorphic by the Implicit Function Theorem and its image is included in $ \mathbb{D}(0,1/10)$ because $\text{Fold}(\mathcal{W}) \subset \mathbb{D}(0,1/10)$. Then by the Cauchy inequality its derivative is smaller than $\frac{1}{10} \cdot 2 \cdot 10^{-6}$ on $\mathbb{D}(0,2)$. Then  $|\mathrm{fold}(\mathcal{V}_{1} \cap \mathcal{W}) -  \mathrm{fold}(\mathcal{V}_{0} \cap \mathcal{W}) |<   1 \cdot \frac{1}{10} \cdot  2 \cdot  10^{-6}  < 10^{-6}$.
\end{proof}
  
The following proposition is important because it says that the image under $g_{j}$ of a 3-folded $(2,3)$-surface is a concentrated folded surface. We use it to prove Corollary \ref{rrrr} and we will also use it in Section 5.

\begin{prop} \label{rp7}
Let $1 \le j \le 4$. Let $\mathcal{W}$ be a 3-folded $(2,3)$-surface such that  we both have $\mathrm{diam}(\mathrm{pr}_{1}(\mathcal{W})) \le 10^{-10}$ and $\mathrm{Fold}(\mathcal{W}) \subset \mathbb{D}(0,1/10)$. Then $g_{j}(\mathcal{W} \cap U_{j})$ is a concentrated 3-folded $(2,3)$-surface satisfying $\mathrm{diam}(\mathrm{Fold}(g_{j}(\mathcal{W}))) \le \frac{1}{2} \cdot  10^{-5}$ for any $f \in \mathcal{F}$.
\end{prop}

\begin{proof} 
Let $\mathcal{V}$ and $\mathcal{V}'$ be two $(1,3)$-quasi planes over $\mathbb{D}^{2}$. Then $f^{2qr}(\mathcal{V} \cap g_{j}(   \mathbb{D}^{3} )  )$ and $f^{2qr}(\mathcal{V}' \cap g_{j}( \mathbb{D}^{3}    ))$ contain $(1,3)$-quasi planes over  $\mathbb{D} \times \mathbb{D}(0,1/2)$ included in $U_{j}$ by  Lemma \ref{presqueenfin}. The sets $\Gamma = f^{2qr}(\mathcal{V} \cap g_{j}(      \mathbb{D}^{3})    ) \cap \mathcal{W}$ and $\Gamma' = f^{2qr}(\mathcal{V}' \cap g_{j}(           \mathbb{D}^{3}   )) \cap \mathcal{W}$ are 3-folded curves. By Lemmas \ref{rp3} and \ref{rp5}, $|\text{fold}(\Gamma') -  \text{fold}(\Gamma) |<10^{-6}$. Lemma \ref{lemb5} implies that $g_{j}( \Gamma)$ is a 3-folded curve with $\text{fold} \big(  g_{j}( \Gamma) \big) \in \mathbb{D}(\frac{1}{\mu^{2qr}} ( \text{fold}(\Gamma) - c_{j}) , (3/2)  \cdot 10^{-6})$ and the analogous for $\Gamma'$. Then $\mathrm{diam}(\text{Fold}(g_{j}(\mathcal{W})) ) \le |\mu^{2qr}|^{-1} \cdot 10^{-6}   + 2 \cdot (3/2) \cdot 10^{-6}<(1/2) \cdot  10^{-5}$ and then $g_{j}(\mathcal{W})$ is a concentrated 3-folded $(2,3)$-surface.
\end{proof}

\begin{corollary}   \label{rrrr}
Let $\mathcal{W}$ be a concentrated 3-folded $(2,3)$-surface such that it holds both $\mathrm{diam}(\mathrm{pr}_{1}(\mathcal{W})) \le 10^{-10}$ and $\mathrm{Fold}(\mathcal{W}) \subset  \mathbb{D}(0,\frac{1}{10})$. Then there exists $j \in \{1,2,3,4\}$ such that $g_{j}(\mathcal{W} \cap U_{j})$ is a concentrated 3-folded $(2,3)$-surface satisfying the inequalities $\mathrm{diam}(\mathrm{pr}_{1}(g_{j}(\mathcal{W} \cap U_{j}) ) \le 10^{-10}$ and $\mathrm{Fold}(g_{j}(\mathcal{W} \cap U_{j})) \subset  \mathbb{D}(0,\frac{1}{10}) $.
\end{corollary}

\begin{proof}
Let $\Gamma$ be a 3-folded curve included in $\mathcal{W}$. By Proposition \ref{rp9}, there exists $j \in \{1,2,3,4\}$ such that  $\frac{1}{\mu^{2qr}} (\mathrm{fold}(\Gamma) -c_{j} ) \in  \mathbb{D}(0,\frac{1}{10}-10^{-5})$. By Lemma \ref{lemb5}, $g_{j}( \Gamma)$ is a 3-folded curve with $\mathrm{fold} (g_{j}(\Gamma)) \in \mathbb{D}(\frac{1}{\mu^{2qr}} (\mathrm{fold}(\Gamma) -c_{j}),(3/2) \cdot 10^{-6})$. By Proposition \ref{rp7}, $g_{j}(\mathcal{W})$ is concentrated with $\mathrm{diam}(\mathrm{Fold}(g_{j}(\mathcal{W}))) \le \frac{1}{2} \cdot  10^{-5}$. Then $\mathrm{Fold}(g_{j}(\mathcal{W} \cap U_{j}))$ is included in  $\mathbb{D}(0 ,\frac{1}{10}-10^{-5}+ (3/2) \cdot 10^{-6}+  \frac{1}{2} \cdot 10^{-5}) \subset  \mathbb{D}(0,\frac{1}{10})$ for $f \in \mathcal{F} $. The inequality $\mathrm{diam}(\mathrm{pr}_{1}(g_{j}(\mathcal{W} \cap U_{j}) ) \le 10^{-10}$ comes from $\mathrm{diam}(\mathrm{pr}_{1}(\mathcal{W})) \le 10^{-10}$ and the forward dilatation of $C^{u}$. The proof is done.
\end{proof}

\begin{proof}[Proof of Proposition \ref{rrrrr}]
By iteration of Corollary \ref{rrrr} there exists a sequence $(j_{n})_{n \ge 1}$ of digits in $\{1,2,3,4\}$ such that the sequence $(\mathcal{W}_{n})_{n \ge 0}$ defined by $\mathcal{W}_{0} = \mathcal{W}$ and $\mathcal{W}_{n+1} = g_{j_{n}} (\mathcal{W}_{n} \cap U_{j_{n}})$ is a sequence of concentrated 3-folded $(2,3)$-surfaces with $\text{Fold}(\mathcal{W}^{n}) \subset  \mathbb{D}(0,\frac{1}{10}) \subset \mathbb{D}$ and $\text{diam}(\text{pr}_{1}(\mathcal{W}^{n}  )) \le 10^{-10}$ for every $n \ge 1$. We define for every $n \ge 1$ $\tilde{\mathcal{W} }_{n}= f^{n}(\mathcal{W}_{n}) \subset \mathcal{W}_{0}$. We have for every $n \ge 1$ the inclusions $\tilde{\mathcal{W} }_{n+1} \subset \tilde{\mathcal{W} }_{n} \subset \tilde{\mathcal{W} }_{0}$. The set $\tilde{\mathcal{W}}_{\infty} = \bigcap_{n \ge 1} \tilde{\mathcal{W} }_{n}$ is non empty since it contains a decreasing sequence of  non empty compact sets. Since $\mathcal{W}_{n}$ is a 3-folded $(2,3)$-surface for every $n \ge 0$, there exists $z_{n} \in \mathcal{W}_{n}$ and a non zero vector $v_{n} \in T_{z_{n}}\mathcal{W}_{n}$ such that $v_{n} \in C^{u}$. We denote for every $n \ge 1$, $\tilde{z}_{n} = f^{n}(z_{n}) \in \tilde{\mathcal{W} }_{n} \subset \tilde{\mathcal{W}}^{0}$ and $\tilde{v}_{n}$ an unitary vector parallel to $D_{z_{n}}f^{n}(v_{n})$. We have $\tilde{v}_{n} \in T_{\tilde{z}_{n}}\tilde{\mathcal{W} }_{n}$. Taking a subsequence if necessary we can suppose $\tilde{z}_{n} \rightarrow \tilde{z}_{\infty} \in \tilde{\mathcal{W}}_{\infty}$ and $\tilde{v}_{n} \rightarrow \tilde{v}_{\infty}$ for some point $\tilde{z}_{\infty}$ and some vector $\tilde{v}_{\infty}$. Clearly $\tilde{z}_{\infty} \in  \mathcal{W}$. By construction the whole forward orbit of $\tilde{\mathcal{W}} _{\infty}$ is in $\mathbb{D}^{3}$. Then $\tilde{z }_{\infty}$ is included in the unstable manifold $W^{u}(\kappa_{f})$ of some point $\kappa_{f}$ of $\mathcal{H}_{f}$ . By construction $\tilde{v}_{\infty} \in T_{\tilde{z}_{\infty}}\mathcal{W}$ and $\tilde{v}_{\infty} \in \bigcap_{n \ge 0} Df^{n}(C^{u}(z_{n} )) = T_{\tilde{z}_{\infty}}W^{u}_{\text{loc}}(\tilde{z}_{\infty})$. Then $\mathcal{W}$ has a point of tangency with the unstable manifold of $\kappa_{f}$.
\end{proof}

\section{Initial heteroclinic tangency}

In this section, we show that:

\begin{enumerate}
\item for every sufficiently small values of $b$ and $\sigma$, we can find $c_{1} = c_{1}(b,\sigma)$ such that $f_{1} = f_{c_{1},b,\sigma}$ has a point of heteroclinic tangency $\tau$ between $W^{s}(\alpha_{f_{1}})$ and $W^{u}(\phi_{f_{1}})$ (where $\phi_{f_{1}}$ is a periodic point in $\mathcal{H}_{f_{1}}$), 
\item we can take iterates of a neighborhood of $\tau$ inside $W^{s}(\alpha_{f_{1}})$ under $f^{-1}_{1}$ in order to create a concentrated 3-folded $(2,3)$-surface inside $W^{s}(\alpha_{f_{1}})$.
\end{enumerate}

\subsection{Initial tangency}

We recall that the disk $\mathcal{C}$ and the integer $m$ were defined in Proposition \ref{marber}.

\begin{prop} \label{orr}
Reducing $\mathcal{C}$ if necessary, there exist $0<b_{1}<b_{0}$, $0< \sigma_{1}<\sigma_{0}$ and an integer $s$ such that for every $0<|b| < b_{1}$ and $0 \le |\sigma| < \sigma_{1}$:
\begin{enumerate}
\item  for every $u$-curve $\mathcal{U}$, $f_{0}^{s+m}(\mathcal{U})$ contains a degree 2 curve over $z_{1}$,
\item  for every holomorphic family of $u$-curves $(\mathcal{U}_{c})_{c \in \mathcal{C}}$, there exists $c_{1} = c_{1}(b,\sigma)$ displaying a quadratic tangency $\tau$ between $W^{s}(\alpha_{f_{1}})$ and $\mathcal{U}_{c_{1}}$ where $f_{1} = f_{c_{1},b,\sigma}$ and $\tau \in \mathcal{U}_{c_{1}}$. Every iterate of $f_{1}^{s+m}(\tau)$ under $f_{1}$ is in $\mathbb{D} \big( \alpha_{c_{0}},10^{-10} \cdot  |w_{m-1}| \cdot (\chi-1)  \big)  \times \mathbb{C}^{2}$ (where the constant $w_{m-1}$ was introduced in Definition \ref{K}) and the mapping $(b,\sigma) \mapsto \tau$ is holomorphic. 
\end{enumerate}
\end{prop}

\begin{proof}
In the following, we are going to reduce several times the bounds $\sigma_{0}$ and $b_{0}$ into bounds $\sigma_{1}$ and $b_{1}$ to satisfy the two items. We begin by taking $\sigma_{1} =  \frac{\sigma_{0}}{2}$ and $b_{1} = \frac{b_{0}}{2}$. Let us take any $u$-curve $\mathcal{U}$ or any holomorphic family of $u$-curves $(\mathcal{U}_{c})_{c \in \mathcal{C}}$. Since $\mathbb{D}'$ intersects the Julia set of $p_{c}$ for every $c \in \mathcal{C}$, reducing $\mathcal{C}$ if necessary, it is possible to find an integer $s$ and a holomorphic map $\beta_{-s} : \mathcal{C} \rightarrow \mathbb{D}$ such that for every $ c \in \mathcal{C}$, we have $p_{c}^{s}(\beta_{-s}(c)) = 0$. In particular, we have $p_{c_{0}}^{s+m}(\beta_{-s}(c_{0})) = \alpha_{c_{0}}$ and the image of $c \rightarrow p_{c}^{s+m}(\beta_{-s}(c)) - \alpha_{c}$ is an open set which contains 0 in its interior. Then reducing sufficiently the bounds $0 \le |\sigma| < \sigma_{1}$ and $0< |b| < b_{1}$ for $\sigma$ and $b$ and $ \mathcal{C}$, we have that for any $c \in \mathcal{C}$, there exists a neighborhood of the point of $\mathcal{U}$ of first coordinate $z_{1} = \beta_{-s}(c_{0})$ inside $\mathcal{U}$ whose image under $f^{s}_{0}$ is of the form $\{(z_{1},u^{2}(z_{1}),u^{3}(z_{1})),z_{1} \in \mathbb{D}(0,\rho) \}$ for some $\rho>0$. Indeed, $\beta_{-s}(c), p_{c}( \beta_{-s}(c)  ), \cdots , p^{s-1}_{c}( \beta_{-s}(c)  )$ are not critical points of $p_{c}$. The image of the curve $\{(z_{1},u^{2}(z_{1}),u^{3}(z_{1})),z_{1} \in \mathbb{D}(0,\rho) \}$ under $f_{0}$ is the curve $\{(p_{c}(z_{1})+b u^{2}(z_{1}) + \sigma  u^{3}(z_{1})  (z_{1}-\alpha_{c_{0}}),z_{1}, \lambda  z_{1}+\mu  u^{3}(z_{1})+ \nu),z_{1} \in \mathbb{D}(0,\rho) \}$. Then it has a point of quadratic tangency with the foliation $z_{1} = C^{st}$ if $b$ and $\sigma$ are sufficiently small. Indeed, reducing $b_{1}$ and $\sigma_{1}$ if necessary, by continuity the derivative of the first coordinate $p_{c}(z_{1})+b  u^{2}(z_{1}) + \sigma u^{3}(z_{1})  (z_{1}-\alpha_{c_{0}})$ vanishes for some value $\overline{z_{1}} \in \mathbb{D}(0,\rho)$. We can iterate this curve $(m-1)$ times under $f_{0}$. Since $p_{c}$ has no other critical point, this will still be a degree 2 curve upon $z_{1}$, reducing $b_{1}$ and $\sigma_{1}$ if necessary. \medskip

 In the case of a holomorphic family of $u$-curves $(\mathcal{U}_{c})_{c \in \mathcal{C}}$, there exists a neighborhood of the point of $\mathcal{U}_{c}$ of first coordinate $z_{1} = \beta_{-s}(c)$ inside $\mathcal{U}_{c}$ which is sent under $f_{0}^{s}$ on a $u$-curve $\{(z_{1},u_{c}^{2}(z_{1}),u_{c}^{3}(z_{1})),z_{1} \in \mathbb{D}(0,\rho) \}$, but using the Cauchy inequality and reducing $\sigma_{1}$ and $b_{1}$ if necessary, $(u^{2}_{c})'$ and $(u^{3}_{c})'$ are uniformly (relatively to $c$) bounded and the conclusion is the same. In particular, this proves the first item of the result. \medskip
 
We call $\mathrm{Tan}_{f_{0}}$ the first coordinate of the point of vertical tangency. The image of the map defined on $\mathcal{C}$ which sends $c$ to $p_{c}^{s+m}(\beta_{-s}(c)) - \alpha_{c}$ is an open set which contains 0 in its interior. Let us denote by $2  l$ the distance of 0 to the image of $\partial \mathcal{C}$. Reducing $b_{1}$ and $\sigma_{1}$ another time if necessary, by continuity $\{  \mathrm{Tan}_{f_{0}} - \alpha_{c} , c \in \partial \mathcal{C} \}$ is a curve in the plane, with 0 is in a bounded connected component of its complement and at distance at least $l$ from $\{  \mathrm{Tan}_{f_{0}} - \alpha_{c} , c \in \partial \mathcal{C} \}$. 
 
\begin{notation}
Let $\kappa_{f} \in \mathcal{H}_{f}$ for $f \in \mathcal{F}$. We denote by $W^{s}_{\mathrm{loc}}(\kappa_{f})$ the connected component of $W^{s}(\kappa_{f}) \cap \mathbb{D}^{3}$ which contains $\kappa_{f}$ and by $W^{u}_{\mathrm{loc}}(\kappa_{f})$ the connected component of $W^{u}(\kappa_{f}) \cap \mathbb{D}^{3}$ which contains $\kappa_{f}$.
\end{notation} 
 
\begin{lemma} \label{utileresult}
Reducing $\mathcal{C}$, $b_{1}$ and $\sigma_{1}$ if necessary, $ W^{s}_{\mathrm{loc}}(\alpha_{f_{0}})$ is a graph over $(z_{2},z_{3}) \in  \mathbb{D}^{2}$ both included in $\mathbb{D} \big( \alpha_{c},l/2  \big)  \times \mathbb{D}^{2}$ and in $\mathbb{D} \big( \alpha_{c_{0}},10^{-10} \cdot  |w_{m-1}| \cdot (\chi-1)  \big)  \times \mathbb{D}^{2}$.
\end{lemma}

\begin{proof}
For $\sigma = 0$, $ W^{s}(\alpha_{f_{0}})$ is the product of $ W^{s}(\alpha_{H})$ by the $z_{3}$ axis where $H$ is the H\'enon map $H : (z_{1},z_{2}) \mapsto (p_{c}(z_{1})+bz_{2},z_{1})$. Moreover, for every $\epsilon>0$, we can reduce $b_{1}$ such that if $ |b|< b_{1}$, then the cone $C^{ss,\epsilon}$ centered at $e_{2}$ of opening $\epsilon/2$ is $H^{-1}$-invariant in some neighborhood of $\{\alpha_{c}\} \times \mathbb{D}$. Then $W^{s}_{\text{loc}}(\alpha_{H})$ is a $ss$-curve included in $\mathbb{D}(\alpha_{c},\epsilon) \times \mathbb{D}$ since all its tangent vectors lie in $C^{ss,\epsilon}$. Then for $\sigma = 0$, the skew-product structure implies that $W_{\mathrm{loc}}^{s}(\alpha_{f_{0}})$ is a $(2,3)$-surface included in $\mathbb{D}(\alpha_{c},\epsilon) \times \mathbb{D}^{2}$ which is the product of $W_{\mathrm{loc}}^{s}(\alpha_{H})$ by the $z_{3}$ axis. Then, by continuity, it is possible to reduce $\sigma_{1}$ such that for every $f_{0}$ with $ 0< |b|<b_{1}$ and $0 \le |\sigma|< \sigma_{1}$, $ W_{\mathrm{loc}}^{s}(\alpha_{f_{0}})$ is a  graph over $(z_{2},z_{3}) \in \mathbb{D}^{2}$ included in $\mathbb{D}(\alpha_{c},\epsilon) \times \mathbb{D}^{2}$. We take $\epsilon = l/2$ to prove the first inclusion. We get the second one by reducing both $\mathcal{C}$ and $\epsilon$.
\end{proof}

We recall that 0 is in a bounded connected component of the complement of $\{  \mathrm{Tan}_{f_{0}} - \alpha_{c} , c \in \partial \mathcal{C} \}$ and at distance at least $l$ from $\{  \mathrm{Tan}_{f_{0}} - \alpha_{c} , c \in \partial \mathcal{C} \}$. In particular, there is a parameter $\overline{c} \in \mathcal{C}$ such that $\text{Tan}_{f_{0}}$ belongs to $\mathbb{D}(\alpha_{\overline{c}},l/2)$. Up to replacing $s+m$ by $s+m+1$, we can suppose that the point of vertical tangency is in $\mathbb{D}^{3}$. For the parameter $\overline{c}$, $f^{s+m}_{0}(\mathcal{U}_{\overline{c}}) \cap (\mathbb{D}(\alpha_{\overline{c}},l/2) \times \mathbb{D}^{2}) $ is not  the union of two graphs upon $z_{1} \in \mathbb{D}(\alpha_{\overline{c}}, l/2)$ and for $c \in \partial \mathcal{C} $, $f^{s+m}_{0}(\mathcal{U}_{c})$ is the union of two graphs over $z_{1} \in \mathbb{D}(\alpha_{c},l/2)$. According to Lemma \ref{utileresult}, $ W_{\mathrm{loc}}^{s}(\alpha_{f_{0}})$ is a graph over $(z_{2},z_{3}) \in  \mathbb{D}^{2}$ included in $\mathbb{D}(\alpha_{c},l/2) \times \mathbb{D}^{2}$. The following lemma is the analogous in dimension 3 of Proposition 8.1 of \cite{dl}. The proof is essentially the same and relies on the continuity of the intersection index of properly intersecting analytic sets of complementary dimensions.
 
\begin{lemma}
Let $(\Gamma_{c})_{c \in \mathcal{C}}$ be a holomorphic family of curves of degree 2 over the first coordinate. We assume that:
\begin{enumerate} 
\item there exists a compact subset $\mathcal{C}' \subset \mathcal{C}$ such that if $c \in  \mathcal{C} \backslash \mathcal{C}'$, $\Gamma_{c}$ is the union of 2 graphs over $z_{1} \in \mathbb{D}(\alpha_{c},l/2 )$, 
\item there exists $\overline{c} \in  \mathcal{C}$  such that $\Gamma_{\overline{c}}$ is not the union of 2 graphs over $z_{1} \in \mathbb{D}(\alpha_{\overline{c}},l/2 )$.
\end{enumerate}
Then, if $(\mathcal{V}_{c})_{c \in  \mathcal{C}}$ is any holomorphic family of graphs over $(z_{2},z_{3}) \in  \mathbb{D}^{2}$ contained in $ \mathbb{D}(\alpha_{c},l/2)   \times  \mathbb{D}^{2}$, there exists $c_{1} \in  \mathcal{C}$ such that $\Gamma_{c_{1}}$ and $\mathcal{V}_{c_{1}}$ admit a point of tangency.
\end{lemma} 

We apply the previous lemma, taking the family $(\Gamma_{c})_{c \in \mathcal{C} }=  (f_{0}^{s+m}(\mathcal{U}_{c}) \cap \mathbb{D}^{3} )_{c \in \mathcal{C} }$ as curves and the family of stable manifolds $W_{\mathrm{loc}}^{s}(\alpha_{f_{0}})    $ as graphs over $(z_{2},z_{3}) \in  \mathbb{D}^{2}$. We can conclude that there is a parameter $c_{1}$ such that there exists a quadratic tangency $f^{s+m}_{1}(\tau)$ between $W_{\text{loc}}^{s}(\alpha_{f_{1}})$ and $f_{1}^{s+m}(\mathcal{U}_{c_{1}})$ where $\tau \in \mathcal{U}_{c_{1}}$ for $f_{1} = f_{c_{1},b,\sigma}$. Then $\tau$ is a point of tangency between $W^{s}(\alpha_{f_{1}})$ and $\mathcal{U}_{c_{1}}$. Since $ W^{s}_{\text{loc}}(\alpha_{f_{1}})$ is included in $\mathbb{D} \big( \alpha_{c_{0}},10^{-10} \cdot  |w_{m-1}| \cdot (\chi-1)  \big)  \times \mathbb{D}^{2}$, all the iterates of $f^{s+m}_{1}(\tau)$ under $f_{1}$ are in $\mathbb{D} \big( \alpha_{c_{0}},10^{-10} \cdot  |w_{m-1}| \cdot (\chi-1)  \big)  \times \mathbb{D}^{2}$. The map $(b,\sigma) \mapsto c_{1}$ is holomorphic  by the Implicit Function Theorem and then $(b,\sigma) \mapsto \tau$ is holomorphic too. This shows item 2 and ends the proof of Proposition \ref{orr}.
\end{proof}

In the next proposition, we show that the tangencies created in the previous result are generically unfolded. Beware that in the next result, the map $f_{1}$ is associated to a family of $u$-curves $(\mathcal{U}_{c})_{c \in \mathcal{C}}$ which can be distinct from the family $(\mathcal{U}'_{f})_{f \in \mathcal{F}'}$. 

\begin{df}
Let $(U_{f})_{f \in \mathcal{F}}$ be a holomorphic family of $u$-curves and $(S_{f})_{f \in \mathcal{F}}$ be a holomorphic family of $s$-surfaces. We suppose that for $f \in \mathcal{F}$, there is a point of quadratic tangency between $U_{f}$ and $S_{f}$. We say that this tangency is generically unfolded if there exists a one-dimensional holomorphic family $(f^{t})_{t \in \mathbb{D}}$ of polynomial automorphisms and a holomorphic family of local biholomorphisms $(\Psi_{t})_{t \in \mathbb{D}}$ with:
\begin{enumerate}
\item $f^{0} = f$ and $f^{t} \in \mathcal{F} $ for every $t \in \mathbb{D}$,
\item $\Psi_{t}(S_{f^{t}})$ is a vertical plane $\{z_{1} = C^{st}\}$ where $C^{st}$ does not depend on $t \in \mathbb{D}$,
\item if we denote by $\mathrm{tan}_{t}$ the first coordinate of the point of tangency of $\Psi_{t}(U_{f^{t}})$ with $\{z_{1} = C^{st}\}$ (there exists exactly one such point because the tangency is quadratic), then $|\frac{ d \mathrm{tan}_{t}}{d t}|$ is uniformly bounded from below by a strictly positive constant for $t \in \mathbb{D}$.
\end{enumerate}
\end{df}

\begin{prop} \label{fd} 
Reducing $\mathcal{C}$, $b_{1}$ and $\sigma_{1}$ if necessary, for every $0<|b| < b_{1}$ and $0 \le |\sigma| < \sigma_{1}$, for every holomorphic family of $u$-curves $(\mathcal{U}_{c})_{c \in \mathcal{C}}$, there exists a neighborhood $\mathcal{F}' \Subset \mathcal{F}$ of the map $f_{1}=f_{c_{1}(b,\sigma),b,\sigma}$ defined in item 2 of Proposition \ref{orr} such that: for every holomorphic family of $u$-curves $(\mathcal{U}'_{f})_{f \in \mathcal{F}'}$, if $f \in \mathcal{F}'$ has a point of tangency $\tau'$ between $W^{s}(\alpha_{f})$ and $\mathcal{U}'_{f}$ such that $f^{s+m}(\tau') \in W^{s}_{\text{loc}}(\alpha_{f})$, then the tangency $\tau'$ is generically unfolded. In particular, the tangency $\tau$ obtained in item 2 of Proposition \ref{orr} is generically unfolded. 
\end{prop}
 
\begin{proof}
For every $f \in \mathcal{F}$, after a holomorphic change of coordinates $\Psi_{f}$, $W^{s}_{\text{loc}}(\alpha_{f})$ is the plane $\{z_{1} = \alpha_{c_{0}}\}$. We first show the result for the map $f_{1}$.  We work in the one-dimensional family $(f_{0})_{c \in \mathcal{C}} = (f_{c,b,\sigma})_{c \in \mathcal{C}}$ where $b$ and $\sigma$ are fixed and we take a  holomorphic family of $u$-curves $(\mathcal{U}'_{f})_{f \in \mathcal{F}'}$. It is a consequence of the proof of Lemma \ref{utileresult} that $\Psi_{f_{0}}$ tends to $\text{Id}$ when $b$ and $\sigma$ tend to 0. When $b$ and $\sigma$ tend to 0, the curve $\Psi_{f_{0}} \circ f^{s+m}_{0}(\mathcal{U}'_{f_{0}})$ tends also to a curve of degree 2 over $z_{1}$ of the form $\{(p_{c}^{s+m}(z_{1}),p_{c}^{s+m-1}(z_{1}),v(z_{1}) ) : z_{1} \in \mathbb{D}(0,\rho)\}$ where $v$ is holomorphic. We call $\mathrm{tan}_{c}$ the first coordinate of the point of vertical tangency of $\Psi_{f_{0}} \circ f^{s+m}_{0}(\mathcal{U}'_{f_{0}})$. We have $p_{c_{0}}^{m}(0) = \alpha_{c_{0}}$, $(p_{c_{0}}^{m})'(0)=0$, $(p_{c_{0}}^{m})''(0) \neq 0$ and at $c= c_{0} $, we have  $ \frac{d  }{dc} \big( p_{c}^{m}(0) - \alpha_{c} \big) \neq 0$ (see Proposition \ref{marber}). Moreover, $c_{1}$ tends to $c_{0}$ when $b$ and $\sigma$ tend to 0. Then for sufficiently small $b$ and $\sigma$, we have $\frac{d \mathrm{tan}_{c}} {dc} \neq 0$ for every $c \in \mathcal{C}$ (reducing $\mathcal{C}$ if necessary).  Since the estimates on the derivatives of the $u$-curves $(\mathcal{U}'_{f})_{f \in \mathcal{F}}$ are uniform, we can reduce uniformly $b_{1}$ and $\sigma_{1}$ so that this inequality is true no matter the choice of $(\mathcal{U}'_{f})_{f \in \mathcal{F}}$. Then, by continuity, there exists some new neighborhood $\mathcal{F}' \Subset \mathcal{F}$ of $f_{1}$ such that every $f \in \mathcal{F}'$ belongs to a one-dimensional family $(f^{t})_{t \in \mathbb{D}}$ such that: $f^{0} = f$, $f^{t} \in \mathcal{F}'$ and $\frac{d \mathrm{tan}_{c}} {dt} \neq 0$ ($t \in \mathbb{D}$). This implies in particular that if $f \in \mathcal{F}'$ has a point of tangency $\tau'$ between $W^{s}(\alpha_{f})$ and $\mathcal{U}'_{f}$ such that $f^{s+m}(\tau') \in W^{s}_{\text{loc}}(\alpha_{f})$, then the tangency $f^{s+m}(\tau')$ is generically unfolded and then also $\tau'$. The proof is over.
\end{proof}

\subsection{A transversality result}

From now on, we construct from this initial heteroclinic tangency a 3-folded $(2,3)$-surface inside a stable manifold with its fold inside $ \mathbb{D}(0,\frac{1}{10})$. In this Subsection, we prove that $W^{s}(\alpha_{f_{1}})$ has some special geometry. We start by choosing a periodic point $\phi_{f_{1}} \in \mathcal{H}_{f_{1}}$. The main point is that its third coordinate is in $\mathbb{D}(0,\frac{1}{10}-10^{-4})$, which will be used later in the proof of Proposition \ref{fm}.

\begin{lemma} \label{defphi}
For every $c \in \mathcal{C}$, $0<|b|<b_{1}$, $0 \le |\sigma|< \sigma_{1}$, there exists a periodic point $\phi_{f_{0}}$ inside $\mathcal{H}_{f_{0}}$ such that $\mathrm{pr}_{3}(\phi_{f_{0}})  \in  \mathbb{D}(0,\frac{1}{10}-10^{-4})$.
\end{lemma}
 
\begin{proof} 
Let us take $\omega \in h^{q}_{1}(\mathbb{D}')$. According to Proposition \ref{blender}, the $ss$-curve $\{\omega\} \times \mathbb{D} \times \{0\}$ intersects the unstable set of a point $\phi_{f_{0}}$ of the horseshoe $\mathcal{H}_{f_{0}}$. Moving slightly  $\{\omega\} \times \mathbb{D} \times \{0\}$ by translation in the third direction (let us say by no more than $\frac{1}{20}$), we can suppose that $\phi_{f_{0}}$ belongs to the preimage $g_{1}(U_{1})$. By density of periodic points in $\mathcal{H}_{f_{0}}$ we can moreover suppose that $\phi_{f_{0}}$ is periodic. We notice that $\omega \in  h^{q}_{1}(\mathbb{D}')$ and $\mathrm{pr}_{1} (  g_{1}(U_{1})  ) \subset h^{q}_{1}(\mathbb{D}')$. Moreover we have $\mathrm{diam}(h^{q}_{1}(\mathbb{D})) < 10^{-11}$ by property 2 of Property \ref{fo}. Since $W^{u}_{\mathrm{loc}}(\phi_{f_{0}})$ is a $u$-curve, this implies that $\mathrm{pr}_{3}(\phi_{f_{0}})$ belongs to $\mathbb{D}(0,  \frac{1}{20} + 10^{-11} \cdot \chi^{-1}  ) \subset \mathbb{D}(0,\frac{1}{10}-10^{-4})$. 
\end{proof}

\begin{notation} 
For every $0<|b|<b_{1}$, $0 \le |\sigma|< \sigma_{1}$, we fix such a periodic point $\phi_{f_{1}}$ for the map $f_{1}$ such that $\mathrm{pr}_{3}(\phi_{f_{1}})  \in  \mathbb{D}(0,\frac{1}{10}-10^{-4})$. Reducing $\mathcal{F}'$ if necessary, we have  $\mathrm{pr}_{3}(\phi_{f})  \in  \mathbb{D}(0,\frac{1}{10}-10^{-4})$ for $f \in \mathcal{F}'$. 
\end{notation}

We now choose the parameter $b$ in function of the parameter $\sigma$: 

\begin{notation} 
In the following, we will take $b= b(\sigma) = \sigma^{2}$. For technical reasons, we reduce $\sigma_{1}$ such that  $\sigma_{1}< \sqrt{b_{1}}$. 
\end{notation} 

In particular, we can use Proposition \ref{orr} for $0 \le |\sigma|< \sigma_{1}$: there exists a map $\sigma \mapsto c_{1}(\sigma)$ such that there is a heteroclinic tangency $\tau$ between the local unstable manifold $\mathcal{U}_{f_{1}}$ of $\phi_{f_{1}}$ and the stable manifold $W^{s}(\alpha_{f_{1}})$ of $\alpha_{f_{1}}$. Moreover, $\sigma \mapsto \tau$ is holomorphic. We have $D(f_{1}^{s})_{\tau}  \cdot  (0,0,1) =  \mu^{s} (\xi_{1}(\sigma),\xi_{2}(\sigma),1)$ with $\xi_{1}(\sigma) = O(\sigma)$ and $\xi_{2}(\sigma) = O(\sigma)$. The maps $\sigma \mapsto \text{pr}_{1}(f^{s}_{1}(\tau)) - \beta_{0}$, $\cdots$, $\sigma \mapsto \text{pr}_{1}(f^{s+m-1}_{1}(\tau)) - \beta_{m-1}$ are holomorphic (remind that the constants $\beta_{i}$ were defined in \ref{bet}). Since they vanish when $\sigma = 0$, they are also $O(\sigma)$. The map $\sigma \mapsto b(\sigma) $ is $O(\sigma^{2})$ because $b(\sigma) = \sigma^{2}$. Then the differentials at $f^{s}_{1}(\tau), \cdots , f^{s+m-1}_{1}(\tau)$ verify the conditions of Corollary \ref{reffin} so from Corollary \ref{reffin} we immediately get the following (remind that $w_{m-1}$ was defined in Definition \ref{K}):

\begin{lemma} \label{lemor}
We have $D(f^{s+m}_{1})_{\tau}  \cdot  ( 0, 0, 1) =  ( \zeta_{1} (\sigma) ,  \zeta_{2}(\sigma)  , \zeta_{3}(\sigma)  )$ where $\zeta_{1}$, $\zeta_{2}$, $\zeta_{3}$ are holomorphic functions, $\zeta_{1}(\sigma) = w_{m-1} \cdot \sigma + O(\sigma^{2})$, $\zeta_{2}(0) = 0+O(\sigma)$ and $\zeta_{3}(\sigma) = \mu^{m}+O(\sigma)$. 
\end{lemma}

\begin{notation} \label{noter}
We reduce $\sigma_{1}$ so that for $0 \le |\sigma|< \sigma_{1}$ it holds $|\zeta_{1}(\sigma)-w_{m-1} \sigma| < \frac{1}{10}  |w_{m-1}|  | \sigma|$, $|\zeta_{2}(\sigma)| <1$ and $|\zeta_{3}(\sigma)-\mu^{m}| <\frac{1}{10}(1-|\mu^{m}|)$. Still reducing $\sigma_{1}$, we have $\sigma_{1}<  10^{-10} \cdot |w_{m-1}| \cdot (\chi-1)$. From now on, we fix $\sigma = \frac{\sigma_{1}}{2}$, $b = b(\sigma) = \sigma^{2}$, $c = c_{1}(\sigma)$ and the associated map $f_{1} = f_{c_{1}(\sigma),b(\sigma),\sigma}$. 
\end{notation}

\begin{lemma}
 \label{or7} For every vector $v^{0}$ such that $|v^{0}_{1} - w_{m-1} \sigma| < \frac{1}{10}  |w_{m-1}|  | \sigma|$, $|v^{0}_{2}  | <1 $  and $|v^{0}_{3} - \mu^{m} | < \frac{1}{10}(1- |\mu|^{m})$, at any point of $W^{s}_{\mathrm{loc}}(\alpha_{f_{1}})$, $v^{0}$ is transverse to $W^{s}_{\mathrm{loc}}(\alpha_{f_{1}})$ for small $\sigma$.
\end{lemma}

\begin{proof}
Let $\psi^{0}$ be a point of $W^{s}_{\text{loc}} (\alpha_{f_{1}})$ and let us consider the sequence of points defined by $\psi^{n} = f^{n}_{1}(\psi^{0})$. According to Lemma \ref{utileresult}, for every $n \ge 0$, $\psi^{n}$ is in $\mathbb{D}(\alpha_{c_{0}}, 10^{-10} |w_{m-1}|(\chi-1)) \times \mathbb{C}^{2}$. Then for every $n \ge 1$, the differential at $\psi^{n}$ of $f_{1}$ is of the form: 
$$I_{n} =  \begin{pmatrix}
   m_{n} & b(\sigma) & \sigma  (z_{1}- \alpha_{c_{0}})\\
   1 & 0& 0 \\
   \lambda & 0 & \mu
\end{pmatrix} \, , $$
where $b(\sigma )= \sigma^{2}$ and $|m_{n}| \ge\chi>1$ for every $n \ge 1$. Since $\sigma_{1}<10^{-10}  |w_{m-1}| (\chi-1)$ we have $|b|  <10^{-10}  |w_{m-1}|  | \sigma|  (\chi-1)$. We denote $v^{n} = (I_{n} \cdot ... \cdot I_{1})(v^{0})$ for every $n \ge 1$. Let us show that there exists an integer $i$ such that $|v^{i}_{1} | \ge  |v^{i}_{3}|$. Let us suppose this is false and we show by induction the following properties for $n \ge 0$: \begin{enumerate}
\item $1 \ge |v^{n+1}_{1}| > \frac{\chi+1}{2} |v^{n}_{1}| \ge  \frac{1}{100}   |w_{m-1}  |  |\sigma  | $ ,
\item $ |v^{n}_{2}| \le1$,
\item $ |v^{n+1}_{3}| \le |v^{n}_{3}| \le 1$.
\end{enumerate}  
For $n = 0$, item 2 is satisfied since  $|v^{0}_{2}  | <1 $. Moreover $ |v^{0}_{1}| <  |v^{0}_{3}|$ by hypothesis. We have $|v^{1}_{3}| \le   |\lambda|  |v^{0}_{1}|+|\mu|  |v^{0}_{3}|<  |\lambda| |v^{0}_{3}|+|\mu|  |v^{0}_{3}|< |v^{0}_{3}| \le 1$. Then item 3 is true. We have $ v^{1}_{1} = m_{1}  v^{0}_{1} + b  v^{0}_{2}+  \sigma ( \psi^{0}_{1} - \alpha_{c_{0}})v_{3}^{0}$. We notice that $|v^{0}_{1}| \ge  \frac{1}{100}   |w_{m-1}  |  |\sigma  | $, $m_{1} \ge \chi$, $|b  v^{0}_{2}| \le  |b|<10^{-10}  |w_{m-1}| | \sigma| (\chi-1)$ and $|   \sigma ( \psi^{0}_{1} - \alpha_{c_{0}}) v^{0}_{3} | < 10^{-10}  |w_{m-1}|  | \sigma|  (\chi-1)$. Then we have $|v^{1}_{1}| > \frac{\chi+1}{2} |v^{0}_{1}| \ge  \frac{1}{100}   |w_{m-1}  |  |\sigma  | $. We also have $|v^{1}_{1}| < |v_{3}^{1}| \le 1$ by hypothesis.  Then item 1 is true and the induction step is true for $n = 0$. \medskip

Suppose  the induction step true for some $n \ge 0$. Then  $|v^{n+1}_{2}| = |v^{n}_{1}| \le 1$ and item 2 is true.  Moreover $|v^{n+2}_{3}| \le  |\lambda| |v^{n+1}_{1}|+|\mu|  |v^{n+1}_{3}|<  |\lambda|  |v^{n+1}_{3}|+|\mu|  |v^{n+1}_{3}|< |v^{n+1}_{3}| \le 1$, which shows item 3. We have $ v^{n+2}_{1} = m_{n+1} v^{n+1}_{1} + b v^{n+1}_{2}+  \sigma( \psi^{n+1}_{1} - \alpha_{c_{0}})v^{n+1}_{3}$. Since we have $|m_{n+1}|>\chi$,  $|v^{n+1}_{1}| \ge  \frac{1}{100}   |w_{m-1}  |  |\sigma  | $, $|b v^{n+1}_{2}| \le |b| <10^{-10} |w_{m-1}|  | \sigma|  (\chi-1)$ and also $| \sigma  ( \psi^{n+1}_{1} - \alpha_{c_{0}})v^{n+1}_{3}  | < 10^{-10}  |w_{m-1}|  | \sigma|  (\chi-1)$, we have $|v^{n+2}_{1}| >\frac{\chi+1}{2} |v^{n+1}_{1}| \ge  \frac{1}{100}   |w_{m-1}  |  |\sigma  | $. We have $ |v^{n+2}_{1}| <  |v^{n+2}_{3}| \le 1$ by hypothesis and then item 1 is true.  \medskip

Then for every $n \ge 0$, we have $|v^{n+1}_{1}| \cdot |v^{n+1}_{3}|^{-1} \ge (\chi+1)/2 \cdot |v^{n}_{1}| \cdot |v^{n}_{3}|^{-1}$. This implies a contradiction: there must exist some integer $i$ such that $|v^{i}_{1} | \ge  |v^{i}_{3}|$. Since we also have $|v^{i}_{2}| = |v^{i-1}_{1}| \le |v^{i}_{1}|$, this implies that $v^{i}$ is transverse to $W^{s}_{\text{loc}}(\alpha_{f_{1}})$ and then $v^{0}$ itself is transverse to $W^{s}_{\text{loc}}(\alpha_{f_{1}})$.
\end{proof}

\begin{prop} \label{orienter}
The vector $(0,0,1)$ is transverse to $W^{s}(\alpha_{f_{1}})$ at $\tau$.
\end{prop}

\begin{proof}
From Lemmas \ref{lemor} and \ref{or7}, we know that the image under $f^{s+m}_{1}$ of a neighborhood of $\tau$ in $W^{s}(\alpha_{f_{1}})$ is transverse to $D(f^{s+m}_{1})_{\tau}  \cdot  ( 0, 0, 1) $. Applying $f^{-s-m}_{1}$, this immediately implies the result. 
\end{proof}

\subsection{Orientation of the fold of $W^{s}(\alpha_{f})$}

 In this subsection, we take iterates of the initial tangency $\tau$ under $f^{-1}_{1}$ in order to create a concentrated 3-folded $(2,3)$-surface inside a stable manifold. We have that $W_{\mathrm{loc}}^{u}(\phi_{f_{1}})$ is a graph $\{(z_{1}, u^{2}(z_{1}),u^{3}(z_{1}) : z_{1} \in \mathbb{D}\}$ over $z_{1} \in \mathbb{D}$ (the periodic point $\phi_{f_{1}}$ was defined in Lemma \ref{defphi}). We consider the biholomorphism $\Psi$ defined by $\Psi(z_{1},z_{2},z_{3}) = (z_{1},z_{2}-u^{2}(z_{1}),z_{3}-u^{3}(z_{1}))$ which sends $W_{\text{loc}}^{u}(\phi_{f_{1}})$ onto the $z_{1}$ axis. We denote $\mathcal{W}_{f}$ the connected component of $ W^{s}(\alpha_{f}) \cap \mathbb{D}^{3}$ which contains $\tau$ for $f \in \mathcal{F}'$. We denote by $\text{Tan}'$ the subset of $\Psi(\mathcal{W}_{f})$ where $\Psi(\mathcal{W}_{f})$ is tangent to some line $\{z_{2} = C^{st}, z_{3} = C^{st}\}$ and by $\text{Tan} = \Psi^{-1}(\text{Tan}')$. When $f = f_{1}$, we denote them by $\text{Tan}^{0}$ and $(\text{Tan}')^{0}$. Let $e^{u}(\tau)$ be a tangent vector of $W_{\text{loc}}^{u}(\phi_{f_{1}})$ at $\tau$.

\begin{lemma} \label{hoppe3}
The curve $\mathrm{Tan}^{0}$ is a complex curve regular at $\tau$ of tangent vector $v_{\text{tan}}$ at $\tau$ such that $v_{\text{tan}} \notin \mathbb{C} \cdot (0,0,1)$ and $v_{\text{tan}} \notin \mathbb{C} \cdot e^{u}(\tau)$. 
\end{lemma} 

\begin{proof} 
We are working in the projectivized tangent bundle $\mathbb{P}T\mathbb{C}^{3} \simeq \mathbb{C}^{3} \times \mathbb{P}^{2}(\mathbb{C})$ of $\mathbb{C}^{3}$ which is of dimension 5. The lift $\hat{\mathcal{W}_{f_{1}}}$ of $\mathcal{W}_{f_{1}}$ to $\mathbb{P}T\mathbb{C}^{3}$ is a complex submanifold of dimension 3. The lift of every complex curve $C_{x,y} = \Psi^{-1}(\mathbb{D} \times \{(x,y) \}) $ to $\mathbb{P}T\mathbb{C}^{3}$ is a complex curve $\hat{C}_{x,y}$. Then $\bigcup_{(x,y) \in \mathbb{D}^{2}} \hat{C}_{x,y}$ is a complex submanifold of dimension 3. Moreover, according to Proposition \ref{orr}, $\mathcal{W}_{f_{1}}$ has a point of quadratic tangency with $C_{0,0}$ and then $\hat{\mathcal{W}_{f_{1}}}$ is transverse to $\bigcup_{(x,y) \in \mathbb{D}^{2}} \hat{C}_{x,y}$. Then $\hat{\text{Tan}^{0}} = \hat{\mathcal{W}_{f_{1}}} \cap \big( \bigcup_{(x,y) \in \mathbb{D}^{2}} \hat{C}_{x,y} \big)$ is a regular complex curve. Then its projection $\text{Tan}^{0}$ is a regular complex curve of tangent vector $v_{\text{tan}}$ at $\tau$. By Proposition \ref{orienter}, $v_{\text{tan}} \notin \mathbb{C} \cdot (0,0,1)$. By construction, $e^{u}(\tau)$ is a tangent vector of $W_{\text{loc}}^{u}(\phi_{f_{1}})$ at $\tau$ so $v_{\text{tan}} \notin \mathbb{C} \cdot e^{u}(\tau)$. The proof is done.
\end{proof}

Since the intersection between $\hat{\mathcal{W}_{f_{1}}}$ and $\bigcup_{(x,y) \in \mathbb{D}^{2}} \hat{C}_{x,y}$ was transverse in the latter proof, by perturbation of the previous result, we get:

\begin{corollary} \label{hoppe5} 
For every $\epsilon>0$, reducing $\mathcal{F}'$ if necessary, there exists $\eta>0$ such that: for every $f \in \mathcal{F}'$, for every  $C^{2}$-foliation $(\mathcal{V}_{x,y})_{(x,y) }$ of some neighborhood of $\tau$ by two-dimensional real manifolds $\mathcal{V}_{x,y}$ such that every $\mathcal{V}_{x,y}$ is $\eta$-close to $W_{\mathrm{loc}}^{u}(\phi_{f})$, the set $\mathrm{Tan}$ of points of $\mathcal{W}_{f}$ where $\mathcal{W}_{f}$ is tangent to some $\mathcal{V}_{x,y}$ is a regular two-dimensional real manifold which has  its direction $\epsilon$-close to $\mathbb{C} \cdot v_{\text{tan}}$ at each point if it is non empty. \end{corollary}

The following result is a technical geometric lemma.

\begin{lemma} \label{eps}
There exist $\epsilon, \rho, \rho_{1},\rho_{2},\rho_{3}>0$ such that for every $t>0$, we have the following property: for every regular two-dimensional real manifold $\Gamma$ going through $\tau+ (\mathbb{D}(0,t\rho))^{3}$, if $\Gamma$ has its direction $\epsilon$-close to $\mathbb{C} \cdot v_{\text{tan}}$ at each point, then $\Gamma$ is horizontal relatively to $\pi_{2}$ in $ \tau + \mathbb{D}(0,t\rho_{1}) e^{u}(\tau)+ \mathbb{D}(0,t\rho_{2}) (0,1,0)+ \mathbb{D}(0,t\rho_{3}) (0,0,1) $.
\end{lemma}

\begin{proof} 
According to Lemma \ref{hoppe3}, $v_{\text{tan}} \notin \mathbb{C} \cdot (0,0,1)$ and $v_{\text{tan}} \notin \mathbb{C} \cdot e^{u}(\tau)$, which easily implies the result. 
\end{proof}

\begin{notation} 
We fix the value of $\epsilon$ given by Lemma \ref{eps} and the associated value of $\eta$ given by Corollary \ref{hoppe5}. We also fix $ \rho, \rho_{1},\rho_{2},\rho_{3}$ given by Lemma \ref{eps}. \end{notation}

\begin{lemma} \label{boi}
There exists $t_{0}>0$ such that for every $t<t_{0}$, there exists an integer $j= j(t)$ such that: 
\begin{enumerate}
\item  $ f_{1}^{-j} \Big(    \tau + \mathbb{D}(0,t\rho_{1}) e^{u}(\tau)+  \partial \mathbb{D}(0,t\rho_{2})(0,1,0)+  \mathbb{D}(0,t\rho_{3}) (0,0,1) \Big)    \cap \mathbb{D}^{3} = \emptyset , $  
\item  $f_{1}^{-j} \Big(   \tau + \mathbb{D}(0,t\rho_{1}) e^{u}(\tau)+ \mathbb{D}(0,t\rho_{2}) (0,1,0)+ \mathbb{D}(0,t\rho_{3}) (0,0,1)  \Big) \cap \partial \mathbb{D}^{3} \Subset \mathbb{D} \times \partial \mathbb{D} \times \mathbb{D}.  $ 
\end{enumerate}
Moreover, $j(t)$ tends to $+\infty$ when $t$ tends to 0. Reducing $\mathcal{F}'$, by continuity this remains true for any $f \in \mathcal{F}'$ for a given $t<t_{0}$. 
\end{lemma}

\begin{proof}
We foliate $\mathbb{D}_{t,\rho_{1},\rho_{2},\rho_{3}} =\tau+ \mathbb{D}(0,t\rho_{1}) e^{u}(\tau)+ \mathbb{D}(0,t\rho_{2}) (0,1,0)+ \mathbb{D}(0,t\rho_{3}) (0,0,1)    $ by disks $\mathcal{L}_{z_{1},z_{3}}$ parallel to the $z_{2}$ axis. To simplify, we can suppose that $\phi_{f_{1}}$ is a fixed point of $f_{1}$, up to replacing $f_{1}$ by one of its iterates. By invariance of $C^{ss}$ under $f_{1}^{-1}$, for every $n \ge 0$, the image of this tridisk under $f^{-n}_{1}$ is foliated by the $ss$-curves $f^{-n}_{1} (\mathcal{L}_{z_{1},z_{3}}) \cap \mathbb{D}^{3}$. We call the length of a $ss$-curve the radius of the maximal (for the inclusion) disk included in it. For every $n \ge 0$, let $l_{n}$ be the minimum of the lengths over all the $s$-curves $f^{-n}_{1} (\mathcal{L}_{z_{1},z_{3}})$. For every $n \ge 0$, we denote by $d_{n}$ the maximum of the diameters of the sets $\mathcal{V} \cap f^{-n}_{1} ( \mathbb{D}_{t,\rho_{1},\rho_{2},\rho_{3}} ) $ where $\mathcal{V}$ varies in the set of $(1,3)$-quasi planes over $\mathbb{D}^{2}$. Every vector in $C^{ss}$ is dilated under $f_{1}^{-1}$ by a factor close to $\Delta = |p'_{c_{1}}(\mathrm{pr}_{1}(\phi_{f_{1}}))|/|b|$.  For every $(1,3)$-quasi plane $\mathcal{V}$ over $\mathbb{D}^{2}$, $f_{1}(\mathcal{V}) \cap \mathbb{D}^{3}$ contains a $(1,3)$-quasi plane. Every tangent vector $u$ to $f_{1}(\mathcal{V}) \cap \mathbb{D}^{3}$ is of the form $u = u^{1} \cdot (1,0,0)+u^{2} \cdot (0,1,0)+u^{3} \cdot (0,0,1)$ with $|u^{2}| \le \mathrm{max}(|u^{1}|,|u^{3}|)$. Then the vector $u$ is dilated under $f^{-1}_{1}$ by less than $\frac{2}{3} \cdot \Delta$. Then $d_{n+1} \le \frac{2}{3} \cdot \Delta \cdot d_{n}$. This implies that  if $t$ is smaller than some value $t_{0}$, there exists $j =j(t)$ such that $(l_{j}  \rho_{3}) \cdot (d_{j} \rho_{2})^{-1} \ge (3/2)^{j} \cdot \rho_{3} \cdot   \rho_{2}^{-1} \ge 10$, $l_{j} \ge 10$ and $d_{j} \le 10^{-10} \cdot |b| \cdot \mathrm{min}\big( 1, \mathrm{dist}(\mathrm{pr}_{1}(\phi_{f_{1}}),\mathbb{D}),\mathrm{dist} (\mathrm{pr}_{3}(\phi_{f_{1}}),\mathbb{D}) \big)$. Increasing a last time $j$ in order to make  $f^{-j}_{1}(\tau)$ closer to $\phi_{f_{1}}$, this implies both items    (1) and (2).  When $t$ tends to 0, the lengths of the disks $\mathcal{L}_{z_{1},z_{3}} $ tend to 0 and then $j=j(t)$ tends to $+\infty$. Reducing $\mathcal{F}'$, by continuity this remains true for $f \in \mathcal{F}'$ for a given $t$. The proof is done.
\end{proof}

Here is a technical lemma which ensures the existence of a foliation with particular properties. In the following, we will say that a graph $\mathcal{V}_{x,y}$ of class $C^{2}$ over $z_{1} \in \mathbb{D}$ is  of slope bounded by $\mathcal{S}<+ \infty$ if  every tangent vector to $\mathcal{V}_{x,y}$ is of the form a multiple of $(1,\varepsilon_{2},\varepsilon_{3})$ with $ |\varepsilon_{2} | \le \mathcal{S}$ and $ |\varepsilon_{3} | \le \mathcal{S}$. For convenience, it will be useful to work with graphs of class $C^{2}$ which are not necessarily holomorphic. This will not be a problem since we will apply later only the Inclination Lemma on them and this does not need the holomorphic assumption.

\begin{lemma} \label{feuille} 
Let $\mathcal{V}$ be a $(1,3)$-quasi plane foliated by $u$-curves $(\mathcal{V}^{x})_{x \in \mathbb{D}}$. Let $\vartheta = ( \vartheta_{1} , \vartheta_{2} , \vartheta_{3}) \in \mathbb{D}^{3}$ be a point such that $\vartheta \notin \mathcal{V}$. Let us take a vector $u \in C^{u}$. Then there exists a foliation of class $C^{2}$ of a neighborhood of $\vartheta$ by graphs $\mathcal{V}_{x,y}$ of class $C^{2}$ over $z_{1} \in \mathbb{D}$ of slope bounded by $\mathcal{S}$, where $0<\mathcal{S}<+\infty$ is independent of $\mathcal{V}$ and $\vartheta$. Moreover we have:
\begin{enumerate} 
\item  for every $x \in \mathbb{D}$, $\mathcal{V}^{x} = \mathcal{V}_{x,0}$, 
\item the leaf going through $\vartheta$ has $\mathbb{C} \cdot u$ as tangent space.
\end{enumerate}
\end{lemma}

\begin{proof}  
Up to multiplying $u$ by a non zero complex number, we can suppose that the first coordinate of $u$ is equal to 1. We first perform a change of coordinates by a biholomorphism $\varphi$, and then construct the graphs $\mathcal{V}_{x,y}$. We construct $\varphi$ such that  $\varphi$ sends $\mathcal{V}$ to the plane $\{z_{2} = 0\}$ and each curve $\mathcal{V}^{x}$ to the line $\{z_{2} = 0,z_{3} = x\}$. For a given $z \in \mathbb{D}^{3}$, we denote by $\mathrm{pr}(z)$ the projection (parallel to $(0,1,0)$) of $z$ on $\mathcal{V}$. We denote by $\varphi_{2}(z) $ the complex number such that $z -\mathrm{pr}(z) =   \varphi_{2}(z)  \cdot (0,1,0)$. Since $\mathcal{V}$ is a $(1,3)$-quasi plane, we have $\|D\varphi_{2}\| \le 1$. We denote by by $\varphi_{3}(z) $ the complex number $x$ such that $\mathrm{pr}(z) \in \mathcal{V}^{x}$. Up to rescaling $x \mapsto \mathcal{V}^{x}$, we can suppose that $\mathcal{V}^{x}$  depends on the third coordinate $x$ of the intersection point of $\mathcal{V}^{x}$ with $\{z_{1} = 0\}$. Then we also have $\|D\varphi_{3}\| \le 1$. We then define the map $\varphi$ by $\varphi(z) = (z_{1}, \varphi_{2}(z), \varphi_{3}(z))$.  In the coordinates given by $\varphi$, we have $\mathcal{V} = \{z_{2} = 0\}$ and $\mathcal{V}^{x}= \{z_{2} = 0,z_{3} = x\}$. We denote $\varphi(\vartheta) = (\tilde{ \vartheta_{1}},\tilde{ \vartheta_{2}} ,  \tilde{  \vartheta_{3}})$ (with $ \vartheta_{1}=\tilde{ \vartheta_{1}}$ and  $|\tilde{ \vartheta_{2} } |>  0 $) and $D_{\vartheta}\varphi(u) = (1,\epsilon_{2},\epsilon_{3})$. Since $u \in C^{u}$, the estimates $\|D\varphi_{2}\| \le 1$ and $\|D\varphi_{3}\| \le 1$ imply $ |\epsilon_{2} | \le 5$ and $ |\epsilon_{3} |\le 5$. \medskip

We now take a graph $W = w(\mathbb{D}^{2})$ of class $C^{2}$ over $(z_{1},z_{3}) \in \mathbb{D}^{2} \subset \mathbb{R}^{4}$ such that $\tilde{\vartheta} \in W$ with $W \cap \{z_{2} = 0\} = \emptyset$ and having $\mathbb{C} \cdot D_{\vartheta}\varphi(u) $ as tangent space at $\tilde{\vartheta}$. Since $D_{\vartheta}\varphi(u) = (1,\epsilon_{2},\epsilon_{3})$ with $ |\epsilon_{2} | \le 5$ and $ |\epsilon_{3} |\le 5$, we can take $w$ such that $\|Dw\| \le 10$. Then for every $y \in \mathbb{D}(0,10|\tilde{\vartheta_{2}}|)$, we define the $C^{2}$ graph $W_{y} = w_{y}(\mathbb{D}^{2})$ over $(z_{1},z_{3}) \in \mathbb{D}^{2} \subset \mathbb{R}^{4}$ by $w_{y} (z_{1},z_{3}) = y \cdot w(z_{1},z_{3})$.  We notice that $W_{0} =\{z_{2} = 0\}$ and $W_{1} = W$. Then for every $(x,y) \in \mathbb{D} \times \mathbb{D}(0,10|\tilde{\vartheta_{2}}|)$, the set $\mathcal{V}_{x,y}$ equal to the image of $z_{1} \in \mathbb{D} \mapsto (z_{1},w_{y}(z_{1},x),x) \in \mathbb{R}^{6}$ is a graph of class $C^{2}$ over $z_{1} \in \mathbb{D}$. Since $\|Dw\| \le 10$, its slope is bounded by $10 \times 10 = 100$ in the coordinates given by $\varphi$. Since $\mathcal{V}$ is a $(1,3)$-quasi plane in the initial coordinates, $D\varphi^{-1}$ is $C^{0}$-bounded and then the slope of any $\mathcal{V}_{x,y}$ in the initial coordinates is bounded  by $\mathcal{S}$, where $0<\mathcal{S}<+\infty$ is independent of $\mathcal{V}$ and $\vartheta$. Since $W \cap \{z_{2} = 0\} = \emptyset$, we have $w(z_{1},z_{2}) \neq 0$ for every $(z_{1},z_{3}) \in \mathbb{D}^{2}$ and then the graphs $\mathcal{V}_{x,y}$ form a foliation of class $C^{2}$ of a neighborhood of $\vartheta$ satisfying the conditions (1) and (2).
\end{proof}

We can reduce $\sigma_{1}$, $b_{1}$ and $\mathcal{F}'$ so that for every $f \in \mathcal{F}'$, $W^{s}_{\mathrm{loc}}(\phi_{f})$ has a point of intersection with every graph $\mathcal{L} \subset \mathbb{D}^{3}$ of class $C^{2}$ of slope bounded by $\mathcal{S}$ over $z_{1} \in \mathbb{D}$, and this intersection is transverse (uniformly in $\mathcal{L}$). The following is an easy consequence of the Inclination Lemma (one can refer to Prop. 6.2.23 page 257 of \cite{caga}).

\begin{lemma} \label{inclin} 
There exists an integer $k$ such that for every $ f \in \mathcal{F}'$, for every graph $\mathcal{L} \subset \mathbb{D}^{3}$ of class $C^{2}$ of slope bounded by $\mathcal{S}$ over $z_{1} \in \mathbb{D}$, $f^{k}(\mathcal{L})$ has a component which is a graph over $z_{1} \in \mathbb{D}$ which is $\eta$-close to $W^{u}_{\mathrm{loc}}(\phi_{f})$. 
\end{lemma}  

\begin{notation}  \label{nottheta} 
We fix the integer $k$ given by Lemma \ref{inclin}. Up to increasing $k$ if necessary, we can fix   $0<t<t_{0}$ such that $j(t) = k $ (see Lemma \ref{boi}). We reduce $\mathcal{F}'$ so that for every $f \in \mathcal{F}'$: there is a point $\vartheta \in W^{s}(\alpha_{f})$ such that $f^{k}(\vartheta) \in ( \tau+ (\mathbb{D}(0,t\rho))^{3})$ and $W^{s}(\alpha_{f})$ has a tangent vector $u \in C^{u}$ at $\vartheta$ (it is possible by continuity since it is the case for $f_{1}$ with the point $f_{1}^{-k}(\tau)$).  
\end{notation}

\includegraphics[width=11cm]{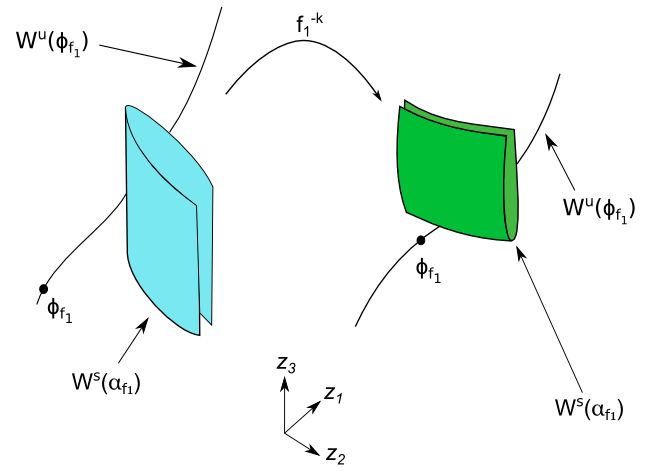}

\begin{center} 
Figure 2: straightening of the initial fold by iterating backwards. The blue and green sets are the connected components $\mathcal{W}_{f}$ and $\mathcal{W}''_{f}$ of $ W^{s}(\alpha_{f}) \cap \mathbb{D}^{3}$. In particular,  the green set $\mathcal{W}''_{f}$ is a  concentrated 3-folded $(2,3)$-surface.
\end{center}

\begin{prop} \label{fm}  
For every $f \in \mathcal{F}'$, $W^{s}(\alpha_{f})$ contains a concentrated 3-folded $(2,3)$-surface $\mathcal{W}''_{f}$ such that $\mathrm{Fold}(\mathcal{W}''_{f}) \subset  \mathbb{D}(0,\frac{1}{10})$ and $\mathrm{diam}(\mathrm{pr}_{1}(\mathcal{W}''_{f})) < 10^{-10}$. 
\end{prop}

\begin{proof}
 Let $f \in \mathcal{F}'$. Let  $\mathcal{V}$ be a $(1,3)$-quasi plane foliated by the $u$-curves $\mathcal{V}^{x} = \mathcal{V} \cap \{z_{3} = x \}$. By Notation \ref{nottheta}, we can take $\vartheta \in W^{s}(\alpha_{f})$ such that $f^{k}(\vartheta) \in (\tau+ (\mathbb{D}(0,t\rho))^{3})$, $ W^{s}(\alpha_{f})$ has a tangent vector in $C^{u}$ at $\vartheta$  and $\vartheta$ does not belong to $\mathcal{V}$. According to Lemma \ref{feuille} it is possible to find a foliation of class $C^{2}$ of a neighborhood of $\vartheta$ by graphs $\mathcal{V}_{x,y}$ of class $C^{2}$ over $z_{1} \in \mathbb{D}$ of slope bounded by $\mathcal{S}$ such that $\mathcal{V}^{x} = \mathcal{V}_{x,0}$ and the leaf going through $\vartheta$  has $\mathbb{C} \cdot u$ as tangent space. In particular, the leaf going through $\vartheta$ is tangent to $W^{s}(\alpha_{f})$. We apply $f^{k}$ to all these curves $\mathcal{V}_{x,y}$. According to Lemma \ref{inclin}, the sets $f^{k}( \mathcal{V}_{x,y} ) $ have components which foliate a neighborhood of $\tau$ by graphs over $z_{1} \in \mathbb{D}$ which are $\eta$-close to $W_{\text{loc}}^{u}(\phi_{f})$. If necessary, we can extend this foliation into a foliation of  $\tau + \mathbb{D}(0,t\rho_{1}) e^{u}(\tau)+  \mathbb{D}(0,t\rho_{2})(0,1,0)+  \mathbb{D}(0,t\rho_{3}) (0,0,1) $ by graphs over $z_{1} \in \mathbb{D}$ $\eta$-close to $W_{\text{loc}}^{u}(\phi_{f})$.    Then according to Lemma \ref{hoppe5}, the set $\text{Tan}$ of points of $f^{k}( \mathcal{V}_{x,y} ) $ where $f^{k}( \mathcal{V}_{x,y} )$ is tangent to $W^{s}(\alpha_{f})$ is a regular two-dimensional real manifold which has  its direction $\epsilon$-close to $\mathbb{C} \cdot v_{\text{tan}}$ at each point and goes through the point $f^{k}(\vartheta) \in \tau+ (\mathbb{D}(0,t\rho))^{3}$. According to Lemma \ref{eps}, the set  $\mathrm{Tan}$ can intersect $\tau + \mathbb{D}(0,t\rho_{1}) e^{u}(\tau)+  \mathbb{D}(0,t\rho_{2})(0,1,0)+  \mathbb{D}(0,t\rho_{3}) (0,0,1) $ only in $\tau + \mathbb{D}(0,t\rho_{1}) e^{u}(\tau)+  \partial \mathbb{D}(0,t\rho_{2})(0,1,0)+  \mathbb{D}(0,t\rho_{3}) (0,0,1) $. By Lemma \ref{boi}, it holds: $$f^{-k} \Big(   \tau + \mathbb{D}(0,t\rho_{1}) e^{u}(\tau)+\partial \mathbb{D}(0,t\rho_{2})(0,1,0)+  \mathbb{D}(0,t\rho_{3}) (0,0,1) \Big)    \subset \mathbb{C}^{3} \backslash\mathbb{D}^{3} \, ,$$  $$f^{-k} \Big(  \tau + \mathbb{D}(0,t\rho_{1}) e^{u}(\tau)+ \mathbb{D}(0,t\rho_{2}) (0,1,0)+ \mathbb{D}(0,t\rho_{3}) (0,0,1)  \Big) \cap \partial \mathbb{D}^{3} \Subset \mathbb{D} \times \partial \mathbb{D} \times \mathbb{D} \, .$$ Then $f^{-k}(   \text{Tan} )$ intersects $\mathcal{V}$. Since both $\mathcal{V}$ and $W^{s}(\alpha_{f})$ are complex manifolds, $f^{-k}(   \text{Tan} )$ intersects $\mathcal{V}$ in exactly one point. We denote by $\mathcal{W}'_{f}$ the continuation of the connected component of $W^{s}(\alpha_{f_{1}}) \cap \mathbb{D}^{3}$ containing $f_{1}^{-k}(\tau)$ for $f \in \mathcal{F}'$. Then $\mathcal{W}'_{f}$ is tangent to exactly one $\mathcal{V}^{x} = \mathcal{V} \cap \{z_{3} = x \}$. This implies that $\text{pr}_{3}$ restricted on $\mathcal{W}'_{f} \cap \mathcal{V}$ is a two-covering with exactly one point of ramification. Since this is true for any $(1,3)$-quasi plane $\mathcal{V}$, by definition, $\mathcal{W}'_{f}$ is a 3-folded $(2,3)$-surface for every $f \in \mathcal{F}'$. We notice that $ \mathrm{pr}_{3} (\phi_{f}) \in  \mathbb{D}(0,\frac{1}{10}-10^{-4})$. Increasing $k$ if necessary, we have $ \mathrm{pr}_{3} (f^{-k}(\tau)) \in  \mathbb{D}(0,\frac{1}{10}-10^{-4})$. Moreover $\vartheta$ can be taken as close to $f^{-k}(\tau)$ as wanted and by the proof of Lemma \ref{boi} the diameter of $\mathrm{pr}_{3} \big( f^{-k}(\mathrm{Tan}) \big) $ is smaller than $10^{-10}$. Then we have  $\mathrm{Fold}(\mathcal{W}''_{f}) \subset  \mathbb{D}(0,\frac{1}{10})$. Since the first direction is dilated under $f$, iterating by $f^{-1}$ if necessary, we can have the additional property that $\mathrm{diam}(\text{pr}_{1}(\mathcal{W}''_{f})) < 10^{-10}$. Then according to Proposition \ref{rp7}, we can iterate a last time $f^{-1}$ if necessary to get a component $\mathcal{W}''_{f}$ of $W^{s}(\alpha_{f})$ which is a concentrated 3-folded $(2,3)$-surface such that $\mathrm{Fold}(\mathcal{W}''_{f}) \subset  \mathbb{D}(0,\frac{1}{10})$ and $\mathrm{diam}(\mathrm{pr}_{1}(\mathcal{W}''_{f})) < 10^{-10}$.  
\end{proof}

\section{Proof of the main results}

\begin{proof}[Proof of the main Theorem and Corollaries 1 and 2] 
By Proposition \ref{fm}, for every $f \in \mathcal{F}'$, $W^{s}(\alpha_{f})$ contains some concentrated 3-folded $(2,3)$-surface $\mathcal{W}^{0}$ with $\mathrm{Fold}(\mathcal{W}^{0}) \subset  \mathbb{D}(0,\frac{1}{10})$ and $\text{diam}(\text{pr}_{1}(\mathcal{W}^{0})) < 10^{-10}$. This implies persistent heteroclinic tangencies. Indeed, by Proposition \ref{rrrrr} there exists a point $\kappa_{f}$ of the horseshoe $\mathcal{H}_{f}$ and a point of quadratic tangency $\tau'$ between $W^{u}(\kappa_{f})$ and  $\mathcal{W}^{0} \subset W^{s}(\alpha_{f})$. The proof of the main Theorem is complete.    \medskip
 
We now prove Corollary 1. We call $\mathcal{U}'_{f}$ the component of $W^{u}(\kappa_{f}) \cap \mathbb{D}^{3}$ containing the point of tangency $\tau'$ and we consider the family $(\mathcal{U}'_{f'})_{f' \in  \mathcal{F}'}$ of holomorphic $u$-curves given by the continuation of $\mathcal{U}'_{f}$. According to Proposition \ref{fd}, the tangency $\tau'$ is generically unfolded.  We apply a standard argument to obtain homoclinic tangencies. Let us take any neighborhood $\mathcal{F}''$ of $f$. Since the tangency $\tau'$ is generically unfolded, we can rely on the following well-know lemma:

\begin{lemma}
There exists $\epsilon''>0$ such that for every holomorphic family $(\mathcal{U}''_{f'})_{f' \in  \mathcal{F}'}$ of $u$-curves $\epsilon''$-close to $(\mathcal{U}'_{f'})_{f' \in  \mathcal{F}'}$ and every holomorphic family $(\mathcal{W}''_{f'})_{f' \in  \mathcal{F}'}$ of complex surfaces $\epsilon''$-close to $(\mathcal{W}^{0})_{f' \in  \mathcal{F}' }$, there exists a point of quadratic tangency between $\mathcal{U}''_{f''}$ and $\mathcal{W}''_{f''}$ where $f'' \in \mathcal{F}''$.
\end{lemma}

We claim that there exists a holomorphic family of $u$-curves $(\mathcal{U}''_{f'})_{f' \in  \mathcal{F}'}$ $\epsilon''$-close to $(\mathcal{U}'_{f'})_{f' \in  \mathcal{F}'}$ such that for every $f' \in  \mathcal{F}'$, $\mathcal{U}''_{f'}$ is a component of $W^{u}(\delta_{f'}) \cap \mathbb{D}^{3}$ and a holomorphic family of complex surfaces $(\mathcal{W}''_{f'})_{f' \in  \mathcal{F}'}$ $\epsilon''$-close to $((\mathcal{W}^{0})_{f'})_{f' \in  \mathcal{F}'}$ such that for every $f' \in  \mathcal{F}'$, $\mathcal{W}''_{f'}$ is a component of $W^{s}(\delta_{f'})$ (we recall that the point $\delta_{f'}$ was defined in Proposition \ref{ferm}).  We prove this fact for the case of $u$-curves, the proof is the same for complex surfaces. According to Proposition \ref{ferm}, $\alpha_{f'}$, and then $\kappa_{f'}$, is in the homoclinic class of $\delta_{f'}$. This implies that there exists a transversal intersection between $W^{s}_{\mathrm{loc}}(\kappa_{f'})$ and $W_{\mathrm{loc}}^{u}(\delta_{f'})$ for every $f' \in \mathcal{F}'$. For a given $f' \in \mathcal{F}'$, by the inclination lemma, there exists an integer $n_{f'}$ such that for every $n \ge n_{f'}$, $(f')^{n_{f'}} (W_{\mathrm{loc}}^{u}(\delta_{f'}))$ contains a $u$-curve $\epsilon''$-close to $\mathcal{U}'_{f'}$. By continuity, $n_{f'}$ can be taken locally constant. Up to reducing $\mathcal{F}'$, we can take  $\mathcal{F}'$ compact (with non empty interior). By compactness, it is then possible to take the maximal value $n$ of $n_{f'}$ on a finite open covering of $\mathcal{F}' $. Then for every $f' \in \mathcal{F}'$, there is a $u$-curve $\mathcal{U}''_{f'}$ $\epsilon''$-close to $\mathcal{U}'_{f'}$ which is a connected component of $W^{u}(\delta_{f'}) \cap \mathbb{D}^{3}$. The proof of Corollary 1 is complete. \medskip

We now prove Corollary 2. The preceding proof shows that there exists $f'' \in \mathcal{F}''$ with a point of homoclinic tangency between $W^{s}(\delta_{f''})$ and $W^{u}(\delta_{f''})$. In particular maps with homoclinic tangencies associated to $\delta_{f}$ are dense in $ \mathcal{F}'$. The point $\delta_{f}$ has the property of being sectionally dissipative for every $f \in  \mathcal{F}'$. Then, by Proposition \ref{finenfin} which gives the creation of sinks from homoclinic tangencies (the proof is given in Appendix for the convenience of the reader) a classical Baire category argument already used in \cite{bb1} allows us to conclude the existence of a residual set of $\text{Aut}_{2}(\mathbb{C}^{3})$ of automorphisms displaying infinitely many sinks. The proof is complete.
\end{proof}

\appendix

\section{From homoclinic tangencies to sinks}

To show the main Theorem, we need the following result. It is known since the work of Newhouse how to get a sink from a homoclinic tangency. The adaptation to the case of $\mathbb{C}^{2}$ was obtained by Gavosto in \cite{est}. Here we adapt her proof to the case of $\mathbb{C}^{3}$. Remind that a generically unfolded tangency is a tangency which is unfolded with a positive speed.

\begin{propo}  \label{finenfin}
Let $(F_{t})_{t \in \mathbb{D}}$ be a family of polynomial automorphisms of $\mathbb{C}^{3}$. We suppose that for the parameter $t = 0$, there is a sectionally dissipative periodic point $P_{0}$ and a generically unfolded homoclinic quadratic tangency $Q \in W^{s}(P_{0}) \cap W^{u}(P_{0})$ between $W^{s}(P_{0})$ and $W^{u}(P_{0})$. We also suppose that the three eigenvalues $\lambda_{0}^{ss}$, $\lambda_{0}^{cs}$ and $\lambda_{0}^{u} $ at $P_{0}$ satisfy $|\lambda_{0}^{ss}| < |\lambda_{0}^{cs}|<1< |\lambda_{0}^{u}|$ (and $|\lambda_{0}^{cs}| \cdot |\lambda_{0}^{u}|<1$ by sectional dissipativity). Then for every neighborhood $\mathcal{Q}$ of $Q$ and every neighborhood $\mathcal{T}$ of 0, there exists $t \in \mathcal{T}$ such that $F_{t}$ admits an attracting periodic point in $\mathcal{Q}$. 
\end{propo}

\begin{proof} Step 1 : Construction of cone fields \medskip 

In the following, iterating if necessary, we will suppose that $P_{0}$ is a fixed point of $F_{0}$. We fix a constant $\eta>0$ such that we have $(|\lambda^{cs}_{0}|+\eta)(|\lambda^{u}_{0}|+\eta)<1$ (the periodic point $P_{0}$ is sectionally dissipative) and $|\lambda_{0}^{ss}|+\eta   < |\lambda_{0}^{cs}| - \eta<|\lambda_{0}^{cs}| + \eta<1$. We can fix a neighborhood $\mathcal{P}$ of $P_{0}$, a neighborhood $\mathcal{T}_{0}$ of 0, cones $C^{u}$, $C^{ss}$ and $C^{cs}$ centered at the three eigenvectors of $DF_{0} (P_{0})$ and an integer $N_{0}$ with the following properties: for every $n \ge N_{0}$, $t \in \mathcal{T}_{0}$, for every matrix $M = M_{n} \cdot \ldots \cdot M_{1}$ where $M_{i}$ is the differential of $F_{t}$ at a point in $\mathcal{P}$, we have: 
\begin{enumerate}
\item $C^{u}$ is invariant under $M$ and there is exactly one eigenvector (up to multiplication) of $M$ in $C^{u}$ of eigenvalue $(|\lambda^{u}_{0}|-\eta)^{n}< |\lambda_{1}|< (|\lambda^{u}_{0}|+\eta)^{n}$, 
\item  $C^{ss}$ is invariant under $M^{-1}$ and there is exactly one eigenvector (up to multiplication) of $M$ in $C^{ss}$ of eigenvalue $(|\lambda^{ss}_{0}|-\eta)^{n}< |\lambda_{2}|< (|\lambda^{ss}_{0}|+\eta)^{n}$, \item for every  vector $v$ which is not in $C^{u} \cup C^{ss} \cup C^{cs}$, we have $M \cdot v \in C^{u}$ or $M^{-1} \cdot v \in C^{ss}$,
\item $   (|\lambda^{u}_{0}|-\eta)^{n}(|\lambda^{cs}_{0}|-\eta)^{n}(|\lambda^{ss}_{0}|-\eta)^{n}<| \det(M)  |<     (|\lambda^{u}_{0}|+\eta)^{n}(|\lambda^{cs}_{0}|+\eta)^{n}(|\lambda^{ss}_{0}|+\eta)^{n}$.
\end{enumerate}
Reducing $\eta$ and increasing $N_{0}$ if necessary, we have for $n \ge N_{0}$: 
 $$   (|\lambda^{u}_{0}|+\eta)^{n}     (|\lambda^{ss}_{0}|+\eta)^{2n}<    (|\lambda^{u}_{0}|-\eta)^{n} (|\lambda^{cs}_{0}|-\eta)^{n} (|\lambda^{ss}_{0}|-\eta)^{n},$$ $$     (|\lambda^{u}_{0}|+\eta)^{n} (|\lambda^{cs}_{0}|+\eta)^{n} (|\lambda^{ss}_{0}|+\eta)^{n}<    (|\lambda^{u}_{0}|-\eta)^{2n}(|\lambda^{ss}_{0}|-\eta)^{n}.$$ 
This implies that $M$ is diagonalizable. Indeed, if it was not the case, $M$ would be triangularizable with a double eigenvalue and we would have:  $$   (|\lambda^{u}_{0}|+\eta)^{n}     (|\lambda^{ss}_{0}|+\eta)^{2n}>| \det(M)  | \text{ or  }  | \det(M)  |>     (|\lambda^{u}_{0}|-\eta)^{2n}(|\lambda^{ss}_{0}|-\eta)^{n}.$$ But by item 4 above, there would be a contradiction with the previous inequalities. Moreover the third eigenvector has to be in $C^{cs}$ according to item 3 above. \newline \newline Step 2 : Local coordinates \medskip

Iterating if necessary, we can suppose that $Q \in \mathcal{P}$. We are going to make several local changes of coordinates so that the stable and unstable manifolds of $P_{t}$ have a simple form near the tangency. We will use this local change of coordinates in Steps 2 and 3. In Step 4, we will mainly use the coordinates in the canonical basis but we will need the local change of coordinates a last time at some point so we will denote it by $\Psi_{\tau}$ in Step 4 to make the distinction between the two systems of coordinates. In the following, to construct the local coordinates, we will keep the notation $z_{1},z_{2},z_{3}$ by simplification.\medskip

We first pick local coordinates such that $Q = 0$ and in the neighborhood of $Q$, the stable manifold is $\{z_{1} = 0\}$. Up to a linear invertible second change of coordinates, we can assure that the tangent vector of $W^{u}(P_{0})$ at $Q$ is $(0,0,1)$. Then we have that locally near 0 the unstable manifold is given by a graph over $z_{3}$ of the form  $\{(w_{1}(z_{3},t),w_{2}(z_{3},t),z_{3}) : z_{3}\}$. Since the tangency is quadratic, for each $t$ in a neighborhood of 0, there exists exactly one $z_{3,t}$ such that $ \frac{\partial w_{1}}{\partial z_{3}}(z_{3,t},t) = 0$. We pick new coordinates a third time by changing the coordinate $z_{3}$ into $z_{3}-z_{3,t}$. The unstable manifold is: $$w_{1}(z_{3}) = w_{1}(z_{3},t) = w_{1}(0,t)+\frac{\partial^{2}w_{1}}{\partial z_{3}^{2}}(0,t)z_{3}^{2}+ h(z_{3},t) \text{ with } h(z_{3},t) = o(z_{3}^{2}) .$$ The stable manifold is still locally equal to $\{z_{1} = 0\}$. Since there is a quadratic tangency for $t = 0$ which is generically unfolded, we have $w_{1}(0,0) = 0$ and $\frac{\partial w_{1}}{\partial t} (0,0) \neq 0 $. Then we change coordinates a fourth time: $z_{1}$ becomes $ \frac{t}{w_{1}(0,t)}z_{1}$ and $z_{2}$ and $z_{3}$ are unchanged. In these new coordinates, the unstable manifold is given by: $$w_{1}(z_{3}) = t +z_{3}^{2}\tilde{h}(z_{3}, t ),$$  where $\tilde{h}(z_{3}, t ) \neq 0$ in a neighborhood of 0. The stable manifold is still locally equal to $\{z_{1} = 0\}$. The fifth change of coordinates is given by $z_{3}$ becoming $z_{3}(\tilde{h}(z_{3},t))^{1/2}$, where $(\tilde{h}(z_{3},t))^{1/2}$ is a complex square root of $\tilde{h}(z_{3},t)$ (which is well defined since $\tilde{h}(z_{3}, t ) \neq 0$ in a neighborhood of 0). We finally get that the unstable manifold is given by $z_{1} = z_{3}^{2}+ t$ and $z_{2} = w_{2}(z_{3},t)$. The stable manifold is still locally equal to $\{z_{1} = 0\}$ and the tangent vector of $W^{u}(P_{0})$ at $Q$ is still $(0,0,1)$. The last change of coordinates is given by $z_{2}-w_{2}(z_{3},t)$. The unstable manifold is given by $z_{1} = z_{3}^{2}+ t$ and $z_{2} = 0$ and the stable manifold by $z_{1} = 0$. The tangent vector of $W^{u}(P_{0})$ at $Q$ is still $(0,0,1)$. \newline \newline Step 3 : Construction of a periodic point \medskip

We take a tridisk $B$ around $Q$ in the coordinates that we just defined : $B = \{(z_{1},z_{2},z_{3}) : |z_{1}| < \delta, |z_{2}|<\delta, |z_{3}|<\delta\}$ where $0<\delta<1$. Since the tangent vector of $W^{u}(P_{0})$ at $Q$ is $(0,0,1)$, reducing $\delta$ if necessary, there exists a neighborhood $\mathcal{T}_{1} \subset \mathcal{T}_{0}$ such that the component of $W^{u}(P_{t}) \cap B$ containing $Q$ (for $t = 0$) or its continuation (for $t \neq 0$) is horizontal in $B$ relatively to the third projection and is included in $\{(z_{1},z_{2},z_{3}) \in B : |z_{2}| < \frac{1}{10}\delta \}$. Using the Inclination Lemma, there exists $N_{1} \ge N_{0}$ such that for $t \in \mathcal{T}_{1}$ and $n \ge N_{1}$ , $F_{t}^{n}(B)$ will intersect $B$ and $F_{t}(B) \cap B$ is horizontal relatively to the third projection. \medskip 

In the following, we show that for every sufficiently high $n$, a periodic point for $F_{t}$ (where $t \in \mathcal{T}_{1}$) is created. Let us now denote $\Delta_{z_{2},z_{3}} = \{(z_{1},z_{2},z_{3}) : |z_{1}|<\delta  \}$ which is a disk, for $|z_{2}|<\delta,|z_{3}|<\delta$. It is possible to increase $N_{2} \ge N_{1}$ such that for any $|z_{2}|<\delta,|z_{3}|<\delta$, for every $n \ge N_{2}$, $F_{t}^{n}(\Delta_{z_{2},z_{3}} ) \cap B$ is horizontal relatively to the third projection (of degree 1) and included in $\{(z_{1},z_{2},z_{3}) \in B : |z_{2}| < \frac{1}{2}\delta \}$. Since $F_{t}^{n}(\Delta_{z_{2},z_{3}} ) \cap B$ is horizontal relatively to the third projection and of degree 1, it intersects exactly one disk $\Delta_{z'_{2},z_{3}}$ where $z'_{2} = z'_{2}(z_{2})$ with $z'_{2} \in \mathbb{D}(0,\delta/2)$. This defines for each fixed $t \in \mathcal{T}_{1}$ a holomorphic map $z_{2} \mapsto z'_{2}$. Then the map $z_{2} \mapsto z_{2}-z'_{2}$ defined on $\mathbb{D}(0,\delta)$ is holomorphic and the image of $\partial \mathbb{D}(0,\delta)$ contains a loop around $\mathbb{D}(0,\delta/2)$. Then by the Argument principle there exists $z_{2}(z_{3})$ such that $\Delta_{z_{2}(z_{3}),z_{3}}$ intersects $F^{n}_{t}(\Delta_{z_{2}(z_{3}),z_{3}})$. \medskip 

We are going to choose $z_{3}$ in order to create a periodic point in $\Delta_{z_{2}(z_{3}),z_{3}}$. There is a point $R_{z_{3}} = (f_{t}(z_{3}),z_{2}(z_{3}),z_{3}) \in \Delta_{z_{2}(z_{3}),z_{3}}$ which is sent on $S_{z_{3}} \in \Delta_{z_{2}(z_{3}),z_{3}}$ where $S_{z_{3}} = (g_{t}(z_{3}),z_{2}(z_{3}),z_{3} )$. We are going to choose $z_{3}$ (in function of $t$) such that $R_{z_{3}} =S_{z_{3}}$.  When $N$ goes to infinity, $f_{t}(z_{3})$ tends to 0 and $g_{t}(z_{3})$ tends to $z_{3}^{2}+t$. Then $g_{t}(z_{3}) - f_{t}(z_{3})$ tends to $z_{3}^{2}+t$. In particular, if $N$ is sufficiently high, for every $t \in \mathbb{D}$, the graph of $z_{3} \mapsto g_{t}(z_{3}) - f_{t}(z_{3})$ is a curve of degree 2 over $z_{3}$ which has exactly one point of tangency with the horizontal foliation.  We take a new bound $N_{3}$ on $n $ such that this is the case, we increase $N_{3}$ in the rest of Step 3 in order to satisfy more assumptions.  The second coordinate of this point of tangency is a holomorphic function of $t$ which tends to $t$ when increasing $N_{3}$. In particular, this implies that there exists exactly one value $t_{0}$ of $t$ for which the equation $g_{t_{0}}(z_{3}) - f_{t_{0}}(z_{3}) = 0$ has one double solution, and for every other value of $t \in \mathbb{D}$, there are two distinct solutions. Then, for $F_{t}$, we have two periodic points which are equal when $t = t_{0}$. We do a reparametrization of the family of maps $(F_{t})_{t \in \mathbb{D}}$ by  taking $t = \tau^{2}$. From now on, we are working with the family of maps $\tilde{F}_{\tau} = F_{\tau^{2}}$. By simplicity, we will simply denote it by $F_{\tau}$. For $F_{\tau}$, we have two periodic points which are equal when $\tau = \tau_{0}$ where $\tau_{0}^{2} = t_{0}$. We can increase $N_{3}$ if necessary so that for each $\tau \in \mathbb{D}  \backslash \mathbb{D}(0,\frac{1}{2})$, the two solutions of $g_{\tau^{2}}(z_{3}) - f_{\tau^{2}}(z_{3}) = 0$ are respectively $\frac{1}{100}$-close to $\pm i \tau$. For $\tau \in \mathbb{D} \backslash \mathbb{D}(0,\frac{1}{2})$, we denote by $R^{\tau}$ the periodic point corresponding to the solution $\frac{1}{100}$-close to $i\tau$. It is clear that the map $\tau \mapsto R^{\tau}$ restricted on $\mathbb{D} \backslash \mathbb{D}(0,\frac{1}{2})$ is continuous. For $\tau = \tau_{0}$, we denote by $R^{\tau_{0}}$ the unique periodic point corresponding to the double solution of $g_{\tau_{0}^{2}}(z_{3}) - f_{\tau_{0}^{2}}(z_{3}) = 0 $. Finally, for $\tau \in \mathbb{D}(0,\frac{1}{2}) $ not equal to $\tau_{0}$, there are two distinct periodic points in $\Delta_{z_{2}(z_{3}),z_{3}}$. We pick any path in $\mathbb{C}$ from $\tau$ to a point $\tau_{1}$ in $\mathbb{D} \backslash    \mathbb{D}(0,\frac{1}{2})$ which does not contain $\tau_{0}$. For $\tau_{1}$, the periodic point $R^{\tau_{1}}$ is defined. Since the path does not contain $\tau_{0}$, there is exactly one of the two distinct periodic points in $\Delta_{z_{2}(z_{3}),z_{3}}$ for $\tau$ which is the continuation of $R^{\tau_{1}}$. We denote it by $R^{\tau}$. Since the map $\tau \mapsto R^{\tau}$ restricted on $\mathbb{D}     \backslash \mathbb{D}(0,\frac{1}{2})$ is continuous, this choice is independent of a particular choice of $\tau_{1}$ in $\mathbb{D} \backslash    \mathbb{D}(0,\frac{1}{2})$. Then we have defined a map $\tau \mapsto R^{\tau}$ on $\mathbb{D}$. It is clear that near any point of $\mathbb{D} \backslash \{\tau_{0}\}$, $R^{\tau}$ is locally the continuation of the same periodic point of $F_{\tau}$, then it is holomorphic. The map $\tau \mapsto R^{\tau}$ is holomorphic on $\mathbb{D} \backslash \{\tau_{0}\}$. It is trivial that it is continuous at $\tau_{0}$. Then it is holomorphic on $\mathbb{D}$. As we already said, from now on, we go back to the coordinates in the canonical basis and we call $\Psi_{\tau}$ the local coordinates we just used near the tangency point $Q$. We will use $\Psi_{\tau}$ a last time at the end of Step 4. \newline \newline  Step 4 : $R^{\tau}$ is a sink \medskip 
 
We now show it is possible to pick $\tau$ such that $R^{\tau}$ is a sink. From now on, we fix a neighborhood $\mathcal{T} \subset \mathcal{T}_{1}$ of 0. Recall that $Q$ belongs to the small neighborhood $\mathcal{P}$ of $P_{0}$ defined in Step 1. We denote by $n = n_{1}+n_{2}$ such that for $k = 1,...,n_{1}$, $F^{k}_{\tau}(R^{\tau})$ is in $\mathcal{P}$ and $F^{n_{1}+1}_{\tau}(R^{\tau}) \notin \mathcal{P}$. We express all matrices in the $\tau$-dependent basis given by the 3 eigenvectors $e'_{1}$,$e'_{2},e'_{3}$ of $D(F^{n_{1}}_{\tau})(R^{\tau})$ (the matrix $DF^{n_{1}}_{\tau}(R^{\tau})$ is diagonalizable according to Step 1). We have that $e'_{1} \in C^{u}$, $e'_{2} \in C^{ss}$ and $e'_{3} \in C^{cs}$. Then, in this basis, the matrix $DF_{\tau}^{n}(R^{\tau})$ is of the form $ DF^{n_{1}+n_{2}}_{\tau}(R^{\tau}) = DF^{n_{2}}_{\tau}(F^{n_{1}}_{\tau}(R^{\tau})) \cdot DF^{n_{1}}_{\tau}(R^{\tau})$. The matrix $DF^{n_{1}}_{\tau}(R^{\tau})$ is of the form:
$$\begin{pmatrix}
   \Lambda_{1}&0&0\\
0&\Lambda_{2}&0 \\
  0& 0&\Lambda_{3}\\
\end{pmatrix} \, .$$ 
We have that $(|\lambda^{u}_{0}| - \eta)^{n_{1}} <   |\Lambda_{1}| < (|\lambda^{u}_{0}| + \eta)^{n_{1}}$, $(|\lambda^{ss}_{0}| - \eta)^{n_{1}} <   |\Lambda_{2}| < (|\lambda^{ss}_{0}| + \eta)^{n_{1}}$ and $(|\lambda^{cs}_{0}| - \eta)^{n_{1}} <   |\Lambda_{3}| < (|\lambda^{cs}_{0}| + \eta)^{n_{1}}$. The matrix $DF^{n_{2}}_{\tau}(F^{n_{1}}_{\tau}(R^{\tau}))$ is of the form:  
$$ \begin{pmatrix}
   A&B&C\\
   D& E&F \\
   G&H& I \\
\end{pmatrix}\, .$$ 
Since $e'_{1} \in C^{u}$, $e'_{2} \in C^{ss}$ and $e'_{3} \in C^{cs}$ and these cones are disjoint with $n_{2}$ bounded, the coefficients $A$,$B$,$C$,$D$,$E$,$F$,$G$,$H$,$I$ are bounded in modulus by some constant $K$ which is independant of $\tau$, $n$ and $n_{1}$. Then we have:
 $$DF^{n_{1}+n_{2}}_{\tau}(R^{\tau}) = \begin{pmatrix}
 A\Lambda_{1} &B\Lambda_{2}&C\Lambda_{3} \\
  D\Lambda_{1} &E\Lambda_{2}&F\Lambda_{3}\\
G\Lambda_{1} &H \Lambda_{2}&I \Lambda_{3}\\  \end{pmatrix} \, .$$
Let us suppose that $A\Lambda_{1} = 0$. There exists $\epsilon>0$ such that if the characteristic polynomial of $DF^{n_{1}+n_{2}}_{\tau}(R^{\tau})$ is $X^{3}+a_{2}X^{2}+a_{1}X +a_{0}$ with $|a_{0}|<\epsilon, |a_{1}|<\epsilon, |a_{2}|<\epsilon$, then the 3 eigenvalues of $DF^{n_{1}+n_{2}}_{\tau}(R^{\tau})$ are of modulus lower than 1 and then $R^{\tau}$ is a sink. We are going to show that it is possible to get a lower bound on  $n_{1}$ so that it is always the case. The coefficients of the characteristic polynomial of $DF^{n_{1}+n_{2}}_{\tau}(R^{\tau})$ are : 
$$a_{2} = - ( A\Lambda_{1}+E\Lambda_{2} + I\Lambda_{3}) , $$
$$ a_{1} =  (EI\Lambda_{2}\Lambda_{3}-FH\Lambda_{2}\Lambda_{3} +AI\Lambda_{1}\Lambda_{3}-CG\Lambda_{1}\Lambda_{3}+AE\Lambda_{1}\Lambda_{2}-BD\Lambda_{1}\Lambda_{2}) , $$
$$ a_{0} = -\det( DF^{n_{2}}_{\tau}(F^{n_{1}}_{\tau}(R^{\tau} )) .$$
We have:  $ |a_{0}| =  | \det( DF^{n_{2}}_{\tau}(F^{n_{1}}_{\tau}(R^{\tau}) ) )|< 6  K^{3}|\Lambda_{1}\Lambda_{2}\Lambda_{3}|<  6 K^{3} \big( (|\lambda^{cs}_{0}|+\eta)(|\lambda^{u}_{0}|+\eta) \big)^{n_{1}}$ which tends to 0 when $n_{1} \rightarrow +\infty $ because  $(|\lambda^{cs}_{0}|+\eta)(|\lambda^{u}_{0}|+\eta)<1$ (remind that $P_{0}$ is a sectionally dissipative periodic point). We increase $n_{1}$ such that $| a_{0} |<\epsilon$. We have that: $$| a_{1}  |  < 6K^{2} \max ( |\Lambda_{1}\Lambda_{2}|,   |\Lambda_{1}\Lambda_{3}|, |\Lambda_{2}\Lambda_{3}|) < 6K^{2}  \big( (|\lambda^{cs}_{0}|+\eta)(|\lambda^{u}_{0}|+\eta) \big)^{n_{1}},$$ which tends to 0 when $n_{1} \rightarrow + \infty $ because $(|\lambda^{cs}_{0}|+\eta)(|\lambda^{u}_{0}|+\eta)<1$. We increase $n_{1}$ further such that $| a_{1} |<\epsilon$. Finally, in the term $a_{2}$, both $E\Lambda_{2}$ and $I\Lambda_{3}$ tend to 0 when $n_{1} \rightarrow + \infty $. The term $A\Lambda_{1}$ is equal to 0 by hypothesis. We increase $n_{1}$ a last time such that $| a_{2} |<\epsilon$. Finally, we pick $N_{4} \ge N_{3}$ such that $n_{1}$ is sufficiently high in order to satisfy the previous inequalities. Finally, the three  eigenvalues of $DF^{n_{1}+n_{2}}_{\tau}(R^{\tau})$ are of modulus lower than 1 and $R^{\tau}$ is a sink. \medskip

It remains us to show that for a given neighborhood $t \in \mathcal{T}$ of 0 and the corresponding neighborhood of 0 for $\tau$ (remind that $t = \tau^{2}$), it is possible to pick a new bound $N_{6}$ on $n$, such that for $n \ge N_{6}$, there is a parameter $\tau=\tau(n)$ such that the differential $DF_{\tau}^{n}(R^{\tau})$ satisfies $A\Lambda_{1} = 0$. To show this, we work again in the coordinates $\Psi_{\tau}$ for the map $F_{\tau}$. We denote by $\Pi$ the projection of the plane $\Psi_{\tau}(W^{s}(P_{\tau^{2}})) = \{z_{1} = 0\}$ into $\mathbb{P}^{2}(\mathbb{C})$. The projection $\Pi$ is a holomorphic curve. We denote by $\Gamma$ the holomorphic curve in $\mathbb{P}^{2}(\mathbb{C})$ given by the tangent directions to the curve $\Psi_{\tau}(W^{u}(P_{0})) = \{z_{1} =  z_{3}^{2}, z_{2} = 0\}$. Since there is a generically unfolded quadratic tangency at $\tau = 0$ between $W^{s}(P_{0})$ and $W^{u}(P_{0})$, there is a transverse intersection between $\Pi$ and $\Gamma$. Moreover, $\mathbb{P}^{2}(\mathbb{C})$ is of dimension 2. We take a small disk $\Delta'$ going through $\Psi_{\tau}(R^{\tau})$ of direction $D\Psi_{\tau}( R^{\tau} ) \cdot e'_{1}$. For a given $n \ge N_{4}$, we call $\Gamma'_{n}$ the complex curve in $\mathbb{P}^{2}(\mathbb{C})$ given by the tangent directions to the curve $\Psi_{\tau}(F^{n}_{\tau}(\Delta'))$ at the periodic point $\Psi_{\tau}(R^{\tau})$ when $\tau$ varies. By the Inclination Lemma, we can increase the bound $N_{5} \ge N_{4}$ on $n$ such that if $n \ge N_{5}$, $\text{pr}_{3}(\Psi_{\tau}(R^{\tau}))$ can be made as close to $i\tau$ as wanted and $\Psi_{\tau}(F^{n}_{\tau}(\Delta'))$ can be taken as close to $\{z_{1} =  z_{3}^{2}+t, z_{2} = 0\}$ in the $C^{1}$ topology. In particular, this shows that the graph $\Gamma'_{n}$ can be taken as close to the graph $\Gamma$ as wanted (in the $C^{1}$-topology) by increasing the bound $N_{5}$. In particular, we can pick $N_{5}$ such that there is a transverse intersection between $\Gamma'_{n}$ and $\Pi$ if $n \ge N_{5}$. There is a last step to get $A\Lambda_{1} = 0$. We denote by $\Pi'$ the projection of the plane $\mathrm{Vect}(D\Psi_{\tau}(R^{\tau}) \cdot e'_{2}, D\Psi_{\tau}(R^{\tau})  \cdot e'_{3})$ in $\mathbb{P}^{2}(\mathbb{C})$, $\Pi'$ is a complex curve. Since $e'_{2}$ and $e'_{3}$ are the two stable eigenvectors of $DF^{n_{1}}_{\tau} (R^{\tau}) $, it is an easy consequence of the Inclination Lemma that if $n$ and then $n_{1}$ tend to infinity, then $\Pi'$ tends to $\Pi$ (locally as graphs in the $C^{1}$-topology). In particular, it is possible to pick a last bound $N_{6} \ge N_{5}$ on $n$ such that if $n \ge N_{6}$, then $\Gamma'_{n}$ and $\Pi'$ have a transverse intersection. Then  there exists $\tau$ with $t = \tau^{2} \in \mathcal{T}$ such that $DF^{n}_{\tau} (R^{\tau}) \cdot e'_{1} \in \mathrm{Vect}(e'_{2},e'_{3})$ which is equivalent to $A\Lambda_{1} = 0$. For this parameter $\tau$, we already saw that $R^{\tau}$ is a sink. \medskip 

This is true for every neighborhood $\mathcal{T}$ of 0 and it is clear that $R^{\tau}$ tends to $Q$ when reducing $\mathcal{T}$. The result is proven: for every neighborhood $\mathcal{Q}$ of $Q$ and every neighborhood $\mathcal{T}$ of 0, we can create a sink in $\mathcal{Q}$ for $F_{t}$ with $t \in \mathcal{T}$. \end{proof}

\bibliographystyle{plain}
\bibliography{/bibli}

\end{document}